\documentclass[11pt]{article}
\usepackage{amsmath,amssymb,color}

\definecolor{mycolorred}{rgb}{1, 0, 0}
\def\eqref #1{\textrm{(\ref#1)}}

\def\N{{\mathbb N}}

\def\<{\langle}
\def\>{\rangle}

\newtheorem{theorem}{Theorem}[section]

\newtheorem{lemma}[theorem]{Lemma}

\newtheorem{proposition}[theorem]{Proposition}
\newtheorem{remark}[theorem]{Remark}

\numberwithin{equation}{section}

\begin{document}

\title{Total variation distance\\
between stochastic polynomials\\
and invariance principles}
\author{ \textsc{Vlad Bally}\thanks{%
Universit\'e Paris-Est, LAMA (UMR CNRS, UPEMLV, UPEC), INRIA, F-77454
Marne-la-Vall\'ee, France. Email: \texttt{bally@univ-mlv.fr}.} \smallskip \\
\textsc{Lucia Caramellino}\thanks{%
Dipartimento di Matematica and INDAM-GNAMPA, Universit\`a di Roma ``Tor
Vergata'', Via della Ricerca Scientifica 1, I-00133 Roma, Italy. Email:
\texttt{caramell@mat.uniroma2.it}}\smallskip\\
}
\maketitle

\begin{abstract}
The goal of this paper is to estimate the total variation distance between
two general stochastic polynomials. As a consequence one obtains an
invariance principle for such polynomials. This generalizes known results
concerning the total variation distance between two multiple stochastic
integrals on one hand, and invariance principles in Kolmogorov distance for
multi-linear stochastic polynomials on the other hand. As an application we
first discuss the asymptotic behavior of U-statistics associated to
polynomial kernels. Moreover we also give an example of CLT associated to
quadratic forms.
\end{abstract}

\bigskip


\medskip


\tableofcontents

\section{Introduction}

This paper deals with stochastic polynomials of the following type: given a
sequence $X=(X_{n})_{n\in {\mathbb{N}}}$ of independent random variables
which have finite moments of any order and, given $N\in {\mathbb{N}}$ and $%
k_{\ast }\in {\mathbb{N}},$ one looks to
\begin{align}
Q_{N,k_{\ast }}(c,X) &=\sum_{m=0}^{N}\Phi _{m}(c,X)\quad \mbox{with}
\label{I1} \\
\Phi _{m}(c,X) &:=\sum_{k_{1},\ldots ,k_{m}=1}^{k_{\ast }}\sum_{n_{1},\ldots
,n_{m}=1}^{\infty }c((n_{1},k_{1}),\ldots
,(n_{m},k_{m}))\prod_{j=1}^{m}(X_{n_{j}}^{k_{j}}-{\mathbb{E}}%
(X_{n_{j}}^{k_{j}})).  \label{I2}
\end{align}%
The coefficients $c$ are symmetric and null on the diagonals (that is, if $%
n_{i}=n_{j}$ for $i\neq j)$ and only a finite number of them are non null,
so the above sum is finite. Let us mention that here, for notation
simplicity, we take $X_{n}\in {\mathbb{R}},$ but in the paper we work with $%
X_{n}=(X_{n,1},\ldots ,X_{n,d_{\ast }})\in {\mathbb{R}}^{d_{\ast }}.$ Note
also that we use the centred random variables $X_{n}^{k}-{\mathbb{E}}%
(X_{n}^{k})$, $k=1,\ldots,k_\ast$, but, if the polynomial is given in terms
of $X_{n}^{k},$ we may always re-write it in terms of centred random
variables.

Our goal is to estimate the total variation distance between the laws of two
such polynomials and moreover to establish an invariance principle, that is
to estimate the error done by changing $Z_{n}=(Z_{n,1},\ldots ,Z_{n,k_{\ast
}}):=(X_{n}-{\mathbb{E}}(X_{n}),\ldots ,X_{n}^{k_{\ast }}-{\mathbb{E}}%
(X_{n}^{k_{\ast }}))$ by a centred Gaussian random variable $%
G_{n}=(G_{n,1},\ldots ,G_{n,k_{\ast }})$ which has the same covariance
matrix as $Z_{n}$. Note that this Gaussian vector does not keep the
structure given by the powers in the original vector $Z_{n}.$

Since the total variation distance concerns measurable functions, a
``regularization effect'' has to be at work. This leads us to make the
following assumption (known as Doeblin's condition): there exists $%
\varepsilon >0,r>0$ and $x_{n}\in {\mathbb{R}},n\in {\mathbb{N}},$ such that
$\sup_n|x_n|<\infty$ and ${\mathbb{P}}(X_{n}\in dx)\geq \varepsilon dx$ on
the ball $B_{r}(x_{n}).$ It is easy to see that this is equivalent with
saying that%
\begin{equation}
{\mathbb{P}}(X_{n}\in dx)=p\psi (x-x_{n})dx+(1-p)\nu _{n}(dx)  \label{I3}
\end{equation}%
where $p\in (0,1],\psi $ is a $C^{\infty }$ probability density with the
support included in $B_{r}(0)$ and $\nu _{n}$ is a probability measure. The
decomposition (\ref{I3}) being given, one constructs three independent
random variable $\chi _{n},V_{n},U_{n}$ with $V_{n}\sim \psi
(x-x_{n})dx,U_{n}\sim \nu _{n}(dx)$ and $\chi _{n}$ Bernoulli with parameter
$p$ and then employs the identity of laws%
\begin{equation}
X_{n}\sim \chi _{n}V_{n}+(1-\chi _{n})U_{n}.  \label{I4}
\end{equation}%
The density $\psi $ may be chosen (see (\ref{N4})) in order that $\ln \psi $
has nice properties and this allows one to built an abstract Malliavin type
calculus based on $V_{n},n\in {\mathbb{N}}$ and to use this calculus in
order to obtain the ``regularization effect'' which is needed. We have
already used this argument in \cite{[BC-CLT],[BR],[BCNon],[BCP]}. In an
independent way, Nourdin and Poly in \cite{[NPy1]} have used similar
arguments in a similar problem: they take $\psi =(1/2r)^{-1}1_{B_{r}(0)}$ so
$V_{n}$ has a uniform distribution, and they use a chaos type decomposition
obtained in \cite{[BGL]}. Note also that hypothesis (\ref{I3}) is in fact
necessary: in his seminal paper \cite{[PROH]} Prohorov proved that (\ref{I3}%
) is (essentially) necessary and sufficient in order to obtain convergence
in total variation distance in the Central Limit Theorem (see \cite{[BC-CLT]}
for details).

The decomposition (\ref{I4}) has been introduced by Nummelin (see \cite{[N]}
and \cite{[LL]}) in order to produce atoms which allow one to use the
renewal theory for studying the convergence to equilibrium for Markov chains
-- this is why it is also known as ``the Nummelin splitting method''. It has
been also used by Poly in his PhD thesis \cite{[Py]} and, to our knowledge,
this is the first place where the idea of using the regularization given by
the noise $V_{n}$ appears.

In order to present our results we have to introduce some more notation.
Given the coefficient $c$ in (\ref{I2}) we denote%
\begin{align*}
\left\vert c\right\vert _{m} &=\Big(\sum_{k_{1},\ldots ,k_{m}=1}^{k_{\ast
}}\sum_{n_{1},\ldots ,n_{m}=1}^{\infty
}c^{2}((n_{1},k_{1}),(n_{2},k_{2}),\ldots ,(n_{m},k_{m}))\Big)^{1/2}, \\
\left\vert c\right\vert _{m,N} &=\Big(\sum_{i=m}^{N}\left\vert c\right\vert
_{i}^{2}\Big)^{1/2},\quad \left\vert c\right\vert =\left\vert c\right\vert
_{0,N}, \\
\delta _{\ast }(c) &=\max_{n}\Big(\sum_{m=0}^{N}\sum_{k_{1},\ldots
,k_{m-1}=1}^{k_{\ast }}\sum_{%
\mbox{\scriptsize{$\begin{array}{c}n_{1},\ldots ,n_{m}=1\\
 \,\exists\, i\,:\, n_i=n\end{array}$}}}^{\infty
}c^{2}((n_1,k_{1}),(n_{2},k_{2}),\ldots ,(n_{m},k_{m}))\Big)^{1/2}.
\end{align*}%
The quantity $\left\vert c\right\vert $ is essentially equivalent (up to a
multiplicative factor) with the variance of $Q_{N,k_{\ast }}(c,X)$ and $%
\delta _{\ast }(c)$ is essentially equivalent with the ``low influence
factor'' as it is defined and used in \cite{[MDO]} (and we follow several
ideas from this paper). These are the quantities which come in, in order to
estimate the errors.

For $f\in C_{b}^{\infty }({\mathbb{R}}^{d})$ we denote by $\left\Vert
f\right\Vert _{k,\infty }$ the supremum norm of $f$ and of its derivatives
of order less or equal to $k,$ and, for two random variables $F$ and $G,$ we
define the distances%
\begin{equation}  \label{TV}
d_{k}(F,G)=\sup \{\left\vert {\mathbb{E}}(f(F))-{\mathbb{E}}%
(f(G))\right\vert :\left\Vert f\right\Vert _{k,\infty }\leq 1\}.
\end{equation}%
For $k=0,$ $d_{0}=d_{\mbox{\rm{\scriptsize{TV}}}}$ is the total variation
distance, and, if $F\sim p_{F}(x)dx$ and $G\sim p_{G}(x)dx$ then $d_{%
\mbox{\rm{\scriptsize{TV}}}}(F,G)=\left\Vert p_{F}-p_{G}\right\Vert _{1}.$ $%
d_{1}$ is the Fortet-Mourier distance which metrizes the convergence in law.
We also consider the Kolmogorov distance%
\begin{equation}
d_{\mbox{\rm{\scriptsize{Kol}}}}(F,G)=\sup_{x\in {\mathbb{R}}}\left\vert {%
\mathbb{P}}(F\leq x)-{\mathbb{P}}(G\leq x)\right\vert.  \label{Kol}
\end{equation}

We are now able to give our first result, Theorem \ref{D}, concerning the
distance between two polynomials $Q_{N,k_{\ast }}(c,X)$ and $Q_{N,k_{\ast
}}(d,Y)$. Assume that $X$ and $Y$ satisfy the Doeblin's condition (see (\ref%
{I3})) and moreover assume that the non degeneracy condition $%
|c|_{m}>0,|d|_{m^{\prime }}>0$ holds for some $m,m^{\prime }\leq N$ and
denote $\overline{m}=m\vee m^{\prime }$. Then we prove (see (\ref{D3})) that
for every $k\in {\mathbb{N}}$ and $\theta \in (\frac{1}{(1+k)^{2}},1),$%
\begin{equation}
\begin{array}{l}
\displaystyle d_{\mbox{\rm{\scriptsize{TV}}}}(Q_{N,k_{\ast
}}(c,X),Q_{N,k_{\ast }}(d,Y))\leq \mathrm{Const}(c,d)\smallskip  \\
\displaystyle\times \big(d_{k}^{\frac{\theta }{2kk_{\ast }\overline{m}+1}%
}(Q_{N,k_{\ast }}(c,X),Q_{N,k_{\ast }}(d,Y))+e^{-|c|_{m}^{2}/C\delta _{\ast
}^{2}(c)}+e^{-|d|_{m^{\prime }}^{2}/C\delta _{\ast }^{2}(d)}+\left\vert
c\right\vert _{m+1,N}^{2\theta /(k_{\ast }\overline{m})}+\left\vert
d\right\vert _{m^{\prime }+1,N}^{2\theta /(k_{\ast }\overline{m})}\big),%
\end{array}
\label{I5}
\end{equation}%
where $\mathrm{Const}(c,d)$ denote a quantity which depends on the
coefficients $c$ and $d$ in an explicit way (see (\ref{D3})). If $m=N$ then $%
\left\vert c\right\vert _{m+1,N}=0$ so this term does no more appear. Theorem %
\ref{D} is the main result in our paper.

In Theorem \ref{K} we give a variant of this result in Kolmogorov distance:
we prove (see (\ref{D5''})) that
\begin{equation}
\begin{array}{l}
d_{\mbox{\rm{\scriptsize{Kol}}}}(Q_{N,k_{\ast }}(c,X),Q_{N,k_{\ast
}}(d,Y))\leq \mathrm{Const}(c,d)\smallskip  \\
\times \big(d_{k\vee 3}^{\theta /(2N(k\vee 3)+1)}(Q_{N,k_{\ast
}}(c,X),Q_{N,k_{\ast }}(d,Y))+\delta _{\ast }^{\theta /(2(k\vee
3)N+1)}(c)+\delta _{\ast }^{\theta /(2(k\vee 3)N+1)}(d)\big),%
\end{array}
\label{I6}
\end{equation}%
$\mathrm{Const}(c,d)$ is again a positive quantity explicitly depending on $c
$ and $d$ (see \ref{D5''}). The estimate (\ref{I6}) holds for general laws
for $X_{n}$ and $Y_{n}$ (without assuming the Doeblin's condition). However
now we have to assume that the covariance matrix of both $(X_{n}^{1},\ldots
,X_{n}^{k_{\ast }})$ and $(Y_{n}^{1},\ldots ,Y_{n}^{k_{\ast }})$ is
invertible. The proof of (\ref{I6}) is a direct consequence of the results
of Mossel \textit{et al.} in \cite{[MDO]}.

In the case $k_{\ast }=1$ (multilinear stochastic polynomials) and if $X_{n}$
and $Y_{n}$ are Gaussian random variables, $\Phi _{N}(c,X)$ and $\Phi
_{N}(d,Y)$ are multiple stochastic integrals. In this special case we may
drop out $e^{-1/C\delta _{\ast }^{2}(c)}$ and $e^{-1/C\delta _{\ast }^{2}(d)}
$ in (\ref{I5}) (see Theorem \ref{mul}). Estimates in total variation for
such integrals are already studied: the inequality (\ref{I5}) for multiple
stochastic integrals (for $k_{\ast }=1$) has been firstly announced in \cite%
{[D]} with the power $\frac{1}{N}$ instead of $\frac{\theta }{2N+1}$ above,
but the proof was only sketched. It has been rigourously proved in \cite%
{[NPy]} with power $\frac{1}{2N+1}$ and recently improved in \cite{[Bo]}
where the power $\frac{1}{N}\times (\ln N)^{d}$ is obtained. So (\ref{I5}) is a generalization of the
above results on multiple stochastic integrals to general polynomials
depending on a general noise. But, as the above discussion suggests, (\ref%
{I5}) is not the best possible estimate (the approach in \cite{[Bo]} does
not seem to work in our general framework, so for the moment we are not able
to improve it).

A second result, given in Theorem \ref{Main}, concerns the invariance
principle. We consider a sequence of independent centred Gaussian random
variables $G_{n}=(G_{n,1},\ldots ,G_{n,k_{\ast }})\in {\mathbb{R}}^{k_{\ast
}}$ and we assume that the covariance matrix of $G_{n}$ coincides with the
covariance matrix of $Z_{n}=(Z_{n,1},\ldots ,Z_{n,k_{\ast }})$ where $%
Z_{n,k}:=X_{n}^{k}-{\mathbb{E}}(X_{n}^{k}).$ We denote by $S_{N}(c,G)$ the
polynomial $Q_{N,k_{\ast }}(c,X)$ in which $Z_{n}=(Z_{n,1},\ldots
,Z_{n,k_{\ast }})$ is replaced by $G_{n}=(G_{n,1},\ldots ,G_{n,k_{\ast }}).$
We stress that $S_{N}(c,G)$ is multi-linear with respect to $%
G_{n,i},i=1,\ldots ,k_{\ast }$ in contrast to $Q_{N,k_{\ast }}(c,X)$ which
is a general polynomial with respect to $X_{n}.$ In Theorem \ref{Main} we
prove that, if $|c|_{m}>0,$ for some $m\leq N,$ then for every $\theta \in (%
\frac{1}{16},1)$,%
\begin{equation}
d_{{\mbox{\rm{\scriptsize{TV}}}}}(Q_{N,k_{\ast }}(c,X),S_{N}(c,G))\leq
\mathrm{Const}(c)\big(\delta _{\ast }^{\theta /(6k_{\ast
}m+1)}(c)+e^{-|c|_{m}^{2}/C\delta _{\ast }^{2}(c)}+\left\vert c\right\vert
_{m+1,N}^{2\theta /(k_{\ast }m)}\big),  \label{I7}
\end{equation}%
$\mathrm{Const}(c)$ being explicitly dependent on $c$ (see (\ref{MR1}). A
result going in the same direction was previously obtained by Nourdin and
Poly in \cite{[NPy1]}. They take $k_{\ast }=1$, so $Q_{N}(c,X)$ is a
multi-linear polynomial, and they assume Doeblin's condition for $X_{i}.$
Then they prove that, if $c_{n},n\in {\mathbb{N}}$ is a sequence of
coefficients such that $\lim_{n}\delta _{\ast }(c_{n})=0,$ then $\lim_{n}d_{{%
\mbox{\rm{\scriptsize{TV}}}}}(Q_{N,k_{\ast }}(c,X),S_{N}(c,G))=0.$ The
progress achieved in our paper consists in the fact that we deal with
general polynomials on one hand and we obtain an estimate of the error on
the other hand.

A similar estimate with $d_{\mbox{\rm{\scriptsize{Kol}}}}$ instead of $d_{%
\mbox{\rm{\scriptsize{TV}}}}$ represents the main result in \cite{[MDO]}
(see Theorem 3.19 therein). Let us be more precise. In \cite{[MDO]} one
considers ``orthonormal ensembles'' which are nothing else than
multi-dimensional random variables $Z_{n}=(Z_{n,1},\ldots ,Z_{n,k_{\ast }})$
such that ${\mathbb{E}}(Z_{n,i})=0$ and ${\mathbb{E}}(Z_{n,i}Z_{n,j})=\delta
_{i,j}$ (the Kronecker delta). One denotes $S_{N}(c,Z)$ the polynomial $%
Q_{N,k_{\ast }}(c,X)$ defined (\ref{I1}) in which $X_{n}^{k}-{\mathbb{E}}%
(X_{n}^{k})$ is replaced by $Z_{n,k}.$ And in \cite{[MDO]} (Theorem 3.19
therein) they prove that if $\left\vert c\right\vert =1,$ then
\begin{equation}
d_{\mbox{\rm{\scriptsize{Kol}}}}(S_{N}(c,Z),S_{N}(c,G))\leq C\times \delta
_{\ast }^{1/(3N+1)}(c).  \label{I8}
\end{equation}%
Note that in this theorem one does not need Doeblin condition to hold true.
Note also that the orthonormality condition for $Z_{n,1},\ldots
,Z_{n,k_{\ast }}$ is not more restrictive than saying that the covariance
matrix $\mathrm{Cov}(Z_{n})$ of $Z_{n}$ is invertible and the lower
eigenvalues $\lambda _{n}$ satisfy $\lambda _{n}\geq \underline{\lambda }>0$
for every $n$ (see the proof of Theorem \ref{MDO}). So, by taking $%
Z_{n,k}:=X_{n}^{k}-{\mathbb{E}}(X_{n}^{k}),$ one obtains also (\ref{I7})
(under the above hypothesis on $\mathrm{Cov}(Z_{n})).$ The difference with
respect to their result is just that we deal with convergence in total
variation distance instead of Kolmogorov distance.

An important consequence of (\ref{I7}) is that it allows to replace the
study of the asymptotic behavior of a sequence $Q_{N,k_{\ast
}}(c_{n},X),n\in {\mathbb{N}}$ of general stochastic polynomials by the
study of $S_{N}(c_{n},G),n\in {\mathbb{N}},$ which are elements of a finite
number of Wiener chaoses. Of course, the central example is the classical
CLT, where $N=1$ and $k_{\ast }=1$, so $S_{1}(c_{n},G)=\sum_{i=1}^{\infty
}c_{n}(i)G_{i}$ is just a Gaussian random variable. But, starting with the
proof of the ``forth moment theorem'' by Nualart and Peccati \cite{[PN]} and
Nourdin and Peccati \cite{[NP1]}, a lot of work has been done in order to
characterize the convergence to normality of elements of a finite number of
Wiener chaoses (see \cite{[NN], [NPRev], [NO],[PT]} or \cite{[NP]} for an overview). Moreover, convergence to a
$\chi _{2}$ distribution has been treated in \cite{[NP1]}. We give the
consequences of these results in Theorem \ref{N} and Theorem \ref{G}.

Finally we give two more applications. The first one concerns U-statistics.
The problem is the following: given a probability law $\mu ,$ an integer $%
N\in {\mathbb{N}},$ and a symmetric kernel $\psi ,$ one wants to estimate%
\begin{equation*}
\theta (\mu )=\int_{{\mathbb{R}}^{N}}\psi (x_{1},\ldots ,x_{N})d\mu
(x_{1})\ldots d\mu (x_{N})
\end{equation*}%
on the basis of a sample $X_{1},\ldots .,X_{n}$ of independent random
variables of law $\mu .$ An un-biased estimator of $\theta (\mu )$ is
constructed by%
\begin{equation*}
U_{n}^{\psi }=\frac{(n-N)!}{n!}\sum_{i_{1},\ldots ,i_{N}=1}^{n}\delta
(i_{1},\ldots ,i_{N})\psi (X_{i_{1}},\ldots ,X_{i_{N}}),
\end{equation*}%
in which $\delta (i_{1},\ldots ,i_{N})=0$ if any two indexes are equal,
otherwise $\delta (i_{1},\ldots ,i_{N})=1$. In the case when $\psi $ is a
polynomial this enters in our framework. This covers an important class of
kernels: for example $\psi (x_{1},x_{2})=(x_{1}-x_{2})^{2}$ gives the
estimator of the variance. But not all: for example $\psi (x_{1},\ldots
,x_{N})=\max_{i=1,N}\left\vert x_{i}\right\vert $ is out of reach. Say that $%
\psi (x_{1},\ldots ,x_{N})=\sum_{k_{1},\ldots ,k_{N}=1}^{k_{\ast }}\delta
(i_{1},\ldots ,i_{N})b(k_{1},\ldots ,k_{N})\prod_{j=1}^{N}x_{j}^{k_{j}}.$
Then%
\begin{equation*}
U_{n}^{\psi }=\frac{(n-N)!}{n!}\sum_{i_{1},\ldots
,i_{N}=1}^{n}\sum_{k_{1},\ldots ,k_{N}=1}^{k_{\ast }}\delta (i_{1},\ldots
,i_{N})b(k_{1},\ldots ,k_{N})\prod_{j=1}^{N}X_{i_{j}}^{k_{j}}.
\end{equation*}%
This fits in (\ref{I1}) except that $X_{j}^{k_{j}}$ is not centred. It turns
out that the procedure which consists in centering $X_{j}^{k_{j}}$
coincides, in this framework, with the Hoeffding's decomposition, which is a
central tool in the U-statistics theory. After doing this one obtains%
\begin{align*}
U_{n}^{\psi }-\theta (\mu )& =\sum_{m=1}^{N}\frac{(n-m)!}{n!}\Phi _{m}(a,X)
\\
& =\sum_{m=1}^{N}\frac{(n-m)!}{n!}\sum_{i_{1},\ldots
,i_{m}=1}^{n}\sum_{k_{1},\ldots ,k_{m}=1}^{k_{\ast }}c((n_{1},k_{1}),\ldots
,(n_{m},k_{m}))\prod_{j=1}^{m}(X_{i_{j}}^{k_{j}}-{\mathbb{E}}%
(X_{i_{j}}^{k_{j}}))
\end{align*}%
for some appropriate coefficients $c((n_{1},k_{1}),\ldots ,(n_{m},k_{m}))$,
and we are back in our framework. In U-statistics theory one says that the
kernel $\psi $ is degenerated at order $m_{0}$ if $\Phi _{m}=0$ for $m\leq
m_{0}-1$ and $\Phi _{m_{0}}\neq 0.$ Then one writes%
\begin{equation*}
n^{m_{0}}(U_{n}^{\psi }-\theta (\mu ))=C\times \Phi _{m_{0}}(a,X)+R_{n}
\end{equation*}%
with $R_{n}\rightarrow 0.$ It follows that the asymptotic behavior of $%
n^{m_{0}}(U_{n}^{\psi }-\theta (\mu ))$ is controlled by $\Phi
_{m_{0}}(a,X). $ Using this decomposition, in Theorem \ref{UU} we
characterizes the limit of $n^{m_{0}}(U_{n}^{\psi }-\theta (\mu ))$ as a
linear combination of multiple stochastic integrals. The limit is considered
both in Kolmogorov distance under general conditions and in total variation
distance under Doeblin condition for $\mu $. Let us mention that number of
results are already known concerning the convergence in Kolmogorov distance
for U-statistics: they represent generalizations of the Berry--Essen theorem
(we refer to \cite{[Lee]} and \cite{[kb]}). But the result in total
variation distance, which generalizes Prohorov's theorem for the CLT, seems
to be new.

Another subject which is very closed, is that of quadratic forms. Here also
the asymptotic behavior in Kolmogorov distance is well understood (see de
Jong \cite{[dJ1],[dJ2]} , Rotar' \textit{et al.} \cite{[GR],[RS]} and G\"{o}%
tze \textit{et al.} \cite{[GT]}) but we have not found results concerning
the convergence in total variation. We do not treat this subject in all
generality but we restrict ourselves to the following interesting example:
for $p\in \lbrack 0,\frac{1}{2}]$ we define%
\begin{equation*}
S_{n,p}=\varepsilon _{p}(n)\sum_{1\leq i<j\leq n}\frac{1}{\left\vert
j-i\right\vert ^{p}}X_{i}X_{j}
\end{equation*}%
where $X_{i},i\in {\mathbb{N}}$ are independent identically distributed
random variables with ${\mathbb{E}}(X_{i})=0$ and ${\mathbb{E}}%
(X_{i}^{2})=1. $ And $\varepsilon _{p}(n)=n^{-(1-p)}$ for $p<\frac{1}{2}$
and $\varepsilon _{1/2}(n)=1/\sqrt{2n\ln n}.$ For $p<\frac{1}{2}$ we prove
that $S_{n,p}\rightarrow \int_{0}^{1}\int_{0}^{t}(t-s)^{-p}dW_{s}dW_{t}$ and
for $p=\frac{1}{2}$ one has $S_{n,p}\rightarrow \Delta $ with $\Delta $ a
standard normal random variable. Thus, there is a change of regime in $p=%
\frac{1}{2}.$ As before, the convergence takes place in Kolmogorov distance
for a general $X$ and in total variation distance under Doeblin's condition.

\medskip The paper is organized as follows. In Section \ref{notations}, we
fix our settings and we give some preliminary results. Section \ref{sect:3}
is devoted to our main results: we first precisely define the Doeblin's
condition and the Nummelin splitting (Section \ref{sect:3.1}); then we
introduce our main result Theorem \ref{D} and its several consequences
(Section \ref{sect:3.2}); finally we analyze the Gaussian and Gamma
approximation (Section \ref{sect:3.3}). The main examples are developed in
Section \ref{sect:examples}: in Section \ref{sect:ustatistics} we study the
asymptotic behavior of U-statistics written on polynomial kernels and in
Section \ref{sect:quadratic} we study the convergence of the above quadratic
CLT result. Finally, Section \ref{sect:doeblin} contains the proof of our
main Theorem \ref{D}, which is given in the last Section \ref{sect:proof}:
in Section \ref{sect:sobolev} we introduce the abstract Malliavin calculus,
in Section \ref{sect:reg} we state the regularization lemma we use in this
paper, Section \ref{sect:est} is devoted to proper estimates of the Sobolev
norms and Section \ref{sect:cov} refers to the non-degeneracy result of the
Malliavin covariance matrix. The paper concludes with two appendixes:
Appendix \ref{hoef} studies an iterated Hoeffding's inequality for
martingales and Appendix \ref{app:norms} gives useful estimates for the
Sobolev norms which are used the Malliavin integration by parts formula.

\bigskip

\noindent \textbf{Acknowledgments.} We thank to Cristina Butucea and to Dan
Timotin for useful discussions.

\section{Notation, basic objects and preliminary results}

\label{notations}

In this section we introduce multi-linear stochastic polynomials based on a
sequence of abstract independent random variables $Z_{n}=(Z_{n,1},\ldots
,Z_{n,m_{\ast }})\in {\mathbb{R}}^{m_{\ast }},$ $n\in {\mathbb{N}}.$ In the
next section, when dealing with general polynomials as in (\ref{I1}), we
will take $Z_{n,k}=X_{n}^{k}-{\mathbb{E}}(X_{n}^{k}).$

\medskip

$\square $ \textbf{The basic noise}. We assume that ${\mathbb{E}}(Z_{n,i})=0$
and that $Z_{n} $ has finite moments of any order: for every $p\geq 1$ there
exists some $M_{p}(Z)\geq 1$ such that for every $n\in {\mathbb{N}}$ and $%
i\in \lbrack m_{\ast }]=\{1,\ldots ,m_{\ast }\}$%
\begin{equation}
\left\Vert Z_{n,i}\right\Vert _{p}\leq M_{p}(Z).  \label{H1}
\end{equation}

$\square $ \textbf{Multi-indexes}. We will use ``double'' multi-indexes $\alpha =(\alpha _{1},\ldots ,\alpha
_{m})$ with $\alpha _{i}=(\alpha _{i}^{\prime },\alpha _{i}^{\prime \prime
})=(n_{i},j_{i})$ with $n_{i}\in {\mathbb{N}}$ and $j_{i}\in \lbrack m_{\ast
}].$ We always assume that $n_{1}<\ldots <n_{m}.$ So we work with "ordered"
multi-indexes. We also denote $\alpha ^{\prime }=(\alpha _{1}^{\prime
},\ldots ,\alpha _{m}^{\prime })=(n_{1},\ldots ,n_{m})$, $\alpha ^{\prime
\prime }=(\alpha _{1}^{\prime \prime },\ldots ,\alpha _{m}^{\prime \prime
})=(j_{1},\ldots ,j_{m})$ and $\left\vert \alpha \right\vert =m.$ The set of
such multi-indexes is denoted by $\Gamma _{m}$ and we set $\Gamma =\cup
_{m}\Gamma _{m}$. We stress that we consider also the void multi-index $%
\alpha =\emptyset $ and in this case we put $\left\vert \alpha \right\vert
=0.$ Moreover, for a sequence $x_{n}=(x_{n,1},\ldots ,x_{n,m_{\ast }})\in {%
\mathbb{R}}^{m_{\ast }},n\in {\mathbb{N}}$ we denote%
\begin{equation*}
x^{\alpha }=\prod_{i=1}^{m}x_{\alpha _{i}},
\end{equation*}%
with $x^{\alpha }=1$ if $\alpha =\emptyset $.

\medskip

$\square $ \textbf{Coefficients}. We consider a Hilbert space $\mathcal{U}$
with norm $\left\vert \cdot \right\vert _{\mathcal{U}}$ and for a $\mathcal{U%
}$ valued random variable $X$, we denote $\left\Vert X\right\Vert _{\mathcal{%
U},p }=({\mathbb{E}}(\left\vert X\right\vert _{\mathcal{U}}^{p})^{1/p}.$ In
a first stage we have just $\mathcal{U}={\mathbb{R}}$ but in Section \ref%
{sect:doeblin}, when considering stochastic derivatives, we have to use some
general space $\mathcal{U}$. We denote $\mathcal{C(U})=\{c=(c(\alpha
))_{\alpha \in \Gamma }:c(\alpha )\in \mathcal{U\}}$. These are the
coefficients we will use. We define%
\begin{equation}
\begin{array}{c}
\displaystyle\left\vert c\right\vert _{\mathcal{U}}=\Big(\sum_{\alpha
}\left\vert c(\alpha )\right\vert _{\mathcal{U}}^{2}\Big)^{1/2},\quad
\left\vert c\right\vert _{\mathcal{U},m}=\Big(\sum_{\left\vert \alpha
\right\vert =m}\left\vert c(\alpha )\right\vert _{\mathcal{U}}^{2}\Big)%
^{1/2}\smallskip \\
\displaystyle\mathcal{N}_{{\mathcal{U}},q}(c,M)=\Big(\sum_{m=0}^{\infty
}m^{q}M^{2m}\left\vert c\right\vert _{\mathcal{U},m}^{2}\Big)%
^{1/2}=\sum_{\alpha }\left\vert \alpha \right\vert ^{m}M^{2\left\vert \alpha
\right\vert }\left\vert c(\alpha )\right\vert _{\mathcal{U}}^{2}%
\end{array}
\label{H2a}
\end{equation}%
and%
\begin{equation}
\delta _{{\mathcal{U}},\ast }(c)=\Big(\sup_{n}(\sum_{\alpha }1_{\{n\in
\alpha ^{\prime }\}}\left\vert c(\alpha )\right\vert _{\mathcal{U}}^{2})\Big)%
^{1/2}.  \label{H2b}
\end{equation}%
The notation $n\in \alpha ^{\prime }$ means that $\alpha _{j}^{\prime }=n$
for some $j\in \lbrack m].$ When ${\mathcal{U}}={\mathbb{R}}$, we shall omit
the subscript ${\mathcal{U}}$, so we simply write $|c|$, $|c|_{m}$, $%
\mathcal{N}_{q}(c,M)$ and $\delta _{\ast }(c)$. For several authors (see
e.g. \cite{[MDO]} or \cite{[NPRein]}), $\delta _{{\mathcal{U}},\ast }^{2}(c)$
is called the ``influence'' factor.

\medskip

$\square $ \textbf{Multi-linear polynomials}. Given $c\in $ $\mathcal{C(U}) $
we define%
\begin{eqnarray}
\Phi _{m}(c,Z) &=&\sum_{\left\vert \alpha \right\vert =m}c(\alpha )Z^{\alpha
}=\sum_{j_{1},\ldots ,j_{m}=1}^{m_{\ast }} \sum_{n_{1}<\cdots
<n_{m}}c((n_{1},j_{1}),\ldots
,(n_{m},j_{m}))\prod_{i=1}^{m}Z_{n_{i},j_{i}},\quad  \label{H2'} \\
S_{N}(c,Z) &=&\sum_{0\leq \left\vert \alpha \right\vert \leq N}c(\alpha
)Z^{\alpha }=\sum_{m=0}^{N}\Phi _{m}(c,Z).  \label{H2''}
\end{eqnarray}

In the sequel we use several times Burkholder's inequality for Hilbert space
valued martingales: if $M_{n}\in \mathcal{U},n\in {\mathbb{N}}$ is a
martingale then for every $p\geq 2$ there exists $b_{p}\geq 1$ such that%
\begin{equation}
\left\Vert M_{n}\right\Vert _{\mathcal{U},p}\leq b_{p}\Big({\mathbb{E}}\Big(%
\Big(\sum_{k=1}^{n-1}\left\vert M_{k+1}-M_{k}\right\vert _{\mathcal{U}}^{2}%
\Big)^{p/2}\Big)\Big)^{1/p}\leq b_{p}\Big(\sum_{k=1}^{n-1}\left\Vert
M_{k+1}-M_{k}\right\Vert _{\mathcal{U},p}^{2}\Big)^{1/2}  \label{H5}
\end{equation}%
the second inequality being obtained by using the triangle inequality with
respect to $\left\Vert \cdot \right\Vert _{\mathcal{U},p/2}.$

Moreover, as an immediate consequence of (\ref{H1}), for every $n\in {%
\mathbb{N}}$ and every $d_{j}\in \mathcal{U},j\in \lbrack m_{\ast }]$ we have%
\begin{equation}
\Big\Vert \sum_{j=1}^{m_{\ast }}d_{j}\times Z_{n,j}\Big\Vert _{\mathcal{U}%
,p}\leq \sqrt{m_{\ast }}M_{p}(Z)\Big(\sum_{j=1}^{m_{\ast }}\left\vert
d_{j}\right\vert _{\mathcal{U}}^{2}\Big)^{1/2}.  \label{H6}
\end{equation}

Using these two inequalities we obtain

\begin{lemma}
\label{B}Suppose that (\ref{H1}) holds and denote $\overline{M}%
_{p}=b_{p}M_{p}(Z)\sqrt{m_{\ast }}.$ Then
\begin{equation}
\left\Vert \Phi _{N}(c,Z)\right\Vert _{\mathcal{U},p}\leq \overline{M}%
_{p}^{N}\left\vert c\right\vert _{\mathcal{U},N}  \label{H7}
\end{equation}%
and%
\begin{equation}
\left\Vert S_{N}(c,Z)-c(\emptyset )\right\Vert _{\mathcal{U},p}\leq \mathcal{%
N}_{{\mathcal{U}},0}(c,\overline{M}_{p}).  \label{H8}
\end{equation}
\end{lemma}

\textbf{Proof}. We proceed by recurrence on $N.$ For $N=0$ we have $\Phi
_{N}(c,Z)=c(\emptyset )$ so (\ref{H7}) is obvious. For $\alpha\in\Gamma $
with $\left\vert \alpha \right\vert =N-1$ we denote
\begin{equation}
c^{n,j}(\alpha )=c(\alpha ,(n,j))1_{\{\alpha _{N-1}^{\prime }<n\}}
\label{H8'}
\end{equation}%
and we write%
\begin{equation}
\Phi _{N}(c,Z)=\sum_{n=N}^{\infty }\sum_{j=1}^{m_{\ast
}}Z_{n,j}\sum_{\left\vert \alpha \right\vert =N-1}c(\alpha
,(n,j))1_{\{\alpha _{N-1}^{\prime }<n\}}Z^{\alpha } =\sum_{n=N}^{\infty
}\sum_{j=1}^{m_{\ast }}Z_{n,j}\Phi _{N-1}(c^{n,j},Z).  \label{H8''}
\end{equation}%
Note that, if $n\geq N$, $Z_{n,j}$ and $\Phi _{N-1}(c^{n,j},Z)$ are
independent. So, using (\ref{H5}) first and (\ref{H6}) then we get%
\begin{equation*}
\left\Vert \Phi _{N}(c,Z)\right\Vert _{\mathcal{U},p}^{2}\leq
b_{p}^{2}\sum_{n=N}^{\infty }\Big\Vert \sum_{j=1}^{m_{\ast }}Z_{n,j}\Phi
_{N-1}(c^{n,j},Z)\Big\Vert _{\mathcal{U},p}^{2}\leq
b_{p}^{2}M_{p}^{2}(Z)m_{\ast }\sum_{n=N}^{\infty }\sum_{j=1}^{m_{\ast
}}\left\Vert \Phi _{N-1}(c^{n,j},Z)\right\Vert _{\mathcal{U},p}^{2}
\end{equation*}%
and by the recurrence hypothesis,
\begin{equation*}
\left\Vert \Phi _{N}(c,Z)\right\Vert _{\mathcal{U},p}^{2}\leq
(b_{p}^{2}M_{p}^{2}(Z)m_{\ast })^{N}\sum_{n=N}^{\infty }\sum_{j=1}^{m_{\ast
}}\left\vert c^{n,j}\right\vert _{\mathcal{U},
N-1}^{2}=(b_{p}^{2}M_{p}^{2}(Z)m_{\ast })^{N}\sum_{\left\vert \alpha
\right\vert =N}\left\vert c(\alpha )\right\vert _{\mathcal{U}}^{2}.
\end{equation*}%
So (\ref{H7}) is proved.

We now prove (\ref{H8}) again by induction. The case $N=1$ follows from (\ref%
{H7}). For $N\geq 2$, we have
\begin{align*}
S_N(c,Z)-c(\emptyset) &=\sum_{m=1}^N\Phi_m(c,Z) =\sum_{m=1}^N
\sum_{n=m}^{\infty }\sum_{j=1}^{m_{\ast }}Z_{n,j}\Phi _{m-1}(c^{n,j},Z) \\
&=\sum_{n=1}^{\infty }\sum_{j=1}^{m_{\ast }}Z_{n,j}\sum_{m=1}^{N\wedge n}
\Phi _{m-1}(c^{n,j},Z) =\sum_{n=1}^{\infty }\sum_{j=1}^{m_{\ast
}}Z_{n,j}S_{N\wedge n-1}(c^{n,j},Z)
\end{align*}
If $n\geq m$, $Z_{n,j}$ and $\Phi _{m-1}(c^{n,j},Z)$ are independent, so $%
Z_{n,j}$ and $S_{N\wedge n-1}(c^{n,j},Z)$ are independent as well. Therefore
we can apply (\ref{H5}) and (\ref{H6}) and we obtain
\begin{align*}
\|S_N(c,Z)-c(\emptyset)\|_{{\mathcal{U}},p}^2 &\leq b_p^2\sum_{n=1}^{\infty }%
\Big\|\sum_{j=1}^{m_{\ast }}Z_{n,j}S_{N\wedge n-1}(c^{n,j},Z)\Big\|_{{%
\mathcal{U}},p}^2 \leq \overline{M}^2_p\sum_{n=1}^{\infty
}\sum_{j=1}^{m_{\ast }}\|S_{N\wedge n-1}(c^{n,j},Z)\|_{{\mathcal{U}},p}^2
\end{align*}
and by the recurrence hypothesis,
\begin{align*}
\|S_N(c,Z)-c(\emptyset)\|_{{\mathcal{U}},p}^2 &\leq \overline{M}%
^2_p\sum_{n=1}^{\infty }\sum_{j=1}^{m_{\ast }}\mathcal{N}_{{\mathcal{U}}%
,0}^2(c^{n,j},\overline{M}_p) \leq \mathcal{N}_{{\mathcal{U}},0}^2(c,%
\overline{M}_p).
\end{align*}
$\square $

\medskip

We give now the basic invariance principle. We take $\mathcal{U}={\mathbb{R}}%
,$ and for $f\in C_{b}^{3}({\mathbb{R}}),$ we denote by $\left\Vert
f\right\Vert _{3,\infty }$ the supremum norm of $f$ and its derivatives up
to order three.

\begin{theorem}
\label{smooth}Let $Z=(Z_{n})_{n\in {\mathbb{N}}},Z_{n}\in {\mathbb{R}}%
^{m_{\ast }}$ be a sequence of centred independent random variables which
verify (\ref{H1}) and let $G=(G_{n})_{n\in {\mathbb{N}}},G_{n}\in {\mathbb{R}%
}^{m_{\ast }}$ be a sequence of independent centred Gaussian random
variables such that ${\mathbb{E}}(G_{n,i}G_{n,j})={\mathbb{E}}%
(Z_{n,i}Z_{n,j}).$ Then, for every $f\in C_{b}^{3}({\mathbb{R}})$
\begin{equation}
\left\vert {\mathbb{E}}(f(S_{N}(c,Z))-{\mathbb{E}}(f(S_{N}(c,G))\right\vert
\leq {\mathcal{K}}_{N,m_{\ast }}(Z)\left\Vert f\right\Vert _{3,\infty
}\times \left\vert c\right\vert ^{2}\times \delta _{\ast }(c)  \label{H11}
\end{equation}%
with%
\begin{equation*}
{\mathcal{K}}_{N,m_{\ast }}(Z)=\frac{2m_{\ast }}{3}%
(M_{3}^{3}(Z)+M_{3}^{3}(G))\overline{M}_{3}^{3N},
\end{equation*}
in which $\overline{M}_{3}=b_{3}\sqrt{m_{\ast }}\, M_{3}(Z)\vee M_{3}(G).$
\end{theorem}

\textbf{Proof}. The proof is based on Lindeberg's method (we follow the
argument from \cite{[MDO]}). We fix $J\geq N,$ we denote $\Gamma
_{N}(J)=\cup _{m=0}^{N}\{\alpha \in \Gamma :\left\vert \alpha \right\vert
=m,\alpha _{m}^{\prime }\leq J\}$ and we define $S_{N,J}(c,Z)=\sum_{\alpha
\in \Gamma _{N}(J)}c(\alpha )Z^{\alpha }.$ For $j=1,\ldots,J+1$ we define
the intermediate sequences $Z^{j}=(Z_{1},\ldots ,Z_{j-1},G_{j},\ldots
,G_{J}) $, with $Z^{1}=(G_{1},\ldots ,G_{J})$ and $Z^{J+1}=(Z_{1},\ldots
,Z_{J})$, and we write%
\begin{equation*}
{\mathbb{E}}(f(S_{N,J}(c,Z))-{\mathbb{E}}(f(S_{N,J}(c,G))=\sum_{j=1}^{J}{%
\mathbb{E}}(f(S_{N,J}(c,Z^{j+1}))-{\mathbb{E}}(f(S_{N,J}(c,Z^{j}))=:%
\sum_{j=1}^{J}I_{j}.
\end{equation*}%
We denote $\Gamma _{N}(j,J)=\{\alpha \in \Gamma _{N}(J):j\notin \alpha
^{\prime }\}$ and, for $\beta \in \Gamma _{N}(j,J)$ with $\left\vert \beta
\right\vert =m$ we define
\begin{eqnarray*}
c_{j,i}(\beta ) &=&\sum_{k=2}^{m}c(\beta _{1},\ldots, \beta
_{k-1},(j,i),\beta _{k},\ldots ,\beta _{m})1_{\{\beta _{k-1}^{\prime
}<j<\beta _{k}^{\prime }\}} \\
&&+c((j,i),\beta _{1},\ldots \beta _{m},(j,i))1_{\{j<\beta _{1}^{\prime
}\}}+c(\beta _{1},\ldots \beta _{m},(j,i))1_{\{\beta _{m}^{\prime }<j\}}.
\end{eqnarray*}%
This means that, if $\beta $ does not contain $j,$ we insert $(j,i)$ in the
convenient position. We put%
\begin{equation*}
A_{j}=\sum_{\alpha \in \Gamma _{N}(j,J)}c(\alpha )(Z^{j})^{\alpha },\quad
B_{j,i}=\sum_{\beta \in \Gamma _{N-1}(j,J)}c_{j,i}(\beta )(Z^{j})^{\beta }
\end{equation*}%
and then
\begin{equation*}
S_{N,J}(c,Z^{j+1})=A_{j}+\sum_{i=1}^{m_{\ast }}Z_{j,i}B_{j,i}.
\end{equation*}%
Moreover, with $f_{j}:{\mathbb{R}}^{m_{\ast }}\rightarrow {\mathbb{R}}$
defined by $f_{j}(x):=f(A_{j}+\sum_{i=1}^{m_{\ast }}x_{i}B_{j,i})$ we get%
\begin{equation*}
I_{j}={\mathbb{E}}(f(S_{N,J}(c,Z^{j+1}))-{\mathbb{E}}(f(S_{N,J}(c,Z^{j}))={%
\mathbb{E}}(f_{j}(Z_{j}))-{\mathbb{E}}(f_{j}(G_{j})).
\end{equation*}

We use now Taylor's expansion of order three around $0$ for both $%
f_{j}(Z_{j})$ and $f_{j}(G_{j})$. Since $Z_{j}$ and $G_{j}$ are independent
of $A_{j}$ and $B_{j,\cdot}$ and the first and second moments of $Z_{j,i}$
and $G_{j,i}$ coincide, the first and second order terms in the Taylor
expansion cancel and we obtain%
\begin{equation*}
\left\vert I_{j}\right\vert \leq \frac{1}{2}\sum_{i_{1},i_{2},i_{3}=1}^{m_{%
\ast }}{\mathbb{E}}\Big(\prod_{r=1}^{3}(\left\vert Z_{j,i_{r}}\right\vert
+\left\vert G_{j,i_{r}}\right\vert )\int_{0}^{1}(1-\lambda )^{2}(|\partial
_{i_{1}i_{2}i_{3}}^{3}f_{j}(\lambda Z_{j})|+|\partial
_{i_{1}i_{2}i_{3}}^{3}f_{j}(\lambda G_{j})|)d\lambda \Big).
\end{equation*}%
We have
\begin{equation*}
|\partial _{i_{1}i_{2}i_{3}}^{3}f_{j}(\lambda Z_{j})|=\vert
f_{j}^{(3)}(\lambda Z_{j})\vert \times \prod_{r=1}^{3}\left\vert
B_{j,i_{r}}\right\vert \leq \left\Vert f\right\Vert _{3,\infty }\times
\prod_{r=1}^{3}\left\vert B_{j,i_{r}}\right\vert .
\end{equation*}%
The same is true for $|\partial _{i_{1}i_{2}i_{3}}^{3}f_{j}(\lambda G_{j})|$%
, so (recall that $Z_{j}$ and $G_{j}$ are independent of $B_{j,\cdot})$%
\begin{equation}
\left\vert I_{j}\right\vert \leq \frac{1}{3}\left\Vert f\right\Vert
_{3,\infty }(M_{3}^{3}(Z)+M_{3}^{3}(G))\sum_{i_{1},i_{2},i_{3}=1}^{m_{\ast }}%
{\mathbb{E}}\Big(\prod_{r=1}^{3}\left\vert B_{j,i_{r}}\right\vert \Big).
\label{H10}
\end{equation}%
Using (\ref{H8}),
\begin{equation*}
\left\Vert B_{j,i}\right\Vert _{3}\leq \overline{M}_{3}^{N}(\sum_{\beta \in
\Gamma _{N-1}(j,J)}\left\vert c_{j,i}(\beta )\right\vert ^{2})^{1/2}\leq
\overline{M}_{3}^{N}\delta _{\ast }(c)
\end{equation*}%
and this gives%
\begin{equation*}
{\mathbb{E}}(\prod_{r=1}^{3}\left\vert B_{j,i_{r}}\right\vert )\leq
\prod_{r=1}^{3}\left\Vert B_{j,i_{r}}\right\Vert _{3}\leq \overline{M}%
_{3}^{N}\delta _{\ast }(c)(\left\Vert B_{j,i_{1}}\right\Vert
_{3}^{2}+\left\Vert B_{j,i_{2}}\right\Vert _{3}^{2}).
\end{equation*}%
We sum over $j$ and we get
\begin{eqnarray*}
\sum_{j=1}^{J}\left\vert I_{j}\right\vert &\leq &\frac{2m_{\ast }}{3}%
\left\Vert f\right\Vert _{3,\infty }(M_{3}^{3}(Z)+M_{3}^{3}(G))\overline{M}%
_{3}^{N}\delta _{\ast }(c)\sum_{j=1}^{J}\sum_{i=1}^{m_{\ast }}\left\Vert
B_{j,i}\right\Vert _{3}^{2} \\
&\leq &\frac{2m_{\ast }}{3}\left\Vert f\right\Vert _{3,\infty
}(M_{3}^{3}(Z)+M_{3}^{3}(G))\overline{M}_{3}^{3N}\delta _{\ast
}(c)\left\vert c\right\vert ^{2}
\end{eqnarray*}%
$\square $

\medskip

We recall now the main result from \cite{[MDO]} concerning the invariance
principle in Kolmogorov distance (defined in (\ref{Kol})).

\begin{theorem}
\label{MDO}Let $Z=(Z_{n})_{n\in {\mathbb{N}}},Z_{n}\in {\mathbb{R}}^{m_{\ast
}}$ be a sequence of centred independent random variables which verify (\ref%
{H1}) and let $\mathrm{Cov}(Z_{n})$ denote the covariance matrix of $Z_{n}.$
We assume that there exists $0<\underline{\lambda }\leq 1$ such that for
every $n\in {\mathbb{N}}$
\begin{equation}
\mathrm{Cov}(Z_{n})\geq \underline{\lambda }.  \label{H2}
\end{equation}%
Let $G=(G_{n})_{n\in {\mathbb{N}}},G_{n}\in {\mathbb{R}}^{m_{\ast }}$ be a
sequence of independent centred Gaussian random variables such that $\mathrm{%
Cov}(Z_{n})=\mathrm{Cov}(G_{n}).$ Then%
\begin{equation}
d_{\mbox{\rm{\scriptsize{Kol}}}}(S_{N}(c,Z),S_{N}(c,G))\leq {\mathcal{K}}%
_{N}(Z)\times \delta _{\ast }^{1/(1+3N)}(c)  \label{H12}
\end{equation}%
with%
\begin{equation*}
{\mathcal{K}}_{N}(Z)=C\times N^{1/(3N+1)}(b_{3}(b_{3}\underline{\lambda }%
^{-m_{\ast }}M_{3}(Z)^{N})^{3N/(3N+1)}\times ((m_{\ast }M_{2}(Z))^{N/(3N+1)}.
\end{equation*}
\end{theorem}

\textbf{Proof}. We denote $A_{n}=\mathrm{Cov}^{1/2}(Z_{n})$ and we define $%
\overline{Z}_{n}=A_{n}^{-1}\times Z_{n},$ so that $\overline{Z}_{n,1},\ldots
,\overline{Z}_{n,m_{\ast }}$ are orthonormal. In the formalism in \cite%
{[MDO]}, $\overline{Z}_{n}$ is called an ``orthonormal ensemble''. Then we
define
\begin{equation}
\overline{c}((n_{1},k_{1}),\ldots ,(n_{N},k_{N}))=\sum_{i_{1},\ldots
,i_{N}=1}^{m_{\ast }}c((n_{1},i_{1}),\ldots
,(n_{N},i_{N}))A_{n_{1}}^{i_{1},k_{1}}\ldots A_{n_{N}}^{i_{N},k_{N}}
\label{H13}
\end{equation}%
and we notice that, with this definition,
\begin{equation}
S_{N}(c,Z)=S_{N}(\overline{c},\overline{Z}).  \label{H14}
\end{equation}%
Moreover one easily checks that%
\begin{equation}
\left\vert \overline{c}\right\vert \leq (m_{\ast }M_{2})^{N}\left\vert
c\right\vert \quad \mbox{and}\quad \delta _{\ast }(\overline{c})\leq
(m_{\ast }M_{2})^{N}\delta _{\ast }(c).  \label{H16}
\end{equation}%
Let us check that $\overline{Z}$ is hypercontractive in the sense of \cite%
{[MDO]}. We notice that $M_{p}(\overline{Z})\leq \underline{\lambda }%
^{-m_{\ast }}M_{p}(Z)$ and we take $\eta ^{-1}=b_{p}(b_{p}\underline{\lambda
}^{-m_{\ast }}M_{p}(Z))^{N}.$ Then, for any coefficients $c\in \mathcal{C}({%
\mathbb{R}}) $ we have (with $p=3)$%
\begin{eqnarray*}
\left\Vert S_{N}(c,\eta \overline{Z})-c(\emptyset )\right\Vert _{p} &\leq
&b_{p}(b_{p}M_{p}(\overline{Z)})^{N}(\sum_{1\leq \left\vert \alpha
\right\vert \leq N}\eta ^{\left\vert \alpha \right\vert }\left\vert
\overline{c}(\alpha )\right\vert ^{2})^{1/2} \\
&\leq &b_{p}(b_{p}\underline{\lambda }^{-m_{\ast }}M_{p}(Z))^{N}(\sum_{1\leq
\left\vert \alpha \right\vert \leq N}\eta ^{\left\vert \alpha \right\vert
}\left\vert \overline{c}(\alpha )\right\vert ^{2})^{1/2} \\
&\leq &(\sum_{1\leq \left\vert \alpha \right\vert \leq N}\left\vert
\overline{c}(\alpha )\right\vert ^{2})^{1/2}=\left\Vert S_{N}(c,\overline{Z}%
)-c(\emptyset )\right\Vert _{2}
\end{eqnarray*}%
and this means, in the formalism from \cite{[MDO]} that $\overline{Z}$ is $%
(2,3,\eta )-$hypercontractive. Now we are able to use Theorem 3.19 in \cite%
{[MDO]} (which is written in terms of $\tau=\delta^2_\ast(c)$), and this
yields (\ref{H12}). $\square $

\section{Main results}

\label{sect:3}

\subsection{Doeblin's condition and splitting}

\label{sect:3.1}

We fix $d_{\ast }\in {\mathbb{N}}$ and $k_{\ast }\in {\mathbb{N}}$, we
denote $m_{\ast }=d_{\ast }\times k_{\ast },$ and we work with a sequence of
independent random variables $X=(X_{n})_{n\in {\mathbb{N}}},$ $%
X_{n}=(X_{n,1},\ldots ,X_{n,d_{\ast }})\in {\mathbb{R}}^{d_{\ast }}.$ We
deal with general polynomials with variables $X_{n,j}$ that is, with linear
combinations of monomials $\prod_{i=1}^{m}X_{n_{i},j_{i}}^{k_{i}},k_{i}\leq
k_{\ast }.$ Because of the powers $k_{i}$, this is no more a multi-linear
polynomial. In order to come back to multi-linear polynomials we define $%
Z_{n}(X)\in {\mathbb{R}}^{m_{\ast }}$ by
\begin{equation}
Z_{n,kd_{\ast }+j}(X)=X_{n,j}^{k+1}-{\mathbb{E}}(X_{n,j}^{k+1})\quad %
\mbox{for}\quad j\in \lbrack d_{\ast }],k\in \{0,1,\ldots ,k_{\ast }-1\}.
\label{H15}
\end{equation}%
With this definition, if $\alpha =((n_{1},l_{1}),\ldots ,(n_{m},l_{m}))$,
with $n_{1}<\cdots <n_{m}$ and $l_{1},\ldots ,l_{m}\in \{1,\ldots ,m_{\ast
}\}$, then
\begin{equation*}
Z^{\alpha }(X)=\prod_{i=1}^{m}(X_{n_{i},j_{i}}^{k_{i}+1}-{\mathbb{E}}%
(X_{n_{i},j_{i}}^{k_{i}+1}))
\end{equation*}%
where $(k_{i},j_{i})=(k(l_{i}),j(l_{i}))$, $i=1,\ldots ,m$, with
\begin{equation}
k(l)=\Big\lfloor\frac{l-1}{d_{\ast }}\Big\rfloor\quad \mbox{and}\quad j(l)=\Big\{\frac{%
l-1}{d_{\ast }}\Big\}d_{\ast }+1,  \label{kj}
\end{equation}%
in which the symbols $\lfloor x\rfloor$ and $\{x\}$ denote the integer and the fractional
part of $x\geq 0$ respectively. We denote
\begin{equation}
Q_{N,k_{\ast }}(c,X)=\sum_{0\leq \left\vert \alpha \right\vert \leq
N}c(\alpha )Z^{\alpha }(X)=S_{N}(c,Z(X)),  \label{H15'}
\end{equation}%
that is
\begin{align*}
Q_{N,k_{\ast }}(c,X)=\sum_{m=0}^{N}\sum_{n_{1}<\cdots
<n_{m}}\sum_{k_{1},\ldots ,k_{m}=1}^{k_{\ast }}\sum_{j_{1},\ldots
,j_{m}=1}^{d_{\ast }}& c\big((n_{1},(k_{1}-1)d_{\ast }+j_{1}),\ldots
,(n_{m},(k_{m}-1)d_{\ast }+j_{m})\big)\times  \\
& \times \prod_{i=1}^{m}\big(X_{n_{i},j_{i}}^{k_{i}}-{\mathbb{E}}%
(X_{n_{i},j_{i}}^{k_{i}})\big),
\end{align*}%
which agrees with (\ref{I1})-(\ref{I2}) in dimension 1 ($d_{\ast }=1$).

The crucial hypothesis in this section is that for every $n\in {\mathbb{N}},$
the law of $X_{n}$ is locally lower bounded by the Lebesgue measure - this
is Doeblin's condition. Let us be more precise.

\medskip

\textbf{Hypothesis ${\mathfrak{D}}(\varepsilon, r, R)$.} \textit{Let $%
\varepsilon>0$, $r>0$ and $R>0$ be fixed. We say that $X=(X_n)_{n\in{\mathbb{%
N}}}$ satisfies hypothesis ${\mathfrak{D}}(\varepsilon,r,R)$ if there exist $%
x_{n}\in {\mathbb{R}}^{d_{\ast }},n\in {\mathbb{N}}$ such that for every
measurable set $A\subset B_{r}(x_{n})$%
\begin{equation}
{\mathbb{P}}(X_{n}\in A)\geq \varepsilon \lambda (A),  \label{N8}
\end{equation}%
$\lambda $ denoting the Lebesgue measure on ${\mathbb{R}}^{d_{\ast }}$, and
\begin{equation}  \label{LV2}
\sup_{n\in{\mathbb{N}}}\left\vert x_{n}\right\vert \leq R.
\end{equation}
}

\medskip

Note that there is no assumption about $X_{n}$, $n\in{\mathbb{N}}$, being
identically distributed, but the fact that the parameters $\varepsilon $, $r$
and $R$ are the same for every $n,$ represents a uniformity assumption. Note
also that this property never holds for $Z_{n}(X).$ This is why we are
obliged to work with $X_{n}$ only.

\medskip

\textbf{Hypothesis ${\mathfrak{M}}(\varepsilon,r,R)$.} \textit{We say that $%
X=(X_n)_{n\in{\mathbb{N}}}$ satisfies hypothesis ${\mathfrak{M}}%
(\varepsilon,r,R)$ if ${\mathfrak{D}}(\varepsilon, r, R)$ holds and if for
every $p\geq 1$ one has $\sup_{n\in {\mathbb{N}}}\left\Vert X_{n}\right\Vert
_{p}<\infty $.}

\medskip

Note that if Assumption ${\mathfrak{M}}(\varepsilon,r,R)$ holds then $%
Z_{n}(X)$ verifies (\ref{H1}).

The interesting point about random variables which verity ${\mathfrak{D}}%
(\varepsilon ,r,R)$ is that one may use a splitting method in order to
obtain a nice representation for $X_{n}$ (in law). We introduce the
auxiliary functions $\theta _{r},\psi _{r}:{\mathbb{R}}\rightarrow {\mathbb{R%
}}_{+}$ defined by
\begin{equation}
\theta _{r}(t)=1-\frac{1}{1-(\frac{t}{r}-1)^{2}}\qquad \psi
_{r}(t)=1_{\{\left\vert t\right\vert \leq r\}}+1_{\{r<\left\vert
t\right\vert \leq 2r\}}e^{\theta _{r}(\left\vert t\right\vert )}  \label{N4}
\end{equation}%
and we denote%
\begin{equation}
\mathfrak{m}_{r}=\int_{{\mathbb{R}}}\psi _{r}(\left\vert z\right\vert
^{2})dz.  \label{N6}
\end{equation}%
Let $V_{n},U_{n}\in {\mathbb{R}}^{d_{\ast }}$ and $\chi _{n}\in \{0,1\}$ be
independent random variables with laws%
\begin{equation}
\begin{array}{l}
\displaystyle{\mathbb{P}}(\chi _{n}=1)=\varepsilon \mathfrak{m}_{r}^{d_{\ast
}},\quad {\mathbb{P}}(\chi _{n}=0)=1-\varepsilon \mathfrak{m}^{d_{\ast
}}\smallskip \\
\displaystyle{\mathbb{P}}(V_{n}\in dx)=\frac{1}{\mathfrak{m}_{r}^{d_{\ast }}}%
\prod_{k=1}^{d_{\ast }}\psi _{r}(\left\vert x_{k}-x_{n,k}\right\vert
^{2})dx_{1}\ldots dx_{d_{\ast }}\smallskip \\
\displaystyle{\mathbb{P}}(U_{n}\in dx)=\frac{1}{1-\mathfrak{m}_{r}^{d_{\ast
}}}\Big({\mathbb{P}}(X_{n}\in dx)-\varepsilon \prod_{k=1}^{d_{\ast }}\psi
_{r}(\left\vert x_{k}-x_{n,k})\right\vert ^{2}\Big)dx_{1}\ldots dx_{d_{\ast
}}.%
\end{array}
\label{laws}
\end{equation}%
Note that the hypothesis ${\mathfrak{D}}(\varepsilon ,r,R)$ ensures that ${%
\mathbb{P}}(X_{n}\in dx)-\varepsilon \prod_{k=1}^{d_{\ast }}\psi
_{r}(\left\vert x_{k}-x_{n,k})\right\vert ^{2})dx\geq 0,$ so that the law of
$U_{n}$ is well defined. It is easy to check that $\chi _{n}V_{n}+(1-\chi
_{n})U_{n}$ has the same law as $X_{n}$. Since all our statements concern
only the law of $X_{n}$, now on we assume that%
\begin{equation}
X_{n}=\chi _{n}V_{n}+(1-\chi _{n})U_{n}.  \label{repr}
\end{equation}%
Let us mention a nice property for the function $\psi _{r}$: it is easy to
check that for each $k\in {\mathbb{N}},p\geq 1$ there exists a universal
constant $C_{k,p}\geq 1$ such that
\begin{equation}
\psi _{r}(t)|\theta _{r}^{(k)}(|t|)|^{p}\leq \frac{C_{k,p}}{r^{kp}}
\label{N5}
\end{equation}%
where $\theta _{r}^{(k)}$ denotes the derivative of order $k$ of $\theta
_{r}.$

Actually, the uniformity property (\ref{LV2}) has not been used so far. We
see now that it gives a ``non degeneracy'' for the powers of the components
of $V_{n}$ uniformly in $n\in{\mathbb{N}}$. More precisely, we define the
random vector $\widetilde V_n=Z_n(V)$ in ${\mathbb{R}}^{m_*}$, that is
\begin{equation}
\widetilde{V}_{n,l}=V_{n,j(l)}^{k(l)+1}-{\mathbb{E}}(V_{n,j(l)}^{k(l)+1}),%
\quad l=1,\ldots,m_\ast,  \label{LV1}
\end{equation}
where $k(l)$ and $j(l)$ are given in (\ref{kj}). Then, one has the following
result.

\begin{lemma}
\label{lambdaR} Let $R>0$ be such that (\ref{LV2}) holds and let $\mathrm{Cov%
}(\widetilde V_n)$ denote the covariance matrix of $\widetilde V_n$.
Then there exists $\lambda _{R}>0$ such that
\begin{equation}
\langle \mathrm{Cov}(\widetilde V_n)\xi ,\xi \rangle \geq \lambda
_{R}\left\vert \xi \right\vert ^{2}  \label{LV3}
\end{equation}
for every $\xi \in {\mathbb{R}}^{m_{\ast }}$ and $n\in {\mathbb{N}}$.
\end{lemma}

\textbf{Proof}. For $y\in {\mathbb{R}}^{d_{\ast }}$ and $\xi \in {\mathbb{R}}%
^{m_{\ast }}$ we define%
\begin{eqnarray*}
e_{l}(y) &=&\frac{1}{\mathfrak{m}_{r}^{d_{\ast }}}\int
x_{j(l)}^{k(l)}\prod_{i=1}^{d_{\ast }}\psi _{r}(\left\vert
x_{i}-y_{i}\right\vert ^{2})dx,\quad l\in \lbrack m_{\ast }],\quad \mbox{and}
\\
I_{\xi }(y) &=&\frac{1}{\mathfrak{m}_{r}^{d\ast }}\int \Big(%
\sum_{l=1}^{m_{\ast }}(x_{j(l)}^{k(l)}-e_{l}(y))\xi _{l}\Big)%
^{2}\prod_{i=1}^{d_{\ast }}\psi _{r}(\left\vert x_{i}-y_{i}\right\vert
^{2})dx
\end{eqnarray*}%
If $I_{\xi }(y)=0$ then $\sum_{l=1}^{m_{\ast }}(x_{j(l)}^{k(l)}-e_{l}(y))\xi
_{l}=0$ for $x$ in an open set, and this imply that $\xi =0.$ Since $\xi
\mapsto I_{\xi }(y)$ is continuous, it follows that $\lambda
(y)=\inf_{\left\vert \xi \right\vert =1}I_{\xi }(y)>0.$ And since $y\mapsto
\lambda (y)$ is continuous it follows that one may find $\lambda _{R}>0$
such that $\inf_{\left\vert y\right\vert \leq R}\lambda (y)\geq \lambda _{R}$%
. Now, we note that $e_{l}(x_{n})={\mathbb{E}}(V_{n,j(l)}^{k(l)})={\mathbb{E}%
}(\widetilde{V}_{n,l})$ and $I_{\xi }(x_{n})=<\mathrm{Cov}(\widetilde{V}%
_{n})\xi ,\xi \>$. Thus, if $|\xi |=1$ we get $\inf_{n}<\mathrm{Cov}(%
\widetilde{V}_{n})\xi ,\xi \>=\inf_{n}\inf_{|\xi |=1}I_{\xi }(x_{n})\geq
\lambda _{R}$, and (\ref{LV3}) follows. $\square $

\medskip

We conclude with an inequality which will be useful later on.

\begin{lemma}
Let $R>0$ be such that (\ref{LV2}) holds and let $\lambda _{R}$ be given in
Lemma \ref{lambdaR}. Then for every $d\in \mathcal{C}({\mathbb{R}})$,
\begin{equation}
{\mathbb{E}}(|S_{N}(d,\widetilde{V})|^{2})\geq \lambda
_{R}^{N}\sum_{m=0}^{N}|d|_{m}^{2}=\lambda _{R}^{N}|d|^{2},  \label{LV4}
\end{equation}%
with $\widetilde{V}=Z(V)$ defined in (\ref{LV1}).
\end{lemma}

\textbf{Proof.} We first fix an integer $m$, $n_{1}<\cdots <n_{m}$ and we
consider $d(l_{1},\ldots ,l_{m})$, $l_{i}\in \lbrack m_{\ast }]$. We prove
that
\begin{equation}
{\mathbb{E}}\Big(\Big(\sum_{l_{1},...,l_{m}=1}^{m_{\ast
}}d(l_{1},...,l_{m})\prod_{i=1}^{m}\widetilde{V}_{n_{i},l_{i}}\Big)^{2}\Big)%
\geq \lambda _{R}^{m}\sum_{l_{1},...,l_{m}=1}^{m_{\ast
}}d^{2}(l_{1},...,l_{m}).  \label{LV4-bis}
\end{equation}%
We define the random variable
\begin{equation*}
\widehat{d}(l_{m})=\sum_{l_{1},\ldots ,l_{m-1}=1}^{m_{\ast }}d(l_{1},\ldots
,l_{m})\prod_{i=1}^{m-1}\widetilde{V}_{n_{i},l_{i}}.
\end{equation*}%
We notice that $\widehat{d}(k),k\in \lbrack m_{\ast }]$ are independent of $%
\widetilde{V}_{n_{m},l},l\in \lbrack m_{\ast }]$ and that
\begin{equation*}
\sum_{l_{1},...,l_{m}=1}^{m_{\ast }}d(l_{1},...,l_{m})\prod_{i=1}^{m}%
\widetilde{V}_{n_{i},l_{i}}=\sum_{l_{m}=1}^{m_{\ast }}\widehat{d}(l_{m})%
\widetilde{V}_{n_{m},l_{m}}.
\end{equation*}%
So,
\begin{align*}
& {\mathbb{E}}\Big(\Big(\sum_{l_{1},\ldots ,l_{m}=1}^{m_{\ast
}}d(l_{1},\ldots ,l_{m})\prod_{i=1}^{m}\widetilde{V}_{n_{i},l_{i}}\Big)^{2}%
\Big)={\mathbb{E}}\Big(\sum_{l_{m},\bar{l}_{m}=1}^{m_{\ast }}\widehat{d}%
(l_{m})\widehat{d}(\bar{l}_{m}){\mathbb{E}}(\widetilde{V}_{n_{m},l_{m}}%
\widetilde{V}_{n_{m},\bar{l}_{m}})\Big) \\
& \geq \lambda _{R}{\mathbb{E}}\Big(\sum_{l_{m}=1}^{m_{\ast }}\widehat{d}%
(l_{m})^{2}\Big)=\lambda _{R}\sum_{l_{m}=1}^{m_{\ast }}{\mathbb{E}}\Big(\Big(%
\sum_{l_{1},\ldots ,l_{m-1}=1}^{m_{\ast }}d(l_{1},\ldots
,l_{m-1},l_{m})\prod_{i=1}^{m-1}\widetilde{V}_{n_{i},l_{i}}\Big)^{2}\Big),
\end{align*}%
the above lower bound following from (\ref{LV3}). By iteration, one gets (%
\ref{LV4}).

Consider now the general case. We recall that, for any two multi-indexes $%
\alpha$ and $\overline{\alpha}$, ${\mathbb{E}}(\widetilde V^\alpha\widetilde
V^{\overline{\alpha}})\neq 0$ if and only if $\alpha^{\prime }=\overline{%
\alpha}^{\prime }$. This gives
\begin{align*}
&{\mathbb{E}}(|S_N(d,\widetilde V|^2) =\sum_{m=0}^N\sum_{|\alpha|=|\overline{%
\alpha}|=m, \alpha^{\prime }=\overline{\alpha}^{\prime }} d(\alpha)d(%
\overline{\alpha}){\mathbb{E}}(\widetilde V^\alpha\widetilde V^{\overline{%
\alpha}}) \\
&=\sum_{m=0}^N\sum_{n_1<\cdots<n_m} {\mathbb{E}}\Big(\Big(%
\sum_{l_1,\ldots,l_m\in[m_*]}d_{n_1,\ldots,n_m}(l_1,\ldots,l_m)\prod_{i=1}^m%
\widetilde V_{n_i,l_i}\Big)^2\Big)
\end{align*}
where, for fixed $n_1<\ldots<n_m$, we have set $d_{n_1,\ldots,n_m}(l_1,%
\ldots,l_m)=d((n_1,l_1),\ldots,(n_m,l_m))$. The statement now follows from (%
\ref{LV4-bis}). $\square $

\subsection{Main results}

\label{sect:3.2}

Our goal is to estimate the total variation distance between two polynomials
of type $Q_{N,k_{\ast }}(c,X)$, which we write as in (\ref{H15'}), that is
\begin{equation*}
Q_{N,k_*}(c,X)=\sum_{0\leq |\alpha|\leq N}c(\alpha)Z^\alpha(X),
\end{equation*}
where $Z(X)$ is defined in (\ref{H15}) and $\alpha=(\alpha^{\prime
},\alpha^{\prime \prime })$ with $\alpha^{\prime \prime }_i\in[m_*]$, $%
m_*=d_*k_*$.

We will use the following quantities related to the coefficients $c.$ We
work first with the Hilbert space $\mathcal{U}={\mathbb{R}}$ (so, we drop $%
\mathcal{U}$ from the notation) and we recall that $\left\vert c\right\vert
=\left\vert c\right\vert _{\mathcal{U}}$ is defined in (\ref{H2a}) and $%
\delta _{\ast }(c)=\delta _{{\mathcal{U}}, \ast }(c)$ is defined in (\ref%
{H2b}). Moreover, for $m\leq N,$ we define%
\begin{equation*}
\left\vert c\right\vert _{m}=\Big(\sum_{\left\vert \alpha \right\vert
=m}c^{2}(\alpha )\Big)^{1/2}\quad \mbox{and}\quad \left\vert c\right\vert
_{m,N}=\Big(\sum_{m\leq \left\vert \alpha \right\vert \leq N}c^{2}(\alpha )%
\Big)^{1/2}.
\end{equation*}%
Finally we assume that $X$ verifies ${\mathfrak{D}}%
(\varepsilon, r, R)$ and we denote
\begin{equation}  \label{e}
e_{m,N}(c) =\exp \Big(-\Big(\frac{\varepsilon\mathfrak{m}_r}{2}\Big)^{2m}%
\frac{|c|_{m}^{2}}{\delta _{\ast}^{2}(c)}\Big).
\end{equation}
Notice that if $X$ and $Y$ satisfy ${\mathfrak{D}}(\varepsilon, r, R)$
respectively ${\mathfrak{D}}(\varepsilon ^{\prime },r^{\prime }, R^{\prime
}),$ then they both satisfy ${\mathfrak{D}}(\varepsilon \wedge \varepsilon
^{\prime },r\wedge r^{\prime }, R\vee R^{\prime })$ so we may assume that $%
\varepsilon $, $r$ and $R$ are the same.

For $k\in {\mathbb{N}}$ we define the distances%
\begin{equation*}
d_{k}(F,G)=\sup \{\left\vert {\mathbb{E}}(f(F))-{\mathbb{E}}%
(f(G))\right\vert :\left\Vert f\right\Vert _{k,\infty }\leq 1\}.
\end{equation*}%
Note that $d_{0}=d_{\mbox{\rm{\scriptsize{TV}}}}$ is the total variation
distance and $d_{1}$ is the Fortet Mourier distance (which metrizes the
convergence in law). We give now our first result:

\begin{theorem}
\label{D} Suppose that $X$ and $Y$ verify Hypothesis ${\mathfrak{M}}%
(\varepsilon ,r,R)$ (that is (\ref{H1}) and ${\mathfrak{D}}(\varepsilon
,r,R) $) and let $c,d\in \mathcal{C}({\mathbb{R}})$ be two families of
coefficients. We fix $k,k_{\ast }$ and $N$ and $m\leq N$ and $m^{\prime
}\leq N$ such that $|c|_{m}>0$ and $|d|_{m^{\prime }}>0$ and we denote $%
\overline{m}=m\vee m^{\prime }.$ We also assume that%
\begin{equation}
d_{k}:=d_{k}(Q_{N,k_{\ast }}(c,X),Q_{N,k_{\ast }}(d,Y))\vee (\left\vert
c\right\vert _{m+1,N}^{2}+\left\vert d\right\vert _{m^{\prime }+1,N}^{2})^{%
\frac{2kk_{\ast }\overline{m}+1}{k_{\ast }\overline{m}}}\leq 1  \label{D3''}
\end{equation}%
Let $\theta \in ((\frac{1}{1+k})^{2},1)$. Then there exist $C>0$ and $a\in (%
\frac{1}{1+k},1]$, which depend on the parameters $\varepsilon
,r,R,k,k_{\ast },N,m,m^{\prime },\theta $ and the moment bounds $M_{p}(X)$, $%
M_{p}(Y)$ for a suitable $p>1,$ but independent of the coefficients $c,d\in
\mathcal{C}({\mathbb{R}})$, such that
\begin{equation}
\begin{array}{ll}
& \displaystyle\left\vert {\mathbb{E}}(f(Q_{N,k_{\ast }}(c,X)))-{\mathbb{E}}%
(f(Q_{N,k_{\ast }}(d,Y))\right\vert \smallskip \\
& \leq \displaystyle C\max \Big(1,\Big(|c|_{m}^{-\frac{2}{k_{\ast }m}%
}+|d|_{m^{\prime }}^{-\frac{2}{k_{\ast }m^{\prime }}}\Big)^{a}\Big)%
\left\Vert f\right\Vert _{\infty }(1+|c|+|d|)^{5k}\times \smallskip \\
& \displaystyle\times \Big(e_{m,N}^{a}(c)+e_{m^{\prime
},N}^{a}(d)+d_{k}(Q_{N,k_{\ast }}(c,X),Q_{N,k_{\ast }}(d,Y))^{\frac{\theta }{%
1+2kk_{\ast }\overline{m}}}+\left\vert c\right\vert _{m+1,N}^{\frac{2\theta
}{k_{\ast }\overline{m}}}+\left\vert d\right\vert _{m^{\prime }+1,N}^{\frac{%
2\theta }{k_{\ast }\overline{m}}}\Big),%
\end{array}
\label{D3}
\end{equation}%
$e_{m,N}(c)$ and $e_{m^{\prime },N}(d)$ being defined in (\ref{e}).
\end{theorem}

In practical situations, one has $|c|^2 _{m+1,N}=|d|^2_{m^{\prime }+1,N}=0$
or both $|c|^2 _{m+1,N}$ and $|d|^2_{m^{\prime }+1,N}$ are very small, so $%
d_k$ in (\ref{D3''}) is actually the $d_k$-distance between $Q_{N,k_{\ast
}}(c,X)$ and $Q_{N,k_{\ast }}(d,Y)$.

\smallskip

The proof of Theorem \ref{D} is done by using a Malliavin type calculus based
on $V_{n},n\in {\mathbb{N}}$ which we present in Section \ref{sect:doeblin},
so we postpone it for Section \ref{sect:proof}. It represents the main
effort in our paper.

As an immediate consequence, we give the following estimate of the total
variation distance between two multiple stochastic integrals. We consider a $%
m_{\ast }$ dimensional Brownian motion $W=(W^{1},\ldots ,W^{m_{\ast }}),$ we
fix $\kappa =(k_{1},\ldots ,k_{m})\in \lbrack m_{\ast }]^{m},$ and, for a
symmetric kernel $f\in L^{2}[0,1]^{m},$ we denote%
\begin{equation*}
I_{\kappa
}(f)=m!\int_{0}^{1}dW_{s_{m}}^{k_{m}}\int_{0}^{s_{m}}dW_{s_{m-1}}^{k_{m-1}}%
\cdots \int_{0}^{s_{2}}f(s_{1},\ldots,s_{m})dW_{s_{1}}^{k_{1}}.
\end{equation*}

\begin{theorem}
\label{mul} Let $f,g\in L^{2p}[0,1]^{m},p>1.$ Then, for every $k,m\in {%
\mathbb{N}}_{\ast }$ and $\theta \in ((\frac{1}{1+k})^{2},1)$ there exist $%
C>0$ and $a\in (\frac{1}{1+k},1)$ (both depending on $\theta ,m$ and $k$)
such that%
\begin{equation}
\begin{array}{ll}
& \displaystyle d_{{\mbox{\rm{\scriptsize{TV}}}}}(I_{\kappa }(f),I_{\kappa
}(g))\smallskip \\
& \leq \displaystyle C(m!)^{5k/2}\max \Big(1,\big(\Vert f\Vert _{2}^{-{2}/{m}%
}+\Vert g\Vert _{2}^{-{2}/{m}}\big)^{a}\Big)(1+\Vert f\Vert _{2}+\Vert
g\Vert _{2})^{5k}\,d_{k}^{\theta /(1+2km)}(I_{\kappa }(f),I_{\kappa }(g)).%
\end{array}
\label{D5}
\end{equation}%
%
%
%
%
\end{theorem}

\begin{remark}
In the case $k=1,$ the above result has first been announced in \cite{[D]}
with the power $\frac{1}{m}$ instead of $\frac{\theta }{2m+1}$ above, but
the proof was only sketched. It has rigourously been proved in \cite{[NPy]}
with power $\frac{1}{2m+1}$ and recently improved in \cite{[Bo]} where the
power $\frac{1}{m}\times (\ln m)^{d}$ is obtained. So (\ref{D5}) is not the
best possible estimate. This also indicates that the power in (\ref{D3}) is
not optimal (but the approach in \cite{[Bo]} does not seem to work in our
general framework, so for the moment we are not able to improve it).
\end{remark}

\begin{remark}
Theorem \ref{mul}, with exactly the same proof, extends to general random
variables which live in a finite sum of Wiener chaoses: let $F$ and $G$ be
two random variables belonging to $\oplus _{m=0}^{N}\mathcal{W}_{m}$ where $%
\mathcal{W}_{m}$ is the chaos of order $m.$ We denote by $P_{m}$ the
projection on $\mathcal{W}_{m}$ and we put $m(F)=\max \{m:P_{m}F\neq 0\}$
and $\alpha (F)=\Vert P_{m(F)}F\Vert _{2}^{-2/{m(F)}}.$ Then, with $%
N=m(F)\vee m(G),$%
\begin{equation}
d_{{\mbox{\rm{\scriptsize{TV}}}}}(F,G)\leq C\max \Big(1,\big(\alpha
(F)+\alpha (G)\big)^{a}\Big)(1+\Vert F\Vert _{2}+\Vert G\Vert
_{2})^{5k}\,d_{k}^{\theta /(1+2kN)}(F,G),  \label{D5'}
\end{equation}%
where $a\in (\frac{1}{1+k},1)$ and $C>0$ depend on $\theta ,k,N$.
\end{remark}

\noindent \textbf{Proof of Theorem \ref{mul}.} Let $n\in \N.$ For $\alpha
^{\prime }=(\alpha _{1}^{\prime },\ldots ,\alpha _{m}^{\prime })\in [
n-1]^{m}$, we denote $I_{\alpha ^{\prime }}=\prod_{j=1}^{m}[\frac{\alpha
_{j}^{\prime }}{n},\frac{\alpha _{j}^{\prime }+1}{n})$ and we define
\begin{equation*}
f_{n}(s)=\sum_{\alpha ^{\prime }}d_{n,f}(\alpha ^{\prime })1_{I_{\alpha
^{\prime }}}(s)\quad \mbox{with}\quad d_{n,f}(\alpha ^{\prime
m}\int_{I_{\alpha ^{\prime }}}f(u)du.
\end{equation*}%
Note that $f_{n}$ is the conditional expectation of $f$ with respect to the
partition $I_{\alpha ^{\prime }}$ and to the uniform law on $[0,1]^{m}.$
Take now $\alpha =(\alpha _{1},\ldots ,\alpha _{n})$ with $\alpha
_{i}=(\alpha _{i}^{\prime },\alpha _{i}^{\prime \prime })$ and $(\alpha
_{1}^{\prime \prime },\ldots ,\alpha _{m}^{\prime \prime })\in \lbrack
m_{\ast }]^{m}$. We denote%
\begin{equation*}
c_{n,f}(\alpha )=m!\,n^{-m/2}\,d_{n,f}(\alpha ^{\prime })1_{\alpha
_{1}^{\prime }<\cdots <\alpha _{m}^{\prime }<n}\prod_{i=1}^{m}1_{\alpha
_{i}^{\prime \prime }=k_{i}},\quad G_{\alpha _{i}^{\prime },\alpha
_{i}^{\prime \prime }}=n^{1/2}\times \Big(W^{\alpha _{i}^{\prime \prime }}%
\Big(\frac{\alpha _{i}^{\prime }+1}{n}\Big)-W^{\alpha _{i}^{\prime \prime }}%
\Big(\frac{\alpha _{i}^{\prime }}{n}\Big)\Big).
\end{equation*}%
so that%
\begin{equation*}
I_{\kappa }(f_{n})=\sum_{\alpha }c_{n,f}(\alpha )G^{\alpha }=\Phi
_{m}(c_{n,f},G).
\end{equation*}%
We are now in the framework of Theorem \ref{D} and we compare $\Phi
_{m}(c_{n,f},G)$ and $\Phi _{m}(c_{n,g},G)$. We take $k_{\ast }=1,d_{\ast
}=m_{\ast }$ and $N=m=m^{\prime }.$ Then $\left\vert c_{n,f}\right\vert
_{m+1,N}=\left\vert c_{n,g}\right\vert _{m+1,N}=0$.
Let us estimate the parameters associated to $c_{n,f}.$ By the convergence
theorem for martingales $\left\vert c_{n,f}\right\vert _{m}^{2}=m!\left\Vert
f_{n}\right\Vert _{2}^{2}\rightarrow m!\left\Vert f\right\Vert _{2}^{2}>0.$
We estimate now $\delta _{\ast }(c_{n,f})$. By using H\"{o}lder's
inequality,
\begin{align*}
\delta _{\ast }^{2}(c_{n,f})& =\max_{i\in \lbrack
n]}\sum_{j=1}^{m}\sum_{\alpha ^{\prime }\,:\,\alpha _{j}^{\prime
}=i}(m!)^{2}n^{-m}\Big(n^{m}\int_{I_{\alpha ^{\prime }}}f(s)ds\Big)^{2} \\
& =\max_{i\in \lbrack n]}(m!)^{2}\sum_{j=1}^{m}\sum_{\alpha ^{\prime }}n^{-m}%
\Big(n^{m}\int_{I_{\alpha ^{\prime }}}f(s)1_{s_{j}\in \lbrack \frac{i}{n},%
\frac{i+1}{n})}ds\Big)^{2} \\
& \leq \max_{i\in \lbrack n]}(m!)^{2}\sum_{j=1}^{m}\sum_{\alpha ^{\prime
}}\int_{I_{\alpha ^{\prime }}}f^{2}(s)1_{s_{j}\in \lbrack \frac{i}{n},\frac{%
i+1}{n})}ds\leq \max_{i\in \lbrack n]}m!\max_{j\in \lbrack
m]}\int_{[0,1]^{m}}f^{2}(s)1_{s_{j}\in \lbrack \frac{i}{n},\frac{i+1}{n})}ds
\\
& \leq m!\Vert f\Vert _{2p}^{2}\frac{1}{n^{1-1/p}}\rightarrow 0
\end{align*}%
so that $e_{m,m}(c_{n,f})\rightarrow 0$ and $e_{m,m}(c_{n,g})\rightarrow 0$
as $n\rightarrow \infty $.

Now (\ref{D3}) gives, for $\theta <1,$ and $n,n^{\prime }\in {\mathbb{N}}$%
\begin{equation}  \label{D5a}
\begin{array}{ll}
& \displaystyle d_{\mbox{\rm{\scriptsize{TV}}}}(I_{\kappa }(f_{n}),I_{\kappa
}(g_{n^{\prime }}))\smallskip \\
& \leq \displaystyle C(m!)^{5k/2}\max \Big(1,\Big(\|f_n\|_2^{-\frac{2}{m}%
}+\|g_{n^{\prime }}\|_2^{-\frac{2}{m}}\Big)^a\Big) (1+\|f_n\|_2+\|g_{n^{%
\prime }}\|_2)^{5k}\times\smallskip \\
& \displaystyle \times\Big(e_{m,m}^{a}(c_{n,f})+e_{m,m}^{a}(c_{n^{\prime
},g})+d_k^{\theta/(1+2km)}(I_{\kappa }(f_{n}),I_{\kappa }(g_{n^{\prime }}))%
\Big),%
\end{array}%
\end{equation}%
where $a\in(\frac 1{1+k},1)$. We take $n^{\prime }>n$ and we notice that $%
d_{k}(I_{\kappa }(f_{n}),I_{\kappa }(f_{n^{\prime }}))\leq \left\Vert
f_{n}-f_{n^{\prime }}\right\Vert _{2}\rightarrow 0$ so that the above
inequality gives $d_{\mbox{\rm{\scriptsize{TV}}}}(I_{\kappa
}(f_{n}),I_{\kappa }(f_{n^{\prime }}))\rightarrow 0$ as $n,n^{\prime
}\rightarrow \infty .$ It follows that the sequences $I_{\kappa }(f_{n})$
and $I_{\kappa }(g_{n}),n\in {\mathbb{N}}$ are Cauchy in $d_{%
\mbox{\rm{\scriptsize{TV}}}}$ and we may pass to the limit in (\ref{D5a}) in
order to obtain (\ref{D5}). $\square $

\medskip

We give now the analogous of Theorem \ref{D} but in terms of Kolmogorov
distance. Here one needs no more Doeblin's condition nor non degeneracy
conditions.

\begin{theorem}
\label{K} Suppose that $X$ and $Y$ verify (\ref{H1}) and are such that $Z(X)$
and $Z(Y)$ both satisfy (\ref{H2}). Let $c,d\in \mathcal{C}({\mathbb{R}})$
be two families of coefficients such that $\left\vert c\right\vert _{N}>0$
and $\left\vert d\right\vert _{N}>0.$
with $\delta_\ast(c),\delta_\ast
(d)\leq 1$. Then, for every $k\in {\mathbb{N}}$ and $\theta \in ((\frac{1}{%
1+k})^{2},1)$ there exist $C>0$ and $a\in (\frac{1}{1+k},1)$ such that
\begin{equation}
\begin{array}{l}
\displaystyle d_{{\mbox{\rm{\scriptsize{Kol}}}}}(Q_{N,k_{\ast
}}(c,X),Q_{N,k_{\ast }}(d,Y))\leq
C(1+|c|_{N}^{-2N}+|d|_{N}^{-2N})(1+|c|+|d|)^{5(k\vee 3)+1}\times \smallskip
\\
\displaystyle\quad \times (\delta _{\ast }^{\theta /(2(k\vee
3)N+1)}(c)+\delta _{\ast }^{\theta /(2(k\vee 3)N+1)}(d)+d_{k\vee 3}^{\theta
/(2(k\vee 3)N+1)}(Q_{N,k_{\ast }}(c,X),Q_{N,k_{\ast }}(d,Y)).%
\end{array}
\label{D5''}
\end{equation}%
where $C>0$ denotes a constant depending on $N$, suitable moments of $X$ and
$Y$ and on the lower bounds $\underline{\lambda }$ in (\ref{H2}) applied to $%
Z(X)$ and $Z(Y)$.
\end{theorem}

\begin{remark}
Note that the estimate (\ref{D5''}) is in terms of $\delta_\ast^{\theta /(2(k\vee 3)N+1)}(c)$ whereas
in (\ref{D3}) it appears $e_{m,N}(c)=\exp (-C\times \frac{|c|_{m}^{2}}{%
\delta _{\ast }^{2}(c)})$ which is much smaller. But we need that $X_{n}$
and $Y_{n}$ satisfy Doeblin's condition ${\mathfrak{D}}(\varepsilon, r, R).$
\end{remark}

\textbf{Proof}. We consider the Gaussian random variables $G_{X}$ and $G_{Y}$
corresponding to $Z(X)$ and $Z(Y)$ respectively and we use Theorem \ref{MDO}
(see (\ref{H12})) in order to obtain
\begin{equation*}
d_{{\mbox{\rm{\scriptsize{Kol}}}}}(Q_{N,k_{\ast }}(c,X),Q_{N,k_{\ast
}}(d,Y))\leq C(\delta _{\ast }^{1/(1+3N)}(c)+\delta _{\ast
}^{1/(1+3N)}(d))+d_{{\mbox{\rm{\scriptsize{Kol}}}}%
}(S_{N}(c,G_{X}),S_{N}(d,G_{Y})).
\end{equation*}%
Using the same argument as in the proof of Theorem \ref{MDO} we may assume
that $G_{X}$ and $G_{Y}$ are standard Gaussian random variables so that $%
S_{N}(c,G_{X})$ and $S_{N}(d,G_{Y})$ are multiple stochastic integrals. By $%
d_{{\mbox{\rm{\scriptsize{Kol}}}}}\leq d_{{\mbox{\rm{\scriptsize{TV}}}}}$
and by (\ref{D5'}) first and (\ref{H11}) (recall that $Q_{N,k_{\ast
}}(c,X)=S_{N}(c,Z_{n}(X))$ then
\begin{align*}
& d_{{\mbox{\rm{\scriptsize{Kol}}}}}(S_{N}(c,G_{X}),S_{N}(d,G_{Y}))\leq d_{{%
\mbox{\rm{\scriptsize{TV}}}}}(S_{N}(c,G_{X}),S_{N}(d,G_{Y})) \\
& \leq C(1+|c|_{N}^{-2N}+|d|_{N}^{-2N})(1+|c|+|d|)^{5(k\vee 3)}d_{k\vee
3}^{\theta /(2(k\vee 3)N+1)}(S_{N}(c,G_{X}),S_{N}(d,G_{Y})) \\
& \leq C(1+|c|_{N}^{-2N}+|d|_{N}^{-2N})(1+|c|+|d|)^{5(k\vee 3)+1}\times \\
& \quad \times (\delta _{\ast }^{\theta /(2(k\vee 3)N+1)}(c)+\delta _{\ast
}^{\theta /(2(k\vee 3)N+1)}(d)+d_{k\vee 3}^{\theta /(2(k\vee
3)N+1)}(Q_{N,k_{\ast }}(c,X),Q_{N,k_{\ast }}(d,Y)).
\end{align*}%
$\square $

\medskip

We give now the invariance principle:

\begin{theorem}
\label{Main}Let $X=(X_{n})_{n\in {\mathbb{N}}}$ be a sequence of independent
${\mathbb{R}}^{d_{\ast }}$ valued random variables which verify Hypothesis ${%
\mathfrak{M}}(\varepsilon,r,R)$
and $G_X=(G_{n,X})_{n\in {\mathbb{N}}},G_{n,X}\in {\mathbb{R}}^{m_{\ast }}$
a sequence of independent and centred Gaussian random variables such that $%
\mathrm{Cov}(G_{n,X})=\mathrm{Cov}(Z_{n}(X)).$ Suppose that for some $m\leq
N $ one has $|c|_{m}>0.$ Let $\theta\in(\frac 1{16},1)$. Then there exist $%
C>0$ and $a\in(\frac 1{4},1]$, which depend on the parameters $\varepsilon,
r, R, k_\ast,N,m,m^{\prime },\theta$ and the moment bounds $M_p(X)$, $M_p(Y)$
for a suitable $p>1$ but independent of the coefficients $c\in\mathcal{C}({%
\mathbb{R}})$, such that
\begin{equation}  \label{MR1}
\begin{array}{rl}
d_{\mbox{\rm{\scriptsize{TV}}}}(Q_{N,k_{\ast }}(c,X),S_{N}(c,G_X))\leq & C
\max(1, |c|_m^{-\frac 2{k_\ast m}})^a (1+|c|)^{19/2}\smallskip \\
& \times \big(\delta _{\ast }^{\frac{\theta}{6k_{\ast }m+1}%
}(c)+e_{m,N}(c)^a+\left\vert c\right\vert _{m+1,N}^{\frac{2\theta }{k_{\ast
}m}}\big).%
\end{array}%
\end{equation}
\end{theorem}

\textbf{Proof.} This is an immediate consequence of Theorem \ref{D} and of
Theorem \ref{smooth}. $\square $

\medskip

In a number of concrete applications (see Theorem \ref{UU} for example), one
takes $S_{N}(c,Z(X))=\sum_{n=m}^{N}\Phi _{n}(c,Z(X))$ and, asymptotically, $%
\Phi _{m}(c,Z(X))$ represents the principal term. Having in mind this we
give the following corollary:

\begin{theorem}
\label{P} Let $c\in \mathcal{C}({\mathbb{R}})$ be such that $c(\alpha )=0$
for $\left\vert \alpha \right\vert \leq m-1$ and $\left\vert c\right\vert
_{m}>0.$ Suppose $|c|_{m+1,N}\leq 1$.

\medskip

\textbf{A}. If $G=(G_n)_{n\in{\mathbb{N}}}$ denote independent centred
Gaussian random variables then, for every $\theta \in(\frac 14,1)$ there
exists $a\in(\frac 12,1]$ such that
\begin{equation}
d_{\mbox{\rm{\scriptsize{TV}}}}(S_{N}(c,G),\Phi _{m}(c,G))\leq
C\max(1,|c|_m^{-\frac 2m})^a(1+|c|)^5 \big(\left\vert c\right\vert
_{m+1,N}^{\frac \theta {2 m+1}}+e_{m,N}(c)^a\big).  \label{MR2}
\end{equation}

\medskip

\textbf{B}. Let $X$ satisfy ${\mathfrak{M}}(\varepsilon,r,R)$
and let $G=(G_{n})_{n\in {\mathbb{N}}},G_{n}\in {\mathbb{R}}^{m_{\ast }},$
be a sequence of independent and centred Gaussian random variables such that
$\mathrm{Cov}(G_{n})=\mathrm{Cov}(Z_{n}(X))$. Then for every $\theta
\in(\frac 14,1)$ there exists $a\in(\frac 12,1]$ such that
\begin{equation}
d_{\mbox{\rm{\scriptsize{TV}}}}(Q_{N,k_{\ast }}(c,X)),\Phi _{m}(c,G)))\leq
C\max(1,|c|_m^{-\frac 2{k_\ast m}})^a(1+|c|)^{\frac {19}2} \big(\delta
_{\ast }^{\frac \theta {6k_{\ast }m+1}}(c)+e_{m,N}(c)^a+\left\vert
c\right\vert _{m+1,N}^{\frac{2\theta}{k_\ast m}\wedge \frac \theta{2m+1}}%
\big).  \label{MR5}
\end{equation}

\textbf{C}. If $Z(X)$ satisfies 
(\ref{H2}) then for every $\theta \in(\frac 14,1)$ there exists $a\in(\frac
12,1]$ such that
\begin{equation}
d_{\mbox{\rm{\scriptsize{Kol}}}}(Q_{N,k_{\ast }}(c,X)),\Phi _{m}(c,G)))\leq
C\max(1,|c|_m^{-\frac 2m})^a(1+|c|)^5 \big(\delta _{\ast }^{\frac 1{1+3N}}+
\left\vert c\right\vert _{m+1,N}^{\frac \theta {2 m+1}}+e_{m,N}(c)^a\big).
\label{MR4}
\end{equation}

\noindent In the above estimates (\ref{MR2}), (\ref{MR5}) and (\ref{MR4}), $%
C>0$ denotes a constant independent of the coefficients $c\in\mathcal{C}({%
\mathbb{R}})$.
\end{theorem}

\textbf{Proof}. One has
\begin{equation*}
d_{1}(S_{N}(c,G),\Phi _{m}(c,G))\leq \left\Vert S_{N}(c,G)-\Phi
_{m}(c,G)\right\Vert _{2}\leq \left\vert c\right\vert _{m+1,N}
\end{equation*}%
so (\ref{MR2}) follows from Theorem \ref{D} (see (\ref{D3})). Using (\ref%
{MR2}) and (\ref{MR1}) we obtain (\ref{MR5}). And (\ref{MR4}) follows from (%
\ref{MR2}) and (\ref{H12}). $\square $

\subsection{Gaussian and Gamma approximation}

\label{sect:3.3}

Theorem \ref{P} has the following interesting application: if one considers
a sequence of coefficients $c_{n}\in \mathcal{C}({\mathbb{R}}),n\in {\mathbb{%
N}},$ the study of the asymptotic behavior of $Q_{N,k_{\ast }}(c_{n},X),n\in
{\mathbb{N}}$ reduces to the study of the asymptotic behavior of $\Phi
_{m}(c_{n},G),n\in {\mathbb{N}}$, where $G=(G_{n})_{n\in {\mathbb{N}}%
},G_{n}\in {\mathbb{R}}^{m_{\ast }},$ is a sequence of independent and
centred Gaussian random variables such that $\mathrm{Cov}(G_{n})=\mathrm{Cov}%
(Z_{n}(X))$. Since $\Phi _{m}(c_{n},G)$ is (nearly) a multiple Wiener
stochastic integral of order $m,$ this problem is already treated at least
in two significant cases: the convergence to normality and the convergence
to a Gamma distribution. In fact, the convergence to normality of the law of
$\Phi _{m}(c_{n},G)$ is controlled by the Forth Moment
Theorem due to Nualart and Peccati \cite{[PN]} and
Nourdin and Peccati \cite{[NP1]}. And the convergence to a Gamma
distribution (and in particular to a $\chi _{2}$ distribution) is treated in
\cite{[NP1]}. In order to give the consequences of these results in our
framework we have to identify the link between the notation in our paper and
in the above mentioned works. Note that the coefficients $c\in \mathcal{C}({%
\mathbb{R}})$ have been defined as $c(\alpha )$ with $\alpha =(\alpha
_{1},\ldots \alpha _{m})$, $\alpha _{i}=(\alpha _{i}^{\prime },\alpha
_{i}^{\prime \prime })$, with $\alpha ^{\prime }$ on the simplex $\alpha
_{1}^{\prime }<\ldots <\alpha _{m}^{\prime }.$ We extend them by symmetry on
the whole $({\mathbb{N}}\times \lbrack m_{\ast }])^{m}$ and we denote by $%
c_{s}$ this extension (with the convention that $c_{s}(\alpha )$ is zero if $%
\alpha _{i}=\alpha _{j}$ for $i\neq j)$. So we will have
\begin{equation*}
\Phi _{m}(c,G)=\sum_{\left\vert \alpha \right\vert =m}c(\alpha )G^{\alpha }=%
\frac{1}{m!}\sum_{\left\vert \alpha \right\vert =m}c_{s}(\alpha )G^{\alpha }.
\end{equation*}%
The second point is to write the sequence of multi-dimensional random
variables $G_{n}=(G_{n,1},\ldots ,G_{n,m_{\ast }})\in {\mathbb{R}}^{m_{\ast
}},$ $n\in {\mathbb{N}}$ as a sequence of one-dimensional random variables $%
\overline{G}_{n}\in {\mathbb{R}},n\in {\mathbb{N}}$ and to re-indicate the
coefficients in a corresponding way. But we have to note first that $%
G_{n,1},\ldots ,G_{n,m_{\ast }}$ are not a priori independent, because $%
\mathrm{Cov}(G_{n})=\mathrm{Cov}(Z_{n}(X))$ is not the identity matrix. So
we assume that $\mathrm{Cov}(Z_{n}(X)$ is invertible and we first use (\ref%
{H14}) in order to write
\begin{equation*}
\Phi _{m}(c,G)=\frac{1}{m!}\sum_{\left\vert \alpha \right\vert =m}\overline{c%
}_{s}(\alpha )\overline{G}^{\alpha }
\end{equation*}%
with $\overline{c}$ defined in (\ref{H13}). Now $\overline{G}_{n,1},\ldots ,%
\overline{G}_{n,m_{\ast }}$ are independent and we are ready to write them
as a sequence. We define $I:{\mathbb{N}}\times \lbrack m_{\ast }]\rightarrow
{\mathbb{N}}$ by $I(n,j)=n\times m_{\ast }+j$. Setting $\lfloor x\rfloor$ and $\{x\}$ the
integer respectively the fractional part of $x$, the inverse function $%
J=I^{-1}:{\mathbb{N}}\rightarrow {\mathbb{N}}\times \lbrack m_{\ast }]$ is
then defined as follows: $J(n)=(\lfloor n/m_{\ast }\rfloor,\{n/m_{\ast }\}m_{\ast })$ if $%
\{n/m_{\ast }\}>0$ and $J(n)=(\lfloor n/m_{\ast }\rfloor-1,m_{\ast })$ if $\{n/m_{\ast
}\}=0$. We extend this definition to multi-indexes: if $\beta =(n_{1},\ldots
,n_{m})\in {\mathbb{N}}^{m}$ then $J(\beta )=(J(n_{1}),\ldots ,J(n_{m}))\in (%
{\mathbb{N}}\times \lbrack m_{\ast }])^{m}.$ And to coefficients: if $f:({%
\mathbb{N}}\times \lbrack m_{\ast }])^{m}\rightarrow {\mathbb{R}}$ we define
$\widehat{f}:{\mathbb{N}}^{m}\rightarrow {\mathbb{R}}^{m}$ by $\widehat{f}%
(\beta )=f(J(\beta )).$ Moreover, we consider the sequence $\widehat{G}_{n}=%
\overline{G}_{J(n)},n\in {\mathbb{N}}.$ Then%
\begin{equation*}
\Phi _{m}(c,G)=\frac{1}{m!}\sum_{\left\vert \alpha \right\vert =m}\overline{c%
}_{s}(\alpha )\overline{G}^{\alpha }=\frac{1}{m!}\sum_{\left\vert \alpha
\right\vert =m}\widehat{c}_{s}(\alpha )\widehat{G}^{\alpha }
\end{equation*}%
with the convention that now we work with the multi-index $\alpha \in {%
\mathbb{N}}^{m}.$ Note that $\Phi _{m}(\widehat{c}_{s},\widehat{G})$ is a
multiple stochastic integral of order $m.$

We introduce now the ``contraction
operators''. For $0\leq r\leq m$ and $\alpha ,\beta \in
\Gamma _{m-r}$ one denotes $\widehat{c}_{s}\otimes _{r}\widehat{c}%
_{s}(\alpha ,\beta )=\sum_{\gamma \in \Gamma _{r}}\widehat{c}_{s}(\alpha
,\gamma )\widehat{c}_{s}(\beta ,\gamma )$ with the convention that for $r=0$
we put $\widehat{c}_{s}\otimes _{0}\widehat{c}_{s}(\alpha ,\beta )=\widehat{c%
}_{s}(\alpha )\widehat{c}_{s}(\beta )$ and for $r=m,$ $\widehat{c}%
_{s}\otimes _{m}\widehat{c}_{s}=\sum_{\gamma \in \Gamma _{m}}\widehat{c}%
_{s}(\gamma )\widehat{c}_{s}(\gamma ).$ Note that, even if $\widehat{c}_{s}$
is symmetric, $\widehat{c}_{s}\otimes _{r}\widehat{c}_{s}$ is not symmetric,
so we introduce $\widehat{c}_{s}\widetilde{\otimes }_{r}\widehat{c}_{s}$ to
be the symmetrization of $\widehat{c}_{s}\otimes _{r}\widehat{c}_{s}.$

We introduce now%
\begin{equation*}
\kappa _{4,m}(\overline{c}_{s})=\sum_{r=1}^{m-1}m!^{2}\left(
\begin{tabular}{l}
$m$ \\
$r$%
\end{tabular}%
\right) ^{2}\{\left\vert \widehat{c}_{s}\otimes _{r}\widehat{c}%
_{s}\right\vert _{2m-r}^{2}+\left(
\begin{tabular}{l}
$2m-2r$ \\
$m-r$%
\end{tabular}%
\right) \left\vert \widehat{c}_{s}\widetilde{\otimes }_{r}\widehat{c}%
_{s}\right\vert _{2m-r}^{2}\}.
\end{equation*}%
It is known (see \cite{[NP1]}) that $\kappa _{4,m}(\widehat{c}_{s})$ is
equal to the forth cumulant of $\Phi _{m}(\widehat{c}_{s},\widehat{G})$ and
moreover, it is proved in \cite{[NP1]} that, if $\mathcal{N}$ is a standard
normal random variable, then%
\begin{equation}
d_{\mbox{\rm{\scriptsize{TV}}}}(\Phi _{m}(\widehat{c}_{s},\widehat{G}),%
\mathcal{N})\leq C\kappa _{4,m}^{1/2}(\widehat{c}_{s}).  \label{MR5a}
\end{equation}%
Using this and Theorem \ref{P} we immediately obtain

\begin{theorem}
\label{N} Let $\mathcal{N}$ be a standard normal random variable.

\medskip

\textbf{A}. If $X$ satisfies ${\mathfrak{M}}(\varepsilon ,r,R)$ and, for
every $n\in \N,$ $\mathrm{Cov}(Z_{n}(X)$ is invertible, then for every $%
\theta \in (\frac{1}{4},1)$ there exists $a\in (\frac{1}{2},1]$ such that
\begin{equation}
\begin{array}{l}
d_{{\mbox{\rm{\scriptsize{TV}}}}}(Q_{N,k_{\ast }}(c,X)),\mathcal{N})\leq
C\max (1,|c|_{m}^{-\frac{2}{k_{\ast }m}})^{a}(1+|c|)^{\frac{19}{2}}\smallskip
\\
\qquad \times \big(\delta _{\ast }^{\frac{\theta }{6k_{\ast }m+1}%
}(c)+e_{m,N}(c)^{a}+\left\vert c\right\vert _{m+1,N}^{\frac{2\theta }{%
k_{\ast }m}\wedge \frac{\theta }{2m+1}}+\kappa _{4,m}^{1/2}(\widehat{c}_{s})%
\big).%
\end{array}
\label{MR7}
\end{equation}

\textbf{B}. If $Z(X)$ satisfies (\ref{H1}) and (\ref{H2}) then for every $%
\theta\in (\frac 14,1)$ there exists $a\in(\frac 12,1]$ such that
\begin{equation}
d_{\mbox{\rm{\scriptsize{Kol}}}}(Q_{N,k_{\ast }}(c,X)),\mathcal{N}))\leq
C\max(1,|c|_m^{-\frac 2m})^a(1+|c|)^5 \big(\delta _{\ast }^{\frac
1{1+3N}}(c)+ \left\vert c\right\vert _{m+1,N}^{\frac \theta {2
m+1}}+e_{m,N}(c)^a +\kappa_{4,m}^{1/2}(\widehat{c}_{s})\big).  \label{MR6}
\end{equation}

\noindent In the above estimates (\ref{MR7}) and (\ref{MR6}), $C>0$ denotes
a constant independent of the coefficients $c\in\mathcal{C}({\mathbb{R}})$.
\end{theorem}

\begin{remark}
This is a generalization of the ``forth moment theorem'' to stochastic
polynomials. However there is a difference because the influence factor $%
\delta _{\ast }(c)$ appears in (\ref{MR7}). One may ask if it is possible to
control the distance between stochastic polynomials and the normal
distribution in terms of $\kappa _{4,m}(\widehat{c}_{s}))$ only. An
affirmative answer has recently been given in the following more particular
framework: assume that $d_{\ast }=k_{\ast }=1$ so that $\Phi _{m}(c,X)$ is a
multi-linear polynomial. Assume also that the random variables $X_{n},n\in {%
\mathbb{N}}$ are identically distributed. Then, if ${\mathbb{E}}%
(X_{1}^{4})\geq 3,$ the convergence to normality is controlled by $\kappa
_{4,m}(\widehat{c}_{s}))$ only (see Theorem 2.3 in \cite{[NPPS]}).
\end{remark}

We discuss now the convergence to a Gamma distribution. For $\nu \geq 1$ we
consider $F(\nu )$ a centred Gamma distribution of parameter $\nu $: $F (\nu
)=2G(\nu /2)-\nu $ where $G(\nu /2)$ has a Gamma law with parameter $\nu /2$
(that is, with density $g_{\nu/2}(x) \varpropto x^{\nu/2-1}e^{-x}1_{x>0}$).
If $\nu $ is integer then $F(\nu )$ is a centred chi-square distribution
with $\nu $ degrees of freedom. We introduce
\begin{eqnarray*}
\eta _{\nu ,m}(\widehat{c}_{s}) &=&(\nu -m!\left\vert \widehat{c}%
_{s}\right\vert _{m}^{2})^{2}+4m!\left\vert \theta _{m}\times \widehat{c}_{s}%
\widetilde{\otimes }_{m/2}\widehat{c}_{s}-\widehat{c}_{s}\right\vert
_{2m-r}^{2} \\
&&+m^{2}\sum_{\substack{ r\in \{1,\ldots ,m-1\}  \\ r\neq m/2}}%
(2m-2r)!(r-1)!^{2}\left(
\begin{array}{c}
m-1 \\
r-1%
\end{array}%
\right) ^{4}\left\vert \widehat{c}_{s}\otimes _{r}\widehat{c}_{s}\right\vert
_{2m-r}^{2}
\end{eqnarray*}%
with $\theta _{m}=\frac{1}{4}(m/2)!\left(
\begin{array}{c}
m \\
m/2%
\end{array}%
\right) .$ Combining Theorem 3.11 and Proposition 3.13 from \cite{[NP1]} one
obtains
\begin{equation*}
d_{1}(\Phi _{m}(c,Z),F(\nu ))\leq C\eta _{\nu ,m}^{1/2}(\widehat{c}_{s}).
\end{equation*}%
If $\nu $ is an integer then $F(\nu )$ has a centred $\chi ^{2}(\nu)$
distribution, so may be represented as a polynomial of degree two of
Gaussian random variables. Then, using Theorem 5.9 in \cite{[Bo]} one
obtains
\begin{equation*}
d_{\mbox{\rm{\scriptsize{TV}}}}(\Phi _{m}(c,Z),F(\nu ))\leq d_{1}^{\frac{1}{%
m+1}}(\Phi _{m}(c,Z),F(\nu ))\leq C\eta _{\nu ,m}^{1/2(m+1)}(\widehat{c}%
_{s}).
\end{equation*}%
Then, using Theorem \ref{P} we obtain

\begin{theorem}
\label{G} Let $\mathcal{X}_{\nu }$ be a random variable with a centred $\chi
^{2}$ distribution with $\nu $ degrees of freedom.

\textbf{A}. If $X$ satisfies ${\mathfrak{M}}(\varepsilon ,r,R)$ and, for
every $n\in \N,$ $\mathrm{Cov}(Z_{n}(X)$ is invertible, then for every $%
\theta \in (\frac{1}{4},1)$ there exists $a\in (\frac{1}{2},1]$ such that
\begin{equation}
\begin{array}{l}
d_{{\mbox{\rm{\scriptsize{TV}}}}}(Q_{N,k_{\ast }}(c,X)),\mathcal{X}_{\nu
})\leq C\max (1,|c|_{m}^{-\frac{2}{k_{\ast }m}})^{a}(1+|c|)^{\frac{19}{2}%
}\smallskip \\
\qquad \times \big(\delta _{\ast }^{\frac{\theta }{6k_{\ast }m+1}%
}(c)+e_{m,N}(c)^{a}+\left\vert c\right\vert _{m+1,N}^{\frac{2\theta }{%
k_{\ast }m}\wedge \frac{\theta }{2m+1}}+\eta _{\nu ,m}^{1/2(m+1)}(\widehat{c}%
_{s})\big).%
\end{array}
\label{MR9}
\end{equation}

\textbf{B}. If $Z(X)$ satisfies (\ref{H1}) and (\ref{H2}) then for every $%
\theta\in (\frac 14,1)$ there exists $a\in(\frac 12,1]$ such that
\begin{equation}
d_{\mbox{\rm{\scriptsize{Kol}}}}(Q_{N,k_{\ast }}(c,X)),\mathcal{N}))\leq
C\max(1,|c|_m^{-\frac 2m})^a(1+|c|)^5 \big(\delta _{\ast }^{\frac
1{1+3N}}(c)+ \left\vert c\right\vert _{m+1,N}^{\frac \theta {2
m+1}}+e_{m,N}(c)^a +\eta _{\nu,m}^{1/2(m+1)}(\widehat{c}_{s})\big).
\label{MR8}
\end{equation}

\noindent In the above estimates (\ref{MR9}) and (\ref{MR8}), $C>0$ denotes
a constant independent of the coefficients $c\in\mathcal{C}({\mathbb{R}})$.
\end{theorem}

\section{Examples}

\label{sect:examples}

\subsection{U-statistics associated to polynomial kernels}

\label{sect:ustatistics}

Let us first shortly recall how U-statistics appear. One considers a class
of distributions $\mathcal{M}$ and aims to estimate a functional $\theta
(\mu )$ with $\mu \in \mathcal{M}.$ In order to do it one has at hand a
sequence of independent random variables $X_{1},\ldots ,X_{n}$ with law $\mu
\in \mathcal{M},$ but does not know which is this law. The goal is to
construct an unbiased estimator, that is a sequence of functions $f_{n}:{%
\mathbb{R}}^{n}\rightarrow {\mathbb{R}},$ such that the estimator $%
U_{n}=f_{n}(X_{1},\ldots ,X_{n})$ converges to $\theta (\mu )$ and moreover $%
{\mathbb{E}}(U_{n})=\theta (\mu )$ for every $\mu \in \mathcal{M}.$ This
means that the estimator is unbiased - and this is the origin of the name
U-statistics. In 1948 Halmos \cite{[Ha]} asked the question if such an
unbiased estimator exists and if it is unique. It turns out that the
necessary and sufficient condition in order to be able to construct such an
estimator is that $\theta (\mu )$ has the following particular form: there
exists $N\in {\mathbb{N}}$ and a measurable function $\psi :{\mathbb{R}}%
^{N}\rightarrow {\mathbb{R}}$ such that
\begin{equation}
\theta (\mu )=\int_{{\mathbb{R}}^{N}}\psi (x_{1},\ldots ,x_{N})d\mu
(x_{1})\ldots d\mu (x_{N}).  \label{U1}
\end{equation}%
In this case one may construct the symmetric unbiased estimator $f_{n}$ (and
if $\mathcal{M}$ is sufficiently large, this estimator is unique in the
class of the symmetric estimators) in the following way:
\begin{equation}
U_{n}^{\psi }=\frac{(n-N)!}{n!}\sum_{(n,N)}\psi (X_{i_{1}},\ldots ,X_{i_{N}})
\label{U2}
\end{equation}%
where the sum $\sum_{(n,N)}$ is taken over all the subsets $\{i_{1},\ldots
,i_{N}\}\subset \{1,\ldots ,n\}$ such that $i_{k}\neq i_{p}$ for $k\neq p$.
It is clear that $\psi $ may be taken to be symmetric (if not one takes its
symmetrization and this change nothing).

When $\psi (x_{1},\ldots ,x_{N})$ is a polynomial, this fits in our
framework and our results apply, but, for example $\psi (x_{1},\ldots
,x_{N})=\max \{\left\vert x_{1}\right\vert ,\ldots ,\left\vert
x_{N}\right\vert \},$ is out of reach. We will treat first two standard
examples.

\medskip

\noindent \textbf{Example 1.} (Variance estimator) We denote $m_{X}={\mathbb{%
E}}(X),v_{X}={\mathbb{E}}((X-{\mathbb{E}}(X))^{2})$ and
$q_{X}=\mathrm{Var}(2m_{X}X-X^{2})$. We take $\psi (x_{1},x_{2})=\frac{1}{2}%
(x_{1}-x_{2})^{2}$ so that
\begin{equation*}
v_{X}=\int \int \psi (x_{1},x_{2})d\mu (x_{1})d\mu (x_{2}).
\end{equation*}%
In order to come back in our framework we write%
\begin{eqnarray*}
(x_{1}-x_{2})^{2} &=&2v_{X}+(x_{1}^{2}-{\mathbb{E}}(X^{2}))+(x_{2}^{2}-{%
\mathbb{E}}(X^{2}))-2(x_{1}-{\mathbb{E}}(X))(x_{2}-{\mathbb{E}}(X)) \\
&&-2m_{X}((x_{1}-{\mathbb{E}}(X))+(x_{2}-{\mathbb{E}}(X))).
\end{eqnarray*}%
It follows that
\begin{equation*}
U_{n}^{\psi }=v_{X}+\frac{1}{n}\sum_{i=1}^{n}(X_{i}^{2}-{\mathbb{E}}%
(X_{i}^{2}))-\frac{2m_{X}}{n}\sum_{i=1}^{n}(X_{i}-{\mathbb{E}}(X_{i}))+\frac{%
1}{n(n-1)}\sum_{i_{1}\neq i_{2}}(X_{i_{1}}-{\mathbb{E}}%
(X_{i_{1}}))(X_{i_{2}}-{\mathbb{E}}(X_{i_{2}})),
\end{equation*}%
thus
\begin{align*}
& \sqrt{n}(U_{n}^{\psi }-v_{X}) \\
& =\frac{1}{\sqrt{n}}\sum_{i=1}^{n}(X_{i}^{2}-{\mathbb{E}}(X^{2}))-\frac{%
2m_{X}}{\sqrt{n}}\sum_{i=1}^{n}(X_{i}-m_{X})+\frac{1}{\sqrt{n}(n-1)}%
\sum_{i_{1}\neq i_{2}}(X_{i_{1}}-m_{X})(X_{i_{2}}-m_{X}).
\end{align*}%
In our notation, we have
\begin{equation*}
\sqrt{n}(U_{n}^{\psi }-v_{X})=Q_{2,2}(c_{n},X)=S_{2}(c_{n},Z(X))
\end{equation*}%
where $c_{n}(\alpha )=0$ if $|\alpha |\neq 1,2$ and
\begin{equation*}
c_{n}(\alpha )=\left\{
\begin{array}{ll}
\displaystyle\Big(-\frac{2m_{X}}{\sqrt{n}}\,1_{\{\alpha _{1}^{\prime \prime
}=1\}}+\frac{1}{\sqrt{n}}\,1_{\{\alpha _{1}^{\prime \prime }=2\}}\Big)%
1_{\{1\leq \alpha _{1}^{\prime }\leq n\}} & \mbox{ if }|\alpha |=1\smallskip
\\
\displaystyle\frac{2}{\sqrt{n}\,(n-1)}\,1_{\{\alpha _{1}^{\prime \prime
}=\alpha _{2}^{\prime \prime }=1\}}1_{\{1\leq \alpha _{1}^{\prime }<\alpha
_{2}^{\prime }\leq n\}} & \mbox{ if }|\alpha |=2.%
\end{array}%
\right.
\end{equation*}%
The quantities which come on in our convergence theorem are%
\begin{equation*}
\delta _{\ast }^{2}(c_{n})=\frac{4m_{X}^{2}+1}{n}+\frac{4}{n(n-1)},\quad
|c_{n}|_{1}^{2}=4m_{X}^{2}+1,\quad |c_{n}|_{2}^{2}=\frac{4}{n-1}%
=|c_{n}|_{2,2}^{2}.
\end{equation*}%
Our invariance principle (Theorem \ref{Main}) says that $Q_{2,2}(c_{n},X)$
is asymptotically equivalent in total variation distance with
\begin{equation*}
-\frac{2m_{X}}{\sqrt{n}}\sum_{j=1}^{n}G_{1,j}+\frac{1}{\sqrt{n}}%
\sum_{j=1}^{n}G_{2,j}=:-2m_{X}G_{1}+G_{2}
\end{equation*}%
where $G_{j}=(G_{1,j},G_{2,j})$ are Gaussian random variables with the same
mean and covariance as $(X-{\mathbb{E}}(X),X^{2}-{\mathbb{E}}(X^{2})).$ Then
$-2m_{X}G_{1}+G_{2}$ is a centred Gaussian random variable with variance $%
q_{X}=\mathrm{Var}(2m_{X}X-X^{2})$ so, if ${\mathfrak{M}}(\varepsilon ,r,R)$
holds, then Theorem \ref{Main} and Theorem \ref{P} yield%
\begin{equation*}
d_{{\mbox{\rm{\scriptsize{TV}}}}}(\sqrt{n}(U_{n}^{\psi }-v_{X}),\sqrt{q_{X}}%
\Delta )\leq \frac{C}{n^{\theta /13}}
\end{equation*}%
for every $\theta <1$, with $\Delta $ a standard normal random variable.

\begin{remark}
Another way to do things, used in U-statistics theory, is the following. One
employs the two dimensional CLT in order to prove that the term normalized
with $1/\sqrt n$ converges in law to $\sqrt{q_{X}}\Delta $ and then one
notes that the remaining term is smaller, so it may be ignored.
\end{remark}

\textbf{Example 2}. We look to the U-statistics associated to $\psi
(x_{1},x_{2})=x_{1}x_{2}.$ We set $m_X={\mathbb{E}}(X)$ and $v_X=\mathrm{Var}%
(X)$. Here $\psi $ is not invariant with respect to translations and we have
two different limits according to the fact that $m_{X}$ is null or not. We
write%
\begin{equation*}
x_{1}x_{2}=(x_{1}-m_{X})(x_{2}-m_{X})+m_{X}((x_{1}-m_{X})+(x_{2}-m_{X}))+m_{X}^{2}
\end{equation*}%
so that
\begin{align*}
U_{n}^{\psi }&=\frac{1}{n(n-1)}\sum_{i_{1}\neq i_{2}}X_{i_{1}}X_{i_{2}} \\
&=m_{X}^{2}+\frac{2m_{X}}{n}\sum_{i=1}^{n}(X_{i}-{\mathbb{E}}(X_{i}))+\frac{1%
}{n(n-1)}\sum_{i_{1}\neq i_{2}}(X_{i_{1}}-{\mathbb{E}}(X_{i_{1}}))(X_{i_{2}}-%
{\mathbb{E}}(X_{i_{2}})).
\end{align*}

\smallskip

\noindent \textbf{Case 1: $m_{X}\neq 0$.} Then
\begin{align*}
\sqrt{n}(U_{n}^{\psi }-m_{X}^{2}) &=\frac{2m_{X}}{\sqrt{n}}%
\sum_{i=1}^{n}(X_{i}-m_{X}) +\frac{1}{\sqrt{n}(n-1)}\sum_{i_{1}\neq
i_{2}}(X_{i_{1}}-m_{X})(X_{i_{2}}-m_{X}) \\
&=Q_{2,1}(c_n,X)=S_2(c_n,X)
\end{align*}
with $c_n(\alpha)=0$ if $|\alpha|\neq 1,2$ and
\begin{equation*}
c_n(\alpha)=\left\{
\begin{array}{ll}
\displaystyle \frac{2m_X}{\sqrt n}1_{\{1\leq \alpha^{\prime }_1\leq n\}} & %
\mbox{ if } |\alpha|=1\smallskip \\
\displaystyle \frac{2}{\sqrt n\,(n-1)}1_{\{1\leq \alpha^{\prime
}_1<\alpha^{\prime }_2\leq n\}} & \mbox{ if } |\alpha|=2%
\end{array}
\right.
\end{equation*}
One has
\begin{equation*}
\delta_\ast^2(c_n) =\frac{4m_X^2}{n}+\frac 4{n(n-1)},\quad
|c_n|^2_1=4m_X^2,\quad |c_n|_2^2=\frac 4{n-1}=|c_n|^2_{2,2}.
\end{equation*}

Using Theorem \ref{Main} and Theorem \ref{P}, the asymptotic behavior of $%
\sqrt{n}(U_{n}-m_{X}^{2})$ is equivalent to the behavior of%
\begin{equation*}
\frac{m_{X}}{\sqrt{n}}\sum_{i=1}^{n}G_{i}=m_{X}\sqrt{v_X}\,\Delta
\end{equation*}%
with $\Delta $ standard normal.

\smallskip

\noindent \textbf{Case 2: $m_{X}=0$.} Then%
\begin{equation*}
nU_{n}^{\psi }=\frac{1}{n-1}\sum_{i_{1}\neq i_{2}}X_{i_{1}}X_{i_{2}}
=Q_{2,1}(c_n,X)=S_2(c_n,X)
\end{equation*}
where $c_n(\alpha)=0$ if $|\alpha|\neq 2$ and
\begin{equation*}
c_n(\alpha)= \frac{2}{n-1}\,1_{\{1\leq \alpha^{\prime }_1<\alpha^{\prime
}_2\leq n\}}\quad \mbox{if}\quad |\alpha|=2.
\end{equation*}
Here,
\begin{equation*}
\delta_*^2(c_n)= \frac{4}{n-1},\quad |c_n|_2^2=\frac {4n}{n-1}.
\end{equation*}
Using the invariance principle (Theorem \ref{Main}) this is close to $\frac{%
v_X}{n-1}\sum_{i_{1}\neq i_{2}}G_{i_{1}}G_{i_{2}}$ with $G_{i},i\in {\mathbb{%
N}}$ independent standard normal random variables. We define $%
D_{n}=[0,1]^{2}\smallsetminus \cup _{i=0}^{n-1}[\frac{i}{n},\frac{i+1}{n})^2$
and $f_{n}(s_{1},s_{2})=\frac{nv_X}{n-1}1_{D_{n}}(s_{1},s_{2}).$ Then the
law of $\frac{v_X}{n-1}\sum_{i_{1}\neq i_{2}}G_{i_{1}}G_{i_{2}}$ coincides
with the law of the double It\^{o} integral $I_{2}(f_{n}).$ Setting $f\equiv
v_X$, we recall that the law of $I_{2}(f)$ coincides with the law of $%
v_X(\Delta ^{2}-1)$ where $\Delta $ is standard normal. Then, using Theorem %
\ref{Main} (with $k_{\ast }=1,N=m=2)$ and Theorem \ref{mul} (with $k=1,m=2)$
one obtains, for every $\theta <1,$%
\begin{eqnarray*}
d_{\mbox{\rm{\scriptsize{TV}}}}(nU_{n}^{\psi },v_X(\Delta ^{2}-1)) &\leq &d_{%
\mbox{\rm{\scriptsize{TV}}}}(nU_{n}^{\psi },I_{2}(f_{n}))+d_{%
\mbox{\rm{\scriptsize{TV}}}}(I_{2}(f_{n}),I_{2}(f)) \\
&\leq &\frac{C}{n^{\theta /13}}+Cd_{1}^{\theta
/5}(I_{2}(f_{n}),I_{2}(f))\leq \frac{C}{n^{\theta /13}}\rightarrow 0.
\end{eqnarray*}

An alternative way to solve the problem is to write%
\begin{equation*}
nU_{n}^{\psi }=\frac{n}{n-1}\Big( \Big(\frac{1}{\sqrt{n}}\sum_{i=1}^{n}X_{i}%
\Big)^{2}-\frac{1}{n}\sum_{i=1}^{n}X_{i}^{2} \Big)
\end{equation*}%
and to use the CLT in order to replace $\frac{1}{\sqrt{n}}%
\sum_{i=1}^{n}X_{i} $ with $\sqrt{v_X}\Delta $ and to say that by the law of
large numbers the last term goes to $v_X$. This gives the convergence in law
of $nU_{n}^{\psi } $ to $v_X(\Delta ^{2}-1).$

\begin{remark}
The above two examples suggest the following rough comparison of the
strategies employed in the U-statistics theory on one hand and in our paper
on the other hand. In the U-statistics theory one tries to make blocks of
terms such that in the end $U_{n}^{\psi }$ appears as a continuous function
of blocks of the form $\frac{1}{\sqrt{n}}\sum_{i=1}^{n}Y_{i}$ or $\frac{1}{n}%
\sum_{i=1}^{n}Y_{i}^{2}$ and then use the CLT, respectively the law of
large numbers, in order to replace them, asymptotically, by a Gaussian
random variable respectively by a constant. Alternatively, in our paper one
begins by using the invariance principle in order to change $X_{i}-{\mathbb{E%
}}(X_{i})$ and $X_{i}^{2}-{\mathbb{E}}(X_{i}^{2})$ by Gaussian random
variables $G_{i,1}$ and $G_{i,2}.$ And then one solves the problem of the
asymptotic behavior in the framework of Wiener chaoses.
\end{remark}

Let us go on and look to general polynomials. We fix $k_{\ast },N\in {%
\mathbb{N}}$, we denote $\mathcal{K}_{N}=\{0,1,\ldots ,k_{\ast }\}^{N},$
and we define
\begin{equation}
\psi (x_{1},\ldots ,x_{N})=\sum_{\kappa \in \mathcal{K}_{N}}a(\kappa
)x^{\kappa }\quad \mbox{with}\quad x^{\kappa }=\prod_{j=1}^{N}x_{j}^{k_{j}}
\label{U3}
\end{equation}%
with symmetric coefficients $a(\kappa )$ which are null on the diagonals$.$
So $\psi $ is a general symmetric polynomial of order $k_{\ast }$ in the
variables $x_{1},\ldots ,x_{N}.$ We associate to $\psi $ the U-statistic $%
U_{n}^{\psi }$ defined in (\ref{U2}):%
\begin{equation}
U_{n}^{\psi }=\frac{(n-N)!}{n!}\sum_{i_{1},\ldots ,i_{N}}\psi
(X_{i_{1}},\ldots ,X_{i_{N}})=\left(
\begin{tabular}{l}
$n$ \\
$N$%
\end{tabular}%
\right) ^{-1}\sum_{i_{1}<\ldots <i_{N}}\sum_{\kappa \in \mathcal{K}%
_{N}}a(\kappa )\prod_{j=1}^{N}X_{i_{j}}^{k_{j}}.  \label{U4}
\end{equation}

The above quantity is linked with the stochastic polynomials defined in the
previous sections in the following way. One takes $d_{\ast }=1$ and $m_{\ast
}=k_{\ast }$ and constructs coefficients $c_{n}$ such that $U_{n}^{\psi
}=Q_{N,k_\ast}(c_{n},X)=S_{N}(c_{n},Z(X))$ with $Z(X)$ associated to $X$ in (%
\ref{H15}): $Z_{i,k}(X)=X_{i}^{k}-{\mathbb{E}}(X_{i}^{k}),k=1,\ldots
,k_{\ast }.$ The problem is that $Z_{i,k}(X)$ is centred whereas $X_{i}^{k},$
which appears in (\ref{U4}), is not. I turns out that the operation which
consists in centering $X_{i}^{k}$ in (\ref{U4}) is exactly the Hoeffding
decomposition, introduced by Hoeffding in \cite{[Hof1],[Hof2]}, and which plays a
crucial role in the theory of U-statistics. Let us recall it. For $1\leq
j\leq N,$ one defines the kernels
\begin{equation*}
h_{j}(x_{1},\ldots ,x_{j})=\int \ldots \int \psi (u_{1},\ldots
,u_{N})\prod_{i=1}^{j}(\delta _{x_{i}}-\mu )(du_{i})\prod_{i=j+1}^{N}\mu
(du_{i}).
\end{equation*}%
Then Hoeffding's decomposition is the following:%
\begin{equation}
U_{n}^{\psi }=\theta (\mu )+\sum_{j=1}^{N}\left(
\begin{tabular}{l}
$N$ \\
$j$%
\end{tabular}%
\right) U_{n}^{h_{j}}  \label{U5}
\end{equation}%
where $U_{n}^{h_{j}}$ is the U-statistic associated to $h_{j}$ in the first
equality from (\ref{U4}) (with $N$ replaced by $j)$. See for example Theorem
1 in Section 1.6 in \cite{[Lee]} for the proof of (\ref{U5}).

We denote $m_{k}={\mathbb{E}}(X^{k})$ and we compute

\begin{equation*}
\int \ldots \int \prod_{l=1}^{N}u_{l}^{k_{l}}\prod_{i=1}^{j}(\delta
_{x_{i}}-\mu )(du_{i})\prod_{i=j+1}^{N}\mu
(du_{i})=\prod_{i=1}^{j}(x_{i}^{k_{i}}-m_{k_{i}})\times
\prod_{i=j+1}^{N}m_{k_{i}}
\end{equation*}%
so we obtain%
\begin{eqnarray*}
h_{j}(x_{1},\ldots ,x_{j}) &=&\sum_{\kappa \in \mathcal{K}_{j}}a_{j}(\kappa
)\prod_{i=1}^{j}(x_{i}^{k_{i}}-m_{k_{i}})\quad \mbox{with} \\
a_{j}(\kappa ) &=&\sum_{k_{j+1},\ldots ,k_{N}=1}^{k_{\ast }}a(\kappa
,k_{j+1},\ldots ,k_{N})\prod_{i=j+1}^{N}m_{k_{i}}.
\end{eqnarray*}%
We conclude that
\begin{equation*}
U_{n}^{\psi }=\theta (\mu )+\sum_{j=1}^{N}\left(
\begin{tabular}{l}
$N$ \\
$j$%
\end{tabular}%
\right) \left(
\begin{tabular}{l}
$n$ \\
$j$%
\end{tabular}%
\right) ^{-1}\sum_{i_{1}<\cdots <i_{j}}a_{j}(\kappa
)\prod_{l=1}^{j}(X_{i_{l}}^{k_{l}}-{\mathbb{E}}(X_{i_{l}}^{k_{l}})).
\end{equation*}

In the theory of U-statistics one says that $U_{n}^{\psi }$ is degenerated
at order $m\in \lbrack N]$ if $h_{j}=0$ for $j\leq m-1$ and $h_{m}\neq 0,$
which amounts to
\begin{equation}  \label{nd1-2}
\sum_{\kappa \in \mathcal{K}_{j}}a_{j}^{2}(\kappa )=0\quad 1\leq j\leq m-1
\quad\mbox{and}\quad \sum_{\kappa \in \mathcal{K}_{m}}a_{m}^{2}(\kappa )>0.
\end{equation}

We assume that (\ref{nd1-2}) holds and we write
\begin{equation*}
V_{m}(n):=n^{m/2}(U_{n}^{\psi }-\theta (\mu ))=\sum_{j=m}^{N}n^{m/2}\left(
\begin{tabular}{l}
$N$ \\
$j$%
\end{tabular}%
\right) U_{n}^{h_{j}}=\sum_{j=m}^{N}\sum_{\left\vert \alpha \right\vert
=j}c_{n}(\alpha )Z^{\alpha }(X)
\end{equation*}%
with
\begin{equation*}
c_{n}((i_{1},k_{1}),\ldots ,(i_{j},k_{j}))=n^{m/2}\left(
\begin{tabular}{l}
$N$ \\
$j$%
\end{tabular}%
\right) \times \left(
\begin{tabular}{l}
$n$ \\
$j$%
\end{tabular}%
\right) ^{-1}a_{j}(k_{1},\ldots ,k_{j}).
\end{equation*}%
By (\ref{nd1-2}), the U-statistic $U_{n}^{\psi }$ is degenerated at order $%
m\in \lbrack N]$ if and only if
\begin{equation*}
|c_n|_j=0\quad\mbox{for}\quad j\leq m-1\quad\mbox{and}\quad |c_n|_m>0,
\end{equation*}
which is the same non-degeneracy condition we are interested in.

We recall that $X_{i}\sim \mu $ and that in (\ref{H2}) we have introduced
the covariance matrix $\mathrm{Cov}(Z(X))=\mathrm{Cov}(\mu)$, that is
\begin{equation*}
\mathrm{Cov}^{i,j} (\mu)={\mathbb{E}}((X^{i}-{\mathbb{E}}(X^{i}))(X^{j}-{%
\mathbb{E}}(X^{j})).
\end{equation*}%
We consider a correlated Brownian motion $W=(W^{1},\ldots ,W^{m})$ with $%
\left\langle W^{i},W^{j}\right\rangle _{t}=C^{i,j}(\mu )t,$ we define the
multiple stochastic integrals%
\begin{equation*}
I_{\kappa }^{\mu
}(1)=\int_{0}^{1}dW_{s_{m}}^{k_{m}}\int_{0}^{s_{m}}dW_{s_{m-1}}^{k_{m-1}}%
\ldots \int_{0}^{s_{2}}1dW_{s_{1}}^{k_{1}}
\end{equation*}%
and we denote%
\begin{equation*}
V_{m}=\left(
\begin{array}{c}
N \\
N-m%
\end{array}%
\right) \sum_{\kappa \in \mathcal{K}_{m}}a_{m}(\kappa )I_{\kappa }^{\mu }(1).
\end{equation*}

\begin{theorem}
\label{UU} \textbf{A}. If $X$ verifies ${\mathfrak{M}}(\varepsilon,r,R)$ and
(\ref{nd1-2}) holds then for every $\theta \in(\frac 14, 1)$%
\begin{equation}
d_{\mbox{\rm{\scriptsize{TV}}}}(V_{m}(n),V_{m})\leq \frac{C}{n^{\theta \beta
(m,k_{\ast })}}\quad \mbox{with}\quad \beta (m,k_{\ast })=\frac 1
{2(6k_{\ast }m+1)}.  \label{U8}
\end{equation}

\textbf{B}. Suppose that $X$ has finite moments of any order and that $%
\mathrm{Cov}(Z(X))=\mathrm{Cov}(\mu )\geq \underline{\lambda }>0.$ If (\ref%
{nd1-2}) holds then, for every $\theta \in(\frac 14, 1)$%
\begin{equation}
d_{\mbox{\rm{\scriptsize{Kol}}}}(V_{m}(n),V_{m})\leq \frac{C}{n^{\theta
\alpha (N)}}\quad \mbox{with}\quad \alpha (N)=\frac{1}{2(3N+1)}.  \label{U7}
\end{equation}
\end{theorem}

\textbf{Proof}. In order to use Theorem \ref{P} we estimate%
\begin{eqnarray*}
\left\vert c_{n}\right\vert _{m+1,N}^{2} &=&\sum_{m+1\leq \left\vert \alpha
\right\vert \leq N}c_{n}^{2}(\alpha )\leq Cn^{m}\times
\sum_{j=m+1}^{N}n^{-2j}\times n^{j}\times \left\Vert a\right\Vert _{\infty
}\leq \frac{C}{n}, \\
\left\vert c_{n}\right\vert _{m}^{2} &=&\sum_{\left\vert \alpha \right\vert
=m}c_{n}^{2}(\alpha )\geq \frac{1}{C}\times n^{m}\times n^{-2m}\times
n^{m}\times \sum_{\kappa \in \mathcal{K}_{m}}a_{m}^{2}(\kappa )=\frac{1}{C}%
\times \sum_{\kappa \in \mathcal{K}_{m}}a_{m}^{2}(\kappa )>0.
\end{eqnarray*}%
Finally we study the influence factor:
\begin{equation*}
\delta _{\ast }(c_{n})=\max_{r}\sum_{m\leq \left\vert \alpha \right\vert
\leq N}c_{n}^{2}(\alpha )1_{\{r\in \alpha ^{\prime }\}}\leq Cn^{m}\times
\sum_{j=m}^{N}n^{-2j}\times n^{j-1}=\frac{C}{n}.
\end{equation*}%
Then (\ref{MR5}) gives
\begin{equation*}
d_{{\mbox{\rm{\scriptsize{TV}}}}}(Q_{N,k_{\ast }}(c,X)),\Phi _{m}(c,G)))\leq
C\big(\Big(\frac{1}{\sqrt{n}}\Big)^{\frac{\theta }{6k_{\ast }m+1}}+\Big(%
\frac{1}{\sqrt{n}}\Big)^{\frac{2\theta }{k_{\ast }m}\wedge \frac{\theta }{%
2m+1}}\big)\leq C\frac{1}{n^{\frac{\theta }{2(6k_{\ast }m+1)}}}.
\end{equation*}%
And by employing (\ref{MR4}) one has
\begin{equation*}
d_{{\mbox{\rm{\scriptsize{Kol}}}}}(Q_{N,k_{\ast }}(c,X)),\Phi
_{m}(c,G)))\leq C\Big(\Big(\frac{1}{\sqrt{n}}\Big)^{\frac{1}{1+3N}}+\Big(%
\frac{1}{\sqrt{n}}\Big)^{\frac{\theta }{2m+1}}\Big)\leq C\frac{1}{n^{\frac{%
\theta }{2(1+3N)}}}.
\end{equation*}%
$\square $

\subsection{A quadratic central limit theorem}

\label{sect:quadratic} For $p\in(0,\frac 12]$, we look to the quadratic form%
\begin{equation*}
S_{n,p}(Z) =\left\{
\begin{array}{ll}
\displaystyle \frac{1}{n^{1-p}}\sum_{i,j=1}^{n}1_{\{i\neq j\}}\frac{1}{%
\left\vert i-j\right\vert ^{p}}Z_{i}Z_{j} & \mbox{ if } 0<p<\frac{1}{2}, \\
\displaystyle \frac{1}{(2n\ln n)^{1/2}}\sum_{i,j=1}^{n}1_{\{i\neq j\}}\frac{1%
}{\left\vert i-j\right\vert ^{1/2}}Z_{i}Z_{j} & \mbox{ if } p=\frac{1}{2},%
\end{array}%
\right.
\end{equation*}
where $Z_{i},i\in {\mathbb{N}}$ are centred independent random variables
which have finite moments of any order. The aim of this section is to prove
that if $p<\frac{1}{2}$ then $S_{n,p}(Z)$ converges to a double stochastic
integral while for $p=\frac{1}{2}$ the limit is a standard Gaussian random
variable. In our notation, we have $d_*=1$, $k_*=1$, $N=2$ and
\begin{equation*}
S_{n,p}(Z)=Q_{2,1}(c_{n,p},Z)=S_2(c_{n,p},Z)
\end{equation*}
where $c_{n,p}(\alpha)=0$ for $|\alpha|\neq 2$ and if $|\alpha|=2$,
\begin{equation}  \label{cnp}
c_{n,p}(\alpha) =\left\{
\begin{array}{ll}
\displaystyle \frac 2{n^{1-p}|\alpha^{\prime }_1-\alpha^{\prime
}_2|^p}\,1_{\{1\leq \alpha^{\prime }_1<\alpha^{\prime }_2\leq n\}} &
\mbox{
if }0<p<\frac 12,\smallskip \\
\displaystyle \frac 2{(2n\ln n)^{1/2}|\alpha^{\prime }_1-\alpha^{\prime
}_2|^{1/2}}\,1_{\{1\leq \alpha^{\prime }_1<\alpha^{\prime }_2\leq n\}} & %
\mbox{ if }p=\frac 12.%
\end{array}
\right.
\end{equation}

\begin{theorem}
\label{Q} Let $Z_{i},i\in {\mathbb{N}}$ be a sequence of independent and
centred random variables, with ${\mathbb{E}}(Z_{i}^{2})=1$ and which have
finite moments of any order.

\medskip

\textbf{A}. Let $p<\frac{1}{2}$. We denote $\psi _{p}(s,t)=\left\vert
s-t\right\vert ^{-p}$ and $I_{2}(\psi _{p})=\int_{0}^{1}\int_{0}^{1}\psi
_{p}(s,t)dW_{s}dW_{t}$, $W$ being a Brownian motion. Then for every $%
\theta\in(\frac 14,1)$ there exists $n_{\ast }$ and $C$ such that for $n\geq
n_{\ast }$%
\begin{equation}
d_{\mbox{\rm{\scriptsize{Kol}}}}(S_{n,p},I_{2}(\psi _{p}))\leq \frac{C}{n^{%
\frac{\theta(1-2p)}{15}}}.  \label{Q1}
\end{equation}%
Suppose moreover that ${\mathfrak{D}}(\varepsilon, r, R)$ holds. Then for
every $\theta\in(\frac 14,1)$ there exists $n_{\ast }$ and $C$ such that for
$n\geq n_{\ast }$%
\begin{equation}  \label{Q1bis}
d_{\mbox{\rm{\scriptsize{TV}}}}(S_{n,p},I_{2}(\psi _{p}))\leq \frac{C}{%
n^{\frac\theta{26}\wedge\frac{\theta(1-2p)}{15}}}.
\end{equation}

\textbf{B}. Let $p=\frac{1}{2}$. We denote $\Delta $ a standard normal
random variable. There exists $n_{\ast }$ and $C$ such that for $n\geq
n_{\ast }$%
\begin{equation}
d_{\mbox{\rm{\scriptsize{Kol}}}}(S_{n,1/2},\Delta )\leq \frac{C}{(\ln
n)^{1/2}}.  \label{Q2}
\end{equation}%
Suppose moreover that ${\mathfrak{D}}(\varepsilon, r, R)$ holds. Then (\ref%
{Q2}) holds with $d_{\mbox{\rm{\scriptsize{TV}}}}$ instead of $d_{%
\mbox{\rm{\scriptsize{Kol}}}}.$
\end{theorem}

\textbf{Proof A}. We extend by symmetry the coefficients $c_{n,p}(\alpha)$
to all indexes $\alpha=(\alpha_1,\alpha_2)$ with $\alpha_1\neq \alpha_2$. We
denote $t_{i}=\frac{i}{n}$ and we define%
\begin{equation*}
\psi _{n,p}(s,t)=\frac{1}{n}1_{i\neq j}\frac{1}{\left\vert
t_{i}-t_{j}\right\vert ^{p}}%
1_{[t_{i},t_{i+1})}(s)1_{[t_{j},t_{j+1})}(t)=c_{n,p}(i,j)1_{[t_{i},t_{i+1})}(s)1_{[t_{j},t_{j+1})}(t).
\end{equation*}%
Let us prove that
\begin{equation}
\int_{0}^{1}\int_{0}^{1}\left\vert \psi _{p}(s,t)-\psi
_{n,p}(s,t)\right\vert ^{2}dsdt\leq \frac{C}{n^{\frac{2}{3}(1-2p)}}.
\label{Q3}
\end{equation}%
We take $q=\frac{2}{3}$ and we write%
\begin{equation*}
\int_{0}^{1}\int_{0}^{1}\left\vert \psi _{p}(s,t)-\psi
_{n,p}(s,t)\right\vert ^{2}dsdt\leq I+J+J^{\prime }
\end{equation*}%
with%
\begin{eqnarray*}
I &=&\int_{\left\vert s-t\right\vert \geq 1/n^{q}}\left\vert \psi
_{p}(s,t)-\psi _{n,p}(s,t)\right\vert ^{2}dsdt, \\
J &=&\int_{\left\vert s-t\right\vert <1/n^{q}}\left\vert \psi
_{p}(s,t)\right\vert ^{2}dsdt,\quad J^{\prime }=\int_{\left\vert
s-t\right\vert <1/n^{q}}\left\vert \psi _{n,p}(s,t)\right\vert ^{2}dsdt.
\end{eqnarray*}%
Note that if $\left\vert s-t\right\vert \geq 1/n^{q}$ then
\begin{equation*}
\left\vert \psi _{p}(s,t)-\psi _{n,p}(s,t)\right\vert \leq \frac{C}{n}\times
\frac{1}{\left\vert s-t\right\vert ^{p+1}}\leq \frac{C}{n^{1-q(p+1)}}
\end{equation*}%
so that
\begin{equation*}
I\leq \frac{C}{n^{2(1-q(p+1))}}.
\end{equation*}%
Moreover%
\begin{equation*}
J=2\int_{0}^{1}dt\int_{0}^{t+\frac{1}{n^{q}}}\frac{ds}{\left\vert
s-t\right\vert ^{2p}}=\frac{C}{n^{q(1-2p)}}.
\end{equation*}%
Finally, by comparing Riemann sums with the corresponding integral,
\begin{eqnarray*}
J^{\prime } &=&\frac{1}{n^{2}}\sum_{i=1}^{n}\sum_{0<\left\vert
t_{i}-t_{j}\right\vert \leq 1/n^{q}}\frac{1}{\left\vert
t_{i}-t_{j}\right\vert ^{2p}} \\
&\leq &\frac{1}{n^{2}}(2n+\sum_{i=1}^{n}\sum_{\substack{ 0<\left\vert
t_{i}-t_{j}\right\vert \leq 1/n^{q}  \\ \left\vert i-j\right\vert \geq 2}}%
\frac{1}{\left\vert t_{i}-t_{j}\right\vert ^{2p}})\leq \frac{1}{n^{2}}%
(2n+J)\leq \frac{C}{n^{q(1-2p)}}.
\end{eqnarray*}%
Since $q=\frac{2}{3}$ we obtain (\ref{Q3}). It follows that, for
sufficiently large $n,$%
\begin{equation*}
\frac{1}{2}\int_{0}^{1}\int_{0}^{1}\left\vert \psi _{p}(s,t)\right\vert
^{2}dsdt\leq \left\vert c_{n}\right\vert
^{2}=\int_{0}^{1}\int_{0}^{1}\left\vert \psi _{n,p}(s,t)\right\vert
^{2}dsdt\leq 2\int_{0}^{1}\int_{0}^{1}\left\vert \psi _{p}(s,t)\right\vert
^{2}dsdt.
\end{equation*}%
And we also have%
\begin{equation*}
\delta _{\ast }^{2}(c_{n,p})=\max_{i\leq n}\frac{1}{n}\sum_{j\neq i}\frac{1}{%
n}\frac{1}{\left\vert t_{i}-t_{j}\right\vert ^{2p}}\leq \frac{C}{n}.
\end{equation*}%
Note that $S_{n,p}(Z)=S_{2}(c_{n},Z)$ and $S_{2}(c_{n},G)=I_{2}(\psi
_{n,p}). $ Using Theorem \ref{MDO} (with $N=2$), Theorem \ref{mul} (see (\ref%
{D5}) with $k=1,m=2,\frac 14<\theta <1$) and (\ref{Q3}) we obtain
\begin{eqnarray*}
d_{\mbox{\rm{\scriptsize{Kol}}}}(S_{n,p}(Z),I_{2}(\psi _{p})) &\leq &d_{%
\mbox{\rm{\scriptsize{Kol}}}}(S_{2}(c_{n,p},Z),S_{2}(c_{n,p},G))+d_{%
\mbox{\rm{\scriptsize{Kol}}}}(I_{2}(\psi _{n,p}),I_{2}(\psi _{p})) \\
&\leq &C(\delta _{\ast }^{1/7}(c_{n,p})+\left\Vert \psi _{p}-\psi
_{n,p}\right\Vert _{2}^{\theta /5}\leq C(\frac{1}{n^{1/14}}+\frac{1}{n^{%
\frac{\theta(1-2p)}{15}}})
\end{eqnarray*}%
so (\ref{Q1}) is proved for $d_{\mbox{\rm{\scriptsize{Kol}}}}.$

We suppose now that $Z$ verifies (\ref{N8}) and we use Theorem \ref{Main}
(see (\ref{MR1}) with $N=2)$ in order to obtain%
\begin{equation*}
d_{\mbox{\rm{\scriptsize{TV}}}}(S_{n,p}(Z),I_{2}(\psi _{p}))\leq C(\delta
_{\ast }^{\theta /13}(c_{n})+\left\Vert \psi _{p}-\psi _{n,p}\right\Vert
_{2}^{\theta /5}\leq C(\frac{1}{n^{1/26}}+\frac{1}{n^{\frac{2}{15}(1-2p)}})
\end{equation*}%
so (\ref{Q2}) is proved for $d_{\mbox{\rm{\scriptsize{TV}}}}$ also.

\medskip

\textbf{B}. We have $S_{n,1/2}(Z)=S_{2}(c_{n},Z)$ with (recall that $%
t_{i}=i/n)$%
\begin{equation*}
c_{n}(i,j)=\frac{1}{\sqrt{2n\ln n}}1_{i\neq j}\frac{1}{\left\vert
i-j\right\vert ^{1/2}}=\frac{1}{\sqrt{2\ln n}}1_{i\neq j}\frac{1}{\left\vert
t_{i}-t_{j}\right\vert ^{1/2}}.
\end{equation*}%
We note first that
\begin{equation*}
\ln i+\ln (n-i)\leq \sum_{j=1}^{n}1_{i\neq j}\left\vert i-j\right\vert
^{-1}\leq 2+\ln i+\ln (n-i).
\end{equation*}%
These inequalities are easily obtained by comparing $\sum_{j=1}^{n}1_{i\neq
j}\left\vert i-j\right\vert ^{-1}$ with $\int_{\{\left\vert
t_{i}-y\right\vert >1/n\}}\left\vert t_{i}-t\right\vert ^{-1}dt.$ It
immediately follows that%
\begin{equation*}
1-\frac{1}{\ln n}\leq \left\vert c_{n}\right\vert ^{2}\leq 1+\frac{1}{\ln n}
\end{equation*}%
and $\delta _{\ast }(c_{n})\leq \frac{\sqrt{2}}{\sqrt{n}}.$ Now, using
Theorem \ref{MDO}%
\begin{equation*}
d_{\mbox{\rm{\scriptsize{Kol}}}}(S_{2}(c_{n},Z),S_{2}(c_{n},G))\leq \frac{C}{%
n^{1/14}}
\end{equation*}%
and, if $Z_{i}$ satisfies ${\mathfrak{D}}(\varepsilon, r, R)$, we use
Theorem \ref{Main} and we obtain
\begin{equation*}
d_{\mbox{\rm{\scriptsize{TV}}}}(S_{2}(c_{n},Z),S_{2}(c_{n},G))\leq \frac{C}{%
n^{\theta/26}}.
\end{equation*}%
Now we have to estimate the total variation distance between $%
S_{2}(c_{n},G)=\Phi_2(c_n,G)$ and the normal random variable $\Delta .$ In
order to do it we use (\ref{MR5a}), so we have to estimate the kurtosis $%
\kappa (c_{n}).$ We denote $a(i,j)=1_{i\neq j}\left\vert i-j\right\vert
^{-1/2}$  and we write
\begin{eqnarray*}
a\otimes _{1}a(i,j) &=&\sum_{k}1_{k\neq i}1_{k\neq j}\frac{1}{\sqrt{%
\left\vert t_{i}-t_{k}\right\vert \left\vert t_{j}-t_{k}\right\vert }}\times
\frac{1}{n} \\
&\leq &2+\sum_{k<\lfloor\frac{i+j}{2}\rfloor}\frac{1}{\sqrt{\left\vert
t_{i}-t_{k}\right\vert \left\vert t_{j}-t_{k}\right\vert }}1_{k\neq
i}1_{k\neq j}+\sum_{k>\lfloor\frac{i+j}{2}\rfloor+1}\frac{1}{\sqrt{\left\vert
t_{i}-t_{k}\right\vert \left\vert t_{j}-t_{k}\right\vert }}1_{k\neq
i}1_{k\neq j} \\
&\leq &2+\int_{0}^{1}\frac{dt}{\sqrt{\left\vert t_{i}-t\right\vert
\left\vert t_{j}-t\right\vert }}.
\end{eqnarray*}%
In order to obtain the last inequality one just looks to the graphs of the
functions $t\mapsto (\left\vert t_{i}-t\right\vert \left\vert
t_{j}-t\right\vert )^{-1/2}$ and to the graph of the step approximation of
this function. And the step approximation is below the function in these
regions. Moreover (see \cite{[BCNon]} Lemma B1 for a complete computation)%
\begin{equation*}
\int_{0}^{1}\frac{dt}{\sqrt{\left\vert t_{i}-t\right\vert \left\vert
t_{j}-t\right\vert }}=\pi +2\ln \frac{\sqrt{1-t_{i}}+\sqrt{1-t_{j}}}{%
\left\vert \sqrt{t_{i}}-\sqrt{t_{j}}\right\vert }.
\end{equation*}%
It follows that
\begin{eqnarray*}
\kappa ^{2}(c_{n}) &=&\left\vert c_{n}\otimes _{1}c_{n}\right\vert ^{2}=%
\frac{1}{4n^{2}\ln ^{2}n}\sum_{i\neq j}(a\otimes _{1}a)^{2}(i,j) \\
&\leq &\frac{2(\pi +2)}{\ln ^{2}n}+\frac{2}{n^{2}\ln ^{2}n}\sum_{i\neq j}\ln
^{2}\frac{\sqrt{1-t_{i}}+\sqrt{1-t_{j}}}{\left\vert \sqrt{t_{i}}-\sqrt{t_{j}}%
\right\vert }\leq \frac{C}{\ln ^{2}n}.
\end{eqnarray*}%
$\square $

\section{Stochastic calculus of variation under the Doeblin's condition}

\label{sect:doeblin}

We assume that the sequence $X=(X_{n})_{n\in {\mathbb{N}}},$ $%
X_{n}=(X_{n,1},\ldots ,X_{n,d_{\ast }})\in{\mathbb{R}}^{d_{\ast }}$, of
independent random variables satisfies Hypothesis ${\mathfrak{M}}%
(\varepsilon,r,R)$, that is the Doeblin's condition ${\mathfrak{D}}%
(\varepsilon, r, R)$ and the moment finiteness one. We strongly use here the
representation (\ref{repr}) discussed in Section \ref{sect:3.1}, that is,
\begin{equation*}
X_n=\chi_nV_n+(1-\chi_n)U_n,\quad n\in{\mathbb{N}},
\end{equation*}
where $\chi_n,V_n,U_n$ are independent with laws given in (\ref{laws}). The
goal of this section is to present a differential calculus based on $%
V_{n},n\in {\mathbb{N}}$ which has been introduced in \cite{[BC-CLT],[BCP]}
(and which is inspired by the Malliavin calculus \cite{bib:[N]}).

\subsection{Abstract Malliavin calculus and Sobolev spaces}

\label{sect:sobolev}

To begin we introduce the space of the simple functionals. We denote by $%
\Lambda _{m}$ the multi-indexes $ \alpha =(\alpha _{1},\ldots ,\alpha _{m})$
with $\alpha _{i}=(n_{i},j_{i})\in {\mathbb{N}}\times \lbrack d_{\ast }]$ (that is, we do not impose that $n_1<\cdots < n_m$).
We consider polynomials with random coefficients
\begin{equation*}
P_{N}(x)=\sum_{m=0}^{N}\sum_{\alpha \in \Lambda _{m}}d(\alpha )x^{\alpha }
\end{equation*}%
where $x=(x_{n})_{n\in \N}$ with $x_{n}=(x_{n,1},...,x_{n,d_{\ast }})\in
R^{d_{\ast }}$ and $x^{\alpha }=\prod_{i=1}^{m}x_{\alpha _{i}}.$ The
coefficients $d(\alpha )\in\mathcal{U}$ are random variables which are measurable with
respect to $\sigma (\chi _{n},U_{n},n\in {\mathbb{N)}}$ and so, in
particular, are independent of $(V_{n})_{n\in \N}.$ And we define $\mathcal{P}%
_{N}(\mathcal{U})$ to be the space of the polynomials computed in $%
x_{n}=V_{n}$ that is $F\in \mathcal{P}_{N}(\mathcal{U})$ if
\begin{equation*}
F=P_{N}(V)=\sum_{m=0}^{N}\sum_{\alpha \in \Lambda _{m}}d(\alpha )V^{\alpha }.
\end{equation*}%
The simple functionals will be $\mathcal{P}(\mathcal{U})=\cup _{N\in {%
\mathbb{N}}}\mathcal{P}_{N}(\mathcal{U}).$ In particular our polynomials $%
Q_{N,k_{\ast }}(c,X)$ belong to $\mathcal{P}_{N}(\mathcal{U}).$ Note that $%
\mathcal{P}(\mathcal{U})$ is dense in $L^{p}(\Omega ,\mathcal{F},P)$ with $%
\mathcal{F}=\sigma (X_{n},n\in {\mathbb{N}})$. So we will define first our
differential operators on $\mathcal{P}(\mathcal{U}),$ and we extend them in
the canonical way to their domains in $L^{p}(\Omega ,\mathcal{F},P)$.

We assume that $\mathcal{U}={\mathbb{R}}^{d}$ (so it is a finite
dimensional Hilbert space). Let $F\in \mathcal{P}(\mathcal{U})$, so $%
F=Q_{N,k_{\ast }}(c,X)$. For $n\in {\mathbb{N}}$ and $i\in \lbrack d_{\ast }]
$ we define the first order derivatives%
\begin{equation*}
D_{n,i}F=\chi _{n}\times \partial _{{n,i}}Q_{N,k_{\ast }}(c,X)=\frac{%
\partial F}{\partial {V_{n,i}}}.
\end{equation*}%
We look to $DF=(D_{n,i}F)_{n\in {\mathbb{N}},i\in \lbrack d_{\ast }]}$ as to
a random element of the following Hilbert space $\mathcal{H(U)}$:
\begin{equation}
\mathcal{H(\mathcal{U})}=\Big\{x\in \otimes _{n=1}^{\infty }\mathcal{U}%
^{d_{\ast }}:\left\vert x\right\vert _{\mathcal{H}}^{2}:=\sum_{n=1}^{\infty
}\sum_{i=1}^{d_{\ast }}\left\vert x_{n,i}\right\vert _{\mathcal{U}%
}^{2}<\infty \Big\}.  \label{HU}
\end{equation}%
So $D:\mathcal{P}_{N}(\mathcal{U})\rightarrow \mathcal{P}_{N-1}(\mathcal{H(%
\mathcal{U})}).$
The Malliavin covariance matrix of $F\in \mathcal{P}(\mathcal{\mathcal{U}}%
)^{d}$ is defined by%
\begin{equation}
\sigma _{F}^{i,j}=\left\langle DF^{i},DF^{j}\right\rangle _{\mathcal{H(U)}%
}=\sum_{n=1}^{\infty }\sum_{l=1}^{d_{\ast }}D_{n,l}F^{i}\times
D_{n,l}F^{j},\quad i,j=1,\ldots ,d.  \label{Cov}
\end{equation}

Moreover we define the higher order derivatives in the following way. Let $%
m\in {\mathbb{N}}$ be fixed and let $\alpha =(\alpha _{1},\ldots ,\alpha
_{m})$ with $\alpha _{i}=(n_{i},j_{i})\in {\mathbb{N}}\times \lbrack d_{\ast
}].$ For $F=Q_{N,k_{\ast }}(c,X)\in \mathcal{P(U)}$, we define
\begin{equation}
D_{\alpha }^{(m)}F=D_{\alpha _{m}}\cdots D_{\alpha _{1}}F=\Big(%
\prod_{j=1}^{m}\chi _{n_{j}}\Big)(\partial _{n_{m},j_{m}}\cdots \partial
_{n_{1},j_{1}}Q_{N})(c,X)=\Big(\prod_{j=1}^{m}\chi _{n_{j}}\Big)\partial
_{\alpha }Q_{N}(c,X).  \label{S4}
\end{equation}%
We look to $D^{(m)}F=(D_{\alpha }^{(m)}F)_{\alpha \in \Gamma _{m}}$ as to a
random element of $\mathcal{H}_{m}:=\mathcal{H}^{\otimes m}(\mathcal{U}),$
so $D^{(m)}:\mathcal{P}_{N}(\mathcal{U})\rightarrow \mathcal{P}_{N-m}(%
\mathcal{H}^{\otimes m}(\mathcal{U}))$. For $m=1$, we have $%
D^{(1)}F=DF$.

We define now the divergence operator%
\begin{eqnarray}
LF &=&-\sum_{n=1}^{\infty }\sum_{i=1}^{d_{\ast
}}(D_{n,i}D_{n,i}F+D_{n,i}F\times \Theta _{n,i})\qquad \mbox{with}\qquad
\label{S6} \\
\Theta _{n,i} &=&2\chi _{n}\theta _{r}^{\prime }(\left\vert
X_{n,i}-x_{n,i}\right\vert ^{2})(X_{n,i}-x_{n,i}).  \label{S6'}
\end{eqnarray}%
Standard integration by parts on ${\mathbb{R}}$ gives the following duality
relation: for every $F,G\in \mathcal{P(U)}$%
\begin{equation}
{\mathbb{E}}(\left\langle DF,DG\right\rangle _{\mathcal{H(U)}})={\mathbb{E}}%
(\left\langle F,LG\right\rangle _{\mathcal{U}})={\mathbb{E}}(\left\langle
G,LF\right\rangle _{\mathcal{U}}).  \label{S7}
\end{equation}

We define now the Sobolev norms. For $q\geq 1$ we set
\begin{equation}
\left\vert F\right\vert _{1,q,\mathcal{U}}=\sum_{n=1}^{q}|D^{(n)}F|_{%
\mathcal{H}^{\otimes n}(\mathcal{U})}\quad \mbox{and}\quad \left\vert
F\right\vert _{q,\mathcal{U}}=\left\vert F\right\vert +\left\vert
F\right\vert _{1,q,\mathcal{U}}.  \label{S8}
\end{equation}%
Moreover we define%
\begin{equation}
\left\Vert F\right\Vert _{1,q,p,\mathcal{U}}=\big({\mathbb{E}}(\left\vert
F\right\vert _{1,q,\mathcal{U}}^{p})\big)^{1/p},\qquad \left\Vert
F\right\Vert _{q,p,\mathcal{U}}=\big({\mathbb{E}}(\left\vert F\right\vert
_{q,\mathcal{U}}^{p})\big)^{1/p}  \label{S9}
\end{equation}%
and%
\begin{equation}
\left\Vert \left\vert F\right\vert \right\Vert _{1,q,p,\mathcal{U}%
}=\left\Vert F\right\Vert _{1,q,p,\mathcal{U}}+\left\Vert LF\right\Vert
_{q-2,p,\mathcal{U}},\qquad \left\Vert \left\vert F\right\vert \right\Vert
_{q,p,\mathcal{U}}=\left\Vert F\right\Vert _{p,\mathcal{U}}+\left\Vert
\left\vert F\right\vert \right\Vert _{1,q,p,\mathcal{U}}.  \label{S3}
\end{equation}

Finally we define the Sobolev spaces%
\begin{equation}
{\mathbb{D}}^{q,p}=\overline{\mathcal{P}}^{\Vert \left\vert \cdot
\right\vert \Vert _{q,p,\mathcal{U}}}(\mathcal{U)},\qquad {\mathbb{D}}%
^{q,\infty }=\cap _{p=1}^{\infty }{\mathbb{D}}^{q,p}\qquad {\mathbb{D}}%
^{\infty }=\cap _{q=1}^{\infty }{\mathbb{D}}^{q,\infty }.  \label{S10}
\end{equation}%
The duality relation (\ref{S7}) implies that the operators $D^{(n)}$ and $L$
are closable so we may extend these operators to ${\mathbb{D}}^{q,p}$ in a
standard way. But in this work we will restrict ourself to $\mathcal{P(U)}$.

We recall now the basic computational rules. For $\phi \in C_{\mathrm{%
{\scriptsize {pol}}}}^{1}({\mathbb{{\mathbb{R}}}}^{M})$ and $F\in \mathcal{%
P(U)}^{M}$ we have%
\begin{equation}
D\phi (F)=\sum_{j=1}^{M}\partial _{j}\phi (F)DF^{j},  \label{FD9}
\end{equation}%
and for $\phi \in C_{\mathrm{{\scriptsize {pol}}}}^{2}({\mathbb{R}}^{M})$
\begin{equation}
L\phi (F)=\sum_{j=1}^{M}\partial _{j}\phi (F)LF^{j}-\frac{1}{2}%
\sum_{i,j=1}^{M}\partial _{i}\partial _{j}\phi (F)\left\langle
DF^{i},DF^{j}\right\rangle _{\mathcal{H}}.  \label{FD11}
\end{equation}%
In particular for $F,G\in {\mathbb{D}}^{2,\infty }$
\begin{equation}
L(FG)=FLG+GLF-\left\langle DF,DG\right\rangle _{\mathcal{H}}.  \label{FD12}
\end{equation}

Let us stress the following fact which is specific in our framework. In
order to establish the integration by parts formula in the classical
Malliavin calculus one needs that $\sigma _{F}$ is almost surely invertible.
And this is always falls here: indeed if $F=\phi (X_{1},\ldots ,X_{n})$ then
$DF=0$ on the set $\{\chi _{1}=\ldots =\chi _{n}=0\}$ which has strictly
positive probability. This is why we have to use a localized version of the
integration by parts formula. Given $\eta >0$ we consider a function $\Phi
_{\eta }:{\mathbb{R}}\rightarrow {\mathbb{R}}_{+}$ such that $%
1_{\{\left\vert x\right\vert \leq \eta \}}\leq \Phi _{\eta }(x)\leq
1_{\{\left\vert x\right\vert \leq 2\eta \}}$ and $|\Phi _{\eta
}^{(k)}(x)|\leq C_{k}\eta ^{-k}$ for every $k\in {\mathbb{N}}.$ Then we
define $\Psi _{\eta }=1-\Phi _{\eta }$ and we notice that on the set $\{\Psi
_{\eta }(\det \sigma _{F})>0\}$ we have $\det \sigma _{F}\geq \eta $, so $%
\sigma _{F}$ is invertible. We denote%
\begin{equation*}
\gamma _{F,\eta }=1_{\{\Psi _{\eta }(\det \sigma _{F})>0\}}\sigma _{F}^{-1}.
\end{equation*}

\begin{theorem}
\label{TH1}Let $F=(F^{1},\ldots ,F^{d}),F_{i}\in {\mathbb{D}}^{2,\infty }$
and $G\in {\mathbb{D}}^{1,\infty }$ and, for $\eta >0,$ we denote $G_{\eta
}=G\times \Psi _{\eta }(\det \sigma _{F}).$ Then for every $\phi \in
C_{p}^{\infty }({\mathbb{R}}^{d})$ and every $i=1,\ldots ,d$%
\begin{equation}
{\mathbb{E}}(\partial _{i}\phi (F)G_{\eta })={\mathbb{E}}(\phi (F)H_{\eta
,i}(F,G))  \label{FD13}
\end{equation}%
with
\begin{equation}
H_{\eta ,i}(F,G_{\eta })=G\gamma _{F,\eta }LF+\left\langle D(G_{\eta }\gamma
_{F,\eta }),DF\right\rangle _{\mathcal{H}}  \label{FD14}
\end{equation}%
Moreover let $m\in {\mathbb{N}},m\geq 2$ and $\alpha =(\alpha _{1},\ldots
,\alpha _{m})\in \{1,\ldots ,d\}^{m}.$ Suppose that $F=(F^{1},\ldots
,F^{d}),F_{i}\in {\mathbb{D}}^{m+1,\infty }$ and $G\in {\mathbb{D}}%
^{m,\infty }.$ Then%
\begin{equation}
{\mathbb{E}}(\partial _{\alpha }\phi (F)G_{\eta })={\mathbb{E}}(\phi
(F)H_{\eta ,\alpha }(F,G))  \label{FD15}
\end{equation}%
with $H_{\eta ,\alpha }(F,G)$ defined by $H_{\eta ,(\alpha _{1},\ldots
,\alpha _{m})}(F,G):=H_{\eta ,\alpha _{m}}(F,H_{\eta ,(\alpha _{1},\ldots
,\alpha _{m-1})}(F,G)).$
\end{theorem}

\textbf{Proof. }The proof is standard so we just sketch it\textbf{. }Using
the chain rule $D\phi (F)=\nabla \phi (F)DF$ so that
\begin{equation*}
\left\langle D\phi (F),DF\right\rangle _{\mathcal{H}}=\nabla \phi
(F)\left\langle DF,DF\right\rangle _{\mathcal{H}}=\nabla \phi (F)\sigma _{F}.
\end{equation*}%
It follows that, on the set $\{\Phi _{\eta }(\det \sigma _{F})>0\},$ one has
$\nabla \phi (F)=\gamma _{F_{\eta }}\left\langle D\phi (F),DF\right\rangle _{%
\mathcal{H}}$. Then, by using (\ref{FD12}) and the duality formula (\ref{S7}%
),
\begin{eqnarray*}
{\mathbb{E}}(G_{\eta }\nabla \phi (F)) &=&{\mathbb{E}}(G_{\eta }\gamma
_{F,\eta }\left\langle D\phi (F),DF\right\rangle _{\mathcal{H}})={\mathbb{E}}%
(G_{\eta }\gamma _{F,\eta }(L(\phi (F)F)-\phi (F)LF+FL\phi (F)) \\
&=&{\mathbb{E}}(\phi (F)(FL(G_{\eta }\gamma _{F,\eta F})+G_{\eta }\gamma
_{F,\eta }LF+L(G_{\eta }\gamma _{F,\eta }F)).
\end{eqnarray*}%
We use once again (\ref{FD12}) in order to obtain $H_{\eta ,i}(F,G)$ in (\ref%
{FD14}). By iteration one obtains the higher order integration by parts
formulae. $\square $

\medskip

We give now useful estimates for the weights which appear in (\ref{FD15}).
For $n,k\in {\mathbb{N}}$ we denote%
\begin{equation}
{\mathcal{K}}_{n,k}(F)=(\left\vert F\right\vert _{1,k+n+1}+\left\vert
LF\right\vert _{k+n})^{n}(1+\left\vert F\right\vert _{1,k+n+1})^{2d(2n+k)}.
\label{FD15'}
\end{equation}

\begin{lemma}
\label{L1} Let $n,k\in {\mathbb{N}}$  and $F\in \mathcal{P}^{d}$ and $G\in
\mathcal{P}.$ There exists a universal constant $C\geq 1$ (depending on $%
d,n,k$ only) such that for every multi index $\alpha $ with $\left\vert
\alpha \right\vert =n$ and every $\eta >0$ one has%
\begin{equation}
\big\vert H_{\alpha}(F,\Psi _{\eta }(\det \sigma _{F})G)\big\vert _{k}\leq
\frac{C}{\eta ^{2n+k}}\times {\mathcal{K}}_{n,k}(F)\times \left\vert
G\right\vert _{k+n}.  \label{FD16}
\end{equation}%
In particular, taking $k=0$ and $G=1$ we have%
\begin{equation}
\left\Vert H_{\alpha }(F,\Psi _{\eta }(\det \sigma _{F}))\right\Vert
_{p}\leq \frac{C}{\eta ^{2n}}\times \left\Vert {\mathcal{K}}%
_{n,0}(F)\right\Vert _{p}  \label{FD4}
\end{equation}
\end{lemma}

The proof is straightforward but technical so we leave it for Appendix \ref%
{app:norms}.

\subsection{Regularization results}

\label{sect:reg}

We deal here with functions and their derivatives on ${\mathbb{R}}^d$. So,
we use a slightly different definition for multi-indexes. Here, for $m\in {%
\mathbb{N}}$, a multi-index of length $m$ is given by $\alpha\in\{1,\ldots,d%
\}^m$ and we set $|\alpha|=m$ its length. For $y=(y_1,\ldots,y_d)\in{\mathbb{%
R}}^d$, we set $y^\alpha=\prod_{i=1}^d y_{\alpha_i}$. We allow the case $%
\alpha=\emptyset$ by setting $|\alpha|=0$ and, for $y\in{\mathbb{R}}^d$, $%
y^\alpha=1$.

We recall that a super kernel $\phi :{\mathbb{R}}^{d}\rightarrow {\mathbb{R}}
$ is a function which belongs to the Schwartz space $\mathbb{S}({\mathbb{R}}%
^d)$ (infinitely differentiable functions which decrease in a polynomial way
to infinity), $\int \phi (x)dx=1,$ and such that for every multi-index $%
\alpha$ with $|\alpha|=m$ one has
\begin{eqnarray}
\int y^{\alpha }\phi (y)dy &=&0 \quad \mbox{and}  \label{kk1} \\
\int \left\vert y\right\vert ^{m}\left\vert \phi (y)\right\vert dy &<&\infty
\quad \forall m\geq 1.  \label{kk2}
\end{eqnarray}%
For $\delta \in (0,1)$ we define $\phi _{\delta }(y)=\delta ^{-d}\phi
(\delta ^{-1}y)$ and for a function $f:{\mathbb{R}}^{d}\rightarrow {\mathbb{R%
}}$ we denote $f_{\delta }=f\ast \phi _{\delta }$, the symbol $\ast$
denoting convolution. For $f\in C_{{\mathrm{{\scriptsize {pol}}}}}^{k}({%
\mathbb{R}}^{d})$ we define $L_{k}(f)$ and $l_{k}(f)$ to be some constants
such that
\begin{equation*}
\sum_{0\leq \left\vert \alpha \right\vert \leq k}\left\vert \partial
^{\alpha }f(x)\right\vert \leq L_{k}(f)(1+\left\vert x\right\vert
)^{l_{k}(f)}.
\end{equation*}%
We give now a ``regularization lemma'' which is an improvement of Lemma 2.5
in \cite{[BC-EJP]}.

\begin{lemma}
\label{R-1} Let $F\in \mathcal{P}({\mathbb{R}})^{d}$ and $q,m\in {\mathbb{N}}%
.$ There exists some constant $C\geq 1,$ depending on $d,m$ and $q $ only,
such that for every $f\in C_{{\mathrm{{\scriptsize {pol}}}}}^{q+m}({\mathbb{R%
}}^{d}),$ every multi index $\gamma $ with $\left\vert \gamma \right\vert =m$%
 and every $\eta ,\delta >0 $%
\begin{eqnarray}
&&\left\vert {\mathbb{E}}(\Psi _{\eta }(\det \sigma _{F})\partial ^{\gamma
}f(F))-{\mathbb{E}}(\Psi _{\eta }(\det \sigma _{F})\partial ^{\gamma
}f_{\delta }(F))\right\vert  \label{1a} \\
&\leq &C\,2^{l_0(f)-1}c_{l_0(f),q} L_{0}(f)\left\Vert F\right\Vert
_{2l_{0}(f)}^{l_{0}(f)}\left\Vert {\mathcal{K}}_{q+m,0}(F)\right\Vert _{2}%
\frac{\delta ^{q}}{\eta ^{2(q+m)}}  \notag
\end{eqnarray}%
with ${\mathcal{K}}_{q+m,0}(F)$ defined in (\ref{FD15'}) and $c_{l,q}=\int
|\phi(y)||y|^q(1+|y|)^{l}dy$.
Moreover, for every $p>1$%
\begin{eqnarray}
&&\left\vert {\mathbb{E}}(\partial ^{\gamma }f(F))-{\mathbb{E}}(\partial
^{\gamma }f_{\delta }(F))\right\vert  \label{1b-1} \\
&\leq &C\left\Vert F\right\Vert _{pl_{0}(f)}^{l_{0}(f)}\Big(%
L_{m}(f)c_{l_m(f),0}{\mathbb{P}}^{(p-1)/p}(\det \sigma _{F}\leq \eta
)+2^{l_0(f)-1}c_{l_0(f),q}\,L_{0}(f)\frac{\delta ^{q}}{\eta ^{2(q+m)}}%
\left\Vert {\mathcal{K}}_{q+m,0}(F)\right\Vert _{2}\Big).  \notag
\end{eqnarray}
\end{lemma}

\textbf{Proof.} Using Taylor expansion of order $q$,
\begin{align*}
\partial ^{\gamma }f(x)-\partial ^{\gamma }f_{\delta }(x)& =\int (\partial
^{\gamma }f(x)-\partial ^{\gamma }f(y))\phi _{\delta }(x-y)dy \\
& =\int I(x,y)\phi _{\delta }(x-y)dy+\int R(x,y)\phi _{\delta }(x-y)dy
\end{align*}%
with
\begin{align*}
I(x,y)& =\sum_{i=1}^{q-1}\frac{1}{i!}\sum_{\left\vert \alpha \right\vert
=i}\partial ^{\gamma }\partial ^{\alpha }f(x)(x-y)^{\alpha }, \\
R(x,y)& =\frac{1}{q!}\sum_{\left\vert \alpha \right\vert
=q}\int_{0}^{1}\partial ^{\gamma }\partial ^{\alpha }f(x+\lambda
(y-x))(x-y)^{\alpha }\lambda ^{q}d\lambda .
\end{align*}%
Using (\ref{kk1}) we obtain $\int I(x,y)\phi _{\delta }(x-y)dy=0$ and by a
change of variable we get
\begin{equation*}
\int R(x,y)\phi _{\delta }(x-y)dy=\frac{1}{q!}\sum_{\left\vert \alpha
\right\vert =q}\int_{0}^{1}\int dz\phi _{\delta }(z)\partial ^{\gamma
}\partial ^{\alpha }f(x+\lambda z)z^{\alpha }\lambda ^{q}d\lambda .
\end{equation*}%
So that%
\begin{align*}
& {\mathbb{E}}(\Psi _{\eta }(\det \sigma _{F})\partial ^{\gamma }f(F))-{%
\mathbb{E}}(\Psi _{\eta }(\det \sigma _{F})\partial ^{\gamma }f_{\delta }(F))
\\
& ={\mathbb{E}}(\int \Psi _{\eta }(\det \sigma _{F})R(F,y)\phi _{\delta
}(F-y)dy) \\
& =\frac{1}{q!}\sum_{\left\vert \alpha \right\vert =q}\int_{0}^{1}\int
dz\phi _{\delta }(z){\mathbb{E}}(\Psi _{\eta }(\det \sigma _{F})\partial
^{\gamma }\partial ^{\alpha }f(F+\lambda z))z^{\alpha }\lambda ^{q}d\lambda .
\end{align*}%
Using integration by parts formula (\ref{FD15}) (with $G=1)$\
\begin{align*}
& \left\vert {\mathbb{E}}(\Psi _{\eta }(\det \sigma _{F})\partial ^{\gamma
}\partial ^{\alpha }f(F+\lambda z))\right\vert =\left\vert {\mathbb{E}}%
(f(F+\lambda z)H_{(\alpha ,\gamma )}(F,\Psi _{\eta }(\det \sigma
_{F}))\right\vert  \\
& \leq L_{0}(f){\mathbb{E}}((1+\left\vert z\right\vert +\left\vert
F\right\vert )^{l_{0}(f)}\left\vert H_{(\alpha ,\gamma )}(F,\Psi _{\eta
}(\det \sigma _{F}))\right\vert ) \\
& \leq C2^{l_{0}(f)-1}(1+\left\vert z\right\vert
)^{l_{0}(f)}L_{0}(f)\left\Vert F\right\Vert _{2l_{0}(f)}^{l_{0}(f)}({\mathbb{%
E}}(\left\vert H_{(\alpha ,\gamma )}(F,\Psi _{\eta }(\det \sigma
_{F}))\right\vert ^{2}))^{1/2}.
\end{align*}%
The upper bound from (\ref{FD4}) (with $p=2)$ gives%
\begin{equation*}
({\mathbb{E}}(\left\vert H_{(\alpha ,\gamma )}(F,\Psi _{\eta }(\det \sigma
_{F}))\right\vert ^{2})^{1/2}\leq \frac{C}{\eta ^{2(q+m)}}\left\Vert {%
\mathcal{K}}_{q+m,0}(F)\right\Vert _{2}
\end{equation*}%
And since
\begin{equation*}
\int dz\left\vert \phi _{\delta }(z)z^{\alpha }\right\vert (1+\left\vert
z\right\vert )^{l_{0}(f)}d\leq \delta ^{q}\int \phi
(y)|y|^{q}(1+|y|)^{l_{0}(f)}dy=\delta ^{q}c_{l_{0}(f),q},
\end{equation*}
we conclude that

\begin{eqnarray*}
&&\left\vert {\mathbb{E}}(\Psi _{\eta }(\det \sigma _{F}))\partial ^{\gamma
}f(F))-{\mathbb{E}}(\Psi _{\eta }(\det \sigma _{F}))\partial ^{\gamma
}f_{\delta }(F))\right\vert \\
&\leq &C\,c_{l_0(f)}L_{0}(f)\left\Vert F\right\Vert
_{2l_{0}(f)}^{l_{0}(f)}\left\Vert {\mathcal{K}}_{q+m,0}(F)\right\Vert _{2}%
\frac{C\delta ^{q}}{\eta ^{2(q+m)}}.
\end{eqnarray*}

In order to prove (\ref{1b}), we write
\begin{eqnarray*}
&&\left\vert {\mathbb{E}}((1-\Psi _{\eta }(\det \sigma _{F})))\partial
^{\gamma }f(F))-{\mathbb{E}}((1-\Psi _{\eta }(\det \sigma _{F})))\partial
^{\gamma }f_{\delta }(F))\right\vert \\
&\leq &2(L_{0}(\partial ^{\gamma }f_{\delta })\vee L_{0}(\partial ^{\gamma
}f)){\mathbb{E}}((1-\Psi _{\eta }(\det \sigma
_{F}))^{p/(p-1)})^{(p-1)/p}(1+\left\vert F\right\vert )^{pl_{0}(\partial
^{\gamma }f_{\delta })\vee l_{0}(\partial ^{\gamma }f)}) \\
&\leq &2(L_{0}(\partial ^{\gamma }f_{\delta })\vee L_{0}(\partial ^{\gamma
}f))\left\Vert F\right\Vert _{2l_{0}(f_{\delta })\vee
l_{0}(f)}^{pl_{0}(f_{\delta })\vee l_{0}(f))}{\mathbb{P}}^{(p-1)/p}(\det
\sigma _{F}\leq \eta ).
\end{eqnarray*}%
So the proof of (\ref{1b}) will be completed as soon as we check that $%
l_{0}(\partial ^{\gamma }f_{\delta })=l_{0}(\partial ^{\gamma }f)\leq
l_{m}(f)$ and $L_{0}(\partial ^{\gamma }f_{\delta })\leq L_{0}(\partial
^{\gamma }f)c_{l_0(\partial^\gamma f),0}\leq L_{m}(f)c_{l_m(f),0}.$ We write
\begin{eqnarray*}
\left\vert \partial ^{\gamma }f_{\delta }(x)\right\vert &=&\left\vert \int
\partial ^{\gamma }f(x-y)\phi _{\delta }(y)dy\right\vert \leq L_{0}(\partial
^{\gamma }f)\int (1+\left\vert x-y\right\vert )^{l_{0}(\partial ^{\gamma
}f)}\left\vert \phi _{\delta }(y)\right\vert dy \\
&\leq &L_{0}(\partial ^{\gamma }f)(1+\left\vert x\right\vert
)^{l_{0}(\partial ^{\gamma }f)}\int (1+\left\vert y\right\vert
)^{l_{0}(\partial ^{\gamma }f)}\left\vert \phi _{\delta }(y)\right\vert
dy\leq L_{m}(f)(1+\left\vert x\right\vert )^{l_{m}(f)} c_{l_m(f),0}.
\end{eqnarray*}%
$\square $

\medskip

As a consequence, we get a regularization result involving functions which
are just continuous and bounded.

\begin{lemma}
\label{R} Let $F\in \mathcal{P}({\mathbb{R}})^{d}$ and $q\in {\mathbb{N}}.$
There exists some constant $C\geq 1,$ depending on $d,m$ and $q $ only, such
that for every $f\in C_{b}({\mathbb{R}}^{d})$, every $\eta ,\delta >0 $ and $%
a<1$,
\begin{equation}  \label{1b}
\left\vert {\mathbb{E}}(f(F))-{\mathbb{E}}(f_{\delta }(F))\right\vert \leq
C\|f\|_\infty \Big({\mathbb{P}}^{a}(\det \sigma _{F}\leq \eta )+\frac{\delta
^{q}}{\eta ^{2q}}\left\Vert {\mathcal{K}}_{q,0}(F)\right\Vert _{2}\Big),
\end{equation}
with ${\mathcal{K}}_{q,0}(F)$ defined in (\ref{FD15'}).
\end{lemma}

\textbf{Proof.} Let $g$ denote the density of the standard $d$-dimensional
normal law and for $\varepsilon>0$, set $g_\varepsilon(x)=\frac
1{\varepsilon^d}g(\frac x\varepsilon)$. We notice that $f\ast g_\varepsilon$%
, $f_\delta\ast g_\varepsilon\in C^\infty_b({\mathbb{R}}^d)$. Moreover, $%
l_0(f\ast g_\varepsilon)=l_0(f_\delta\ast g_\varepsilon)=0$ and $L_0(f\ast
g_\varepsilon)=L_0(f_\delta\ast g_\varepsilon)=\|f\|_\infty$, for every $%
\varepsilon>0$. So, we can apply (\ref{1b-1}) with $|\gamma|=0$ and we
obtain
\begin{equation*}
\left\vert {\mathbb{E}}(f\ast g_\varepsilon(F))-{\mathbb{E}}(f_{\delta }\ast
g_\varepsilon(F))\right\vert \leq C\|f\|_\infty \Big({\mathbb{P}}%
^{(p-1)/p}(\det \sigma _{F}\leq \eta )+\frac{\delta ^{q}}{\eta ^{2q}}%
\left\Vert {\mathcal{K}}_{q,0}(F)\right\Vert _{2}\Big).
\end{equation*}
We now let $\varepsilon$ tend to 0 and obtain (\ref{1b}). $\square$

\subsection{Estimates of the Sobolev norms}

\label{sect:est}

Through this section we assume that $X$ verifies ${\mathfrak{M}}%
(\varepsilon,r,R)$ (that is (\ref{H1}) and ${\mathfrak{D}}(\varepsilon, r,
R))$ and we estimates the Sobolev norms of $Q_{N}(c,X)$ and of $%
LQ_{N}(c,X). $ We will give our estimates in terms of the norms $\mathcal{N}%
_{{\mathcal{U}},q}(c,M)$ defined in (\ref{H2a}).

\begin{proposition}
\label{M0} Let $p\geq 2$ and $N,q\in {\mathbb{N}}$ be given and let $%
\overline{M}_{p}=b_{p}M_{p}\sqrt{k_{\ast }d_{\ast }}$ with $%
M_{p}=M_{p}(Z(X)).$ Then
\begin{equation}
\left\Vert Q_{N,k_{\ast }}(c,X)\right\Vert _{\mathcal{U},q,p}\leq
(q+1)2^q(1+k_{\ast}^{3/2}(1+M_{k_{\ast }}))^{q}\mathcal{N}_{{\mathcal{U}}%
,q}(c,\overline{M}_{p}).  \label{BIS1}
\end{equation}
\end{proposition}

\begin{remark}
\label{M1} (\ref{BIS1}) says in particular that if $\lim_{N\to\infty}%
\mathcal{N}_{{\mathcal{U}},q}(c,\overline{M}_{p})<\infty $ (recall that $%
\mathcal{N}_{{\mathcal{U}},q}(c,\overline{M}_{p})$ is a sum up to $N$, see (%
\ref{H2a})) then the infinite series $Q_{\infty,k_{\ast }}(c,X)$ belongs to $%
{\mathbb{D}}^{q,p}.$ Let us compare this result with the corresponding one
for functionals on the Wiener space. We take $k_{\ast }=1,d_{\ast }=1, {%
\mathcal{U}}={\mathbb{R}}$ and $X_{n} $ to be standard normal distributed.
Then $\Phi _{m}(c,X)$ is a multiple integral of order $m$ associated to the
kernel $f_{c,m}$ which is constant on cubes and equal to the corresponding $%
c(\alpha ).$ So $Q_{\infty ,1}(c,X)=\sum_{m=0}^{\infty }c(\alpha )X^{\alpha
}=\sum_{m=0}^{\infty }J_{m}(f_{c,m})=\sum_{m=0}^{\infty }\frac{1}{m!}%
I_{m}(f_{c,m})$ where $J_{m}$ denotes the iterated stochastic integral and $%
I_{m}$ is the multiple stochastic integral. Note that $b_{2}=1$ and $M_{2}=1$
so $\overline{M}_{2}=1. $ So we have%
\begin{equation*}
\mathcal{N}_{q}^{2}(c,\overline{M}_{2})=\sum_{m=0}^{sI}m^{q}\left\vert
c\right\vert _{\mathcal{U},m}^{2}=\sum_{m=0}^{\infty }m^{q}\frac{1}{m!}%
\left\Vert f_{c,m}\right\Vert _{L^{2}({\mathbb{R}}_{+}^{m})}^{2}.
\end{equation*}%
It is known that $Q_{\infty ,1}(c,X)$ is $q$ time differentiable in $L^{2}$
in Malliavin sense if and only if the quantity in the right hand side is
finite. And this is the same in our framework. But in our calculus we need
estimates for a large $p>2$ and then $\overline{M}_{p}>1.$ This is why we
give up in this paper the case of infinite series and we restrict ourself to
finite sums.
\end{remark}

\textbf{Proof.} \textbf{Step 1}. For simplicity of notation, we set here $%
Z=Z(X)$. For fixed $n_0\in{\mathbb{N}}$, $j_0\in[d_*]$ and $m\in{\mathbb{N}}$
we set $\Lambda_{n_0,j_0}(m,k)$ as the set of the multi-indexes of length $m$
which do not contain the pair $(n_0,kd_*+j_0)$, the case $m=0$ giving the
set $\Lambda_{n_0,j_0}(0,k)$ made just by the null multi-index. Then, by
observing that $\chi_nV_n^k=\chi_nX_n^k$ for every $n$ and $k$, one has
\begin{equation*}
D_{n_{0},j_{0}}S_N(c,Z)=D_{n_{0},j_{0}}\sum_{m=0}^{N}\sum_{|\alpha|=m
}c(\alpha )Z^{\alpha
}=\sum_{m=0}^{N-1}\sum_{k=0}^{k_*-1}\sum_{\beta\in%
\Lambda_{n_0,j_0}(m,k)}(Dc)_{n_{0},j_{0},k}(\beta )\chi
_{n_{0}}V_{n_0,j_0}^kZ^{\beta }
\end{equation*}
where $(Dc)_{n_{0},j_{0},k}(\beta)=c((n_0,kd_*+j_0))$ if $|\beta|=0$ and for
$|\beta|=m\geq 1$,
\begin{equation*}
\begin{array}{rl}
\displaystyle (Dc)_{n_{0},j_{0},k}(\beta) = & \displaystyle %
\sum_{i=1}^{m-1}c(\beta _{1},\ldots ,\beta _{i},(n_{0},kd_*+j_{0}),\beta
_{i+1},\ldots ,\beta _{m})1_{\{\beta _{i}^{\prime }<n_{0}<\beta
_{i+1}^{\prime }\}} \smallskip \\
& \displaystyle +c((n_{0},kd_*+j_{0}),\beta _{1},\ldots ,\beta
_{m})1_{\{n_{0}<\beta _{1}^{\prime }\}}+c(\beta _{1},\ldots ,\beta
_{m},(n_{0},kd_*+j_{0}))1_{\{n_{0}>\beta _{m}^{\prime }\}}.%
\end{array}%
\end{equation*}
It can be easily checked that
\begin{equation}  \label{Dnj-1}
D_{n_{0},j_{0}}S_N(c,Z) =\chi_{n_0}S_{N}((Tc)_{n_{0},j_{0}},Z).
\end{equation}%
%
%
%
%
where, for $|\beta|=m=0,1,\ldots,N$,
\begin{equation}  \label{Tc}
(Tc)_{n_{0},j_{0}}(\beta )= \overline{c}^{n_{0},j_{0}}(\beta)1_{m=0}+
1_{m\geq 1}\sum_{i=1}^{m}\sum_{k=1}^{k_{\ast }-1}(k+1)(\widehat{c}%
_{k,i}^{n_{0},j_{0}}(\beta ) +d_{k,i}^{n_{0},j_{0}}(\beta )1_{m\leq N-1})
\end{equation}%
and the above coefficients are
\begin{equation}  \label{c}
\begin{array}{rl}
\overline{c}^{n_{0},j_{0}}(\emptyset ) = & \displaystyle c((n_0,j_0))
+\sum_{k=1}^{k_*-1}(k+1)c((n_0,kd_\ast+j_0)){\mathbb{E}}(X_{n_0,j_0}^k)%
\smallskip \\
\displaystyle \widehat{c}_{k,i}^{n_{0},j_{0}}(\beta )= & c(\beta _{1},\ldots
,\beta _{i-1},(n_{0},kd_{\ast }+j_{0}),\beta _{i+1},\ldots ,\beta
_{m})1_{\{\beta _{i}=(n_{0},(k-1)d_{\ast }+j_{0})\}},\smallskip \\
d_{k,i}^{n_{0},j_{0}}(\beta ) = & \displaystyle {\mathbb{E}}%
(X_{n_{0},j_{0}}^{k})c((\beta _{1},\ldots ,\beta _{i-1},(n_{0},kd_{\ast
}+j_{0}),\beta _{i+1},\ldots ,\beta _{m}))%
\end{array}%
\end{equation}%
We study $\mathcal{N}_{\mathcal{H(U)},q}(Tc,M)$. First,
\begin{equation*}
|\overline{c}|^2_{\mathcal{H(U)},m} =|\overline{c}|^2_{\mathcal{H(U)},0}\leq
k_*^3M_{k_*}^2|c|^2_{\mathcal{U},1}.
\end{equation*}
Moreover, for $m\geq 1$,
\begin{eqnarray*}
\Big\vert \sum_{i=1}^{m}\sum_{k=1}^{k_{\ast }-1}(k+1)\widehat{c}_{k,i}%
\Big\vert _{\mathcal{H(U)},m}^{2} &\leq &mk_{\ast
}^{3}\sum_{i=1}^{m}\sum_{k=1}^{k_{\ast }-1}\left\vert \widehat{c}%
_{k,i}\right\vert _{\mathcal{H(U)},m}^{2} \\
&=&mk_{\ast }^{3}\sum_{i=1}^{m}\sum_{k=1}^{k_{\ast
}-1}\sum_{n_{0},j_{0}}\sum_{\left\vert \beta \right\vert =m}\left\vert
\widehat{c}_{k,i}^{n_{0},j_{0}}(\beta )\right\vert _{\mathcal{U}}^{2}\leq
mk_{\ast }^{3}\left\vert c\right\vert _{\mathcal{U},m}^{2}
\end{eqnarray*}%
and similarly,
\begin{equation*}
\Big\vert \sum_{i=1}^{m}\sum_{k=1}^{k_{\ast }-1}(k+1)d_{k,i}\Big\vert _{%
\mathcal{H(U)},m}^{2}\leq mk_{\ast }^{3}M_{k_{\ast }}^{2}\left\vert
c\right\vert _{\mathcal{U},m+1}^{2}
\end{equation*}%
We put all this together and we obtain
\begin{equation*}
\mathcal{N}_{\mathcal{H(U)},q}(Tc,M)\leq 2(1+k_{\ast }^{3/2}(1+M_{k_{\ast
}}))\mathcal{N}_{{\mathcal{U}},q+1}(c,M).
\end{equation*}

\textbf{Step 2}. Starting from formula (\ref{Dnj-1}), we use Burkholder's
inequality (\ref{H8}) in order to obtain
\begin{equation*}
\left\Vert D^{q}Q(c,X)\right\Vert _{\mathcal{H}^{\otimes q}\mathcal{(U)}%
,p}\leq\left\Vert Q(T^{q}c,X)\right\Vert _{\mathcal{H}^{\otimes q}\mathcal{%
(U)},p}\leq \mathcal{N}_{\mathcal{H}^{\otimes q},0}(T^{q}c,\overline{M}_{p})
\leq 2^q(1+k_{\ast}^{3/2}(1+M_{k_{\ast }}))^{q}\mathcal{N}_{{\mathcal{U}}%
,q}(c,\overline{M}_{p}).
\end{equation*}

$\square $

\medskip

In order to treat $LQ_{N,k_{\ast }}(c,X)$ we need the following auxiliary
lemma:

\begin{proposition}
\textbf{A}. Let $B_{n},\Lambda _{n}\in \mathcal{U}$ be random variables
such that $B_{n},\Lambda _{n}\in \mathcal{P(U)}$ for every $n$ and $B_{n}$
is $\sigma (X_{1},\ldots ,X_{n})$ measurable. We fix $j\in \lbrack d_{\ast
}],k\in \lbrack k_{\ast }]$ and we consider the process
\begin{equation}
Y_{J}=\sum_{n=1}^{J}B_{n-1}LX_{n,j}^{k}+\Lambda _{J}.  \label{bis1}
\end{equation}%
For every $q\in {\mathbb{N}}$ and $p\geq 2$ there exists a universal
constant $C\geq 1$ depending on $k_{\ast }$ and on $p$ only, such that
\begin{equation}
\max_{n\leq J}\left\Vert Y_{n}\right\Vert _{\mathcal{U},q,p}\leq q\,\frac{(C%
\widehat{M}_{p})^{q+1}}{r^{q+1}}\times {\mathcal{K}}_{q,p}(B,\Lambda )
\label{bis2}
\end{equation}%
with
\begin{align}
\widehat{M}_{p}&=b_{p}M_{2k_{\ast }p}^{k_{\ast }}(X)\sqrt{k_{\ast }d_{\ast }}
\label{bis2'} \\
{\mathcal{K}}_{q,p}(B,\Lambda )&=\Big(\sum_{k=1}^{J}\left\Vert
B_{k}\right\Vert _{\mathcal{U},q,p}^{2}\Big)^{1/2}+\max_{m\leq J}\left\Vert
\Lambda _{m}\right\Vert _{\mathcal{U},q,p}.  \label{bis3}
\end{align}%
\textbf{B}. If
\begin{equation*}
U_{J}=\sum_{n=1}^{J}B_{n-1}(X_{n,j}^{k}-{\mathbb{E}}(X_{n,j}^{k}))+\Lambda
_{J}.
\end{equation*}%
then
\begin{equation}
\max_{n\leq J}\left\Vert U_{n}\right\Vert _{\mathcal{U},q,p}\leq q\,\frac{(C%
\widehat{M}_{p})^{q+1}}{r^{q+1}}\times {\mathcal{K}}_{q,p}(B,\Lambda )
\label{bis5}
\end{equation}
\end{proposition}

\textbf{Proof. } In the following $C\geq 1$ denotes a constant depending on $%
k_{\ast }$ and on $p$ only and which may change from a line to another.

\textbf{Step 1.} We will use the following facts. First, by the duality
formula ${\mathbb{E}}(LX_{n,j}^{k})={\mathbb{E}}(\langle
DX_{n,j}^{k},D1\rangle )=0.$ Moreover using the computational rules (see (%
\ref{FD11}))%
\begin{equation*}
LX_{n,j}^{k}=kX_{n,j}^{k-1}LX_{n,j}+2k(k-1)X_{n,j}^{k-2}\left\langle
DX_{n,j},DX_{n,j}\right\rangle
=kX_{n,j}^{k-1}LX_{n,j}+2k(k-1)X_{n,j}^{k-2}\chi_n.
\end{equation*}%
It follows that%
\begin{equation*}
\Vert LX_{n,j}^{k}\Vert _{q,p}\leq k\Vert X_{n,j}^{k-1}\Vert _{q,2p}\Vert
LX_{n,j}\Vert _{q,2p}+2k(k-1)\Vert X_{n,j}^{k-2}\Vert _{q,p}
\end{equation*}%
It is easy to check that $\Vert X_{n,j}^{k-1}\Vert _{q,2p}\leq
(k-1)!M_{2k_{\ast }p}^{k_{\ast }}(X)$ and a similar estimates holds for $%
\Vert X_{n,j}^{k-2}\Vert _{q,2p}.$ Moreover it is proved in Lemma 3.2 in
\cite{[BC-CLT]} that there exists a universal constant $C$ such that$%
\left\Vert LX_{n,j}\right\Vert _{q,2p}\leq \frac{C}{r^{q+1}}$ so that
\begin{equation}
\Vert LX_{n,j}^{k}\Vert _{q,2p}\leq \frac{C}{r^{q+1}}M_{2k_{\ast
}p}^{k_{\ast }}(X).  \label{bis6}
\end{equation}

\textbf{Step 2}. Let $q=0,$ so that $\left\Vert Y_{J}\right\Vert _{\mathcal{U%
},q,p}=\left\Vert Y_{J}\right\Vert _{\mathcal{U},p}.$ We have to check that%
\begin{equation}
\max_{n\leq J}\left\Vert Y_{n}\right\Vert _{\mathcal{U},p}\leq \frac{C}{r}%
\times {\mathcal{K}}_{0,p}(B,\Lambda ).  \label{bis7}
\end{equation}%
Since $B_{n-1}$ is $\sigma (X_{1},\ldots ,X_{n-1})$ measurable and ${\mathbb{%
E}}(LX_{n,j}^{k})=0,$ it follows that $M_{m}=%
\sum_{n=1}^{m}B_{n-1}LX_{n,j}^{k}$ is a martingale. By (\ref{H5})%
\begin{equation*}
\left\Vert M_{m}\right\Vert _{\mathcal{U},p}\leq b_{p}\Big(%
\sum_{n=1}^{m}\Vert LX_{n,j}^{k}B_{n-1}\Vert _{\mathcal{U},p}^{2}\Big)^{1/2}.
\end{equation*}%
Since $LX_{n,j}^{k}$ and $B_{n-1}$ are independent,
\begin{equation*}
\Vert LX_{n,j}^{k}B_{n-1}\Vert _{\mathcal{U},p}^{2}=\Vert LX_{n,j}^{k}\Vert
_{p}^{2}\Vert B_{n-1}\Vert _{\mathcal{U},p}^{2}\leq \frac{C\widehat{M}%
_{p}^{2}}{r^{2}}\left\Vert B_{n-1}\right\Vert _{\mathcal{U},p}^{2}.
\end{equation*}%
From $Y_{m}=M_{m}+\Lambda _{m}$, we conclude that%
\begin{equation*}
\Vert Y_{m}\Vert _{\mathcal{U},p}\leq \Vert M_{m}\Vert _{\mathcal{U}%
,p}+\Vert \Lambda _{m}\Vert _{\mathcal{U},p}\leq \frac{1}{r}C\widehat{M}_{p}%
\Big(\Big(\sum_{k=1}^{m}\Vert B_{k}\Vert _{\mathcal{U},p}^{2}\Big)%
^{1/2}+\Vert \Lambda _{m}\Vert _{\mathcal{U},p}\Big)
\end{equation*}%
so the statement holds for $q=0$.

\textbf{Step 3}. We estimate the derivatives of $Y_{m}$. We have%
\begin{equation*}
\overline{Y}_{m}:=DY_{m}=\sum_{n=1}^{m}\overline{B}_{n-1}LX_{n,j}^{k}+%
\overline{\Lambda }_{m}.
\end{equation*}%
where $\overline{B}_{n}=DB_{n}$ is $\sigma(X_1,\ldots,X_{n-1})$-measurable
and $\overline{\Lambda }_{m}=\sum_{k=1}^{m}DLX_{n,j}^{k}B_{n-1}+D\Lambda
_{m}.$ Notice that $\overline{Y}_{m}$, $\overline{B}_{k}$ and $\overline{%
\Lambda }_{m}$ take values in $\mathcal{H}(\mathcal{U})$ (defined in (\ref%
{HU})). So, by applying the step above, we get
\begin{equation*}
\max_{n\leq J}\Vert DY_{n}\Vert _{\mathcal{H}(\mathcal{U}),p}\leq \frac{C%
\widehat{M}_{p}}{r}{\mathcal{K}}_{0,p}(\overline{B},\overline{\Lambda }),
\end{equation*}%
where
\begin{equation*}
{\mathcal{K}}_{0,p}(\overline{B},\overline{\Lambda })=\Big(%
\sum_{k=1}^{J}\left\Vert \overline{B}_{k}\right\Vert _{\mathcal{H}(\mathcal{U%
}),p}^{2}\Big)^{1/2}+\max_{m\leq J}\left\Vert \overline{\Lambda }%
_{m}\right\Vert _{\mathcal{H}(\mathcal{U}),p}.
\end{equation*}%
If we prove that
\begin{equation}
{\mathcal{K}}_{0,p}(\overline{B},\overline{\Lambda })\leq \frac{{\mathcal{K}}%
_{k_{\ast },p}\widehat{M}_{p}}{r}\times {\mathcal{K}}_{1,p}(B,\Lambda )
\label{bis8}
\end{equation}%
then we obtain
\begin{equation*}
\max_{m\leq J}\Vert Y_{m}\Vert _{\mathcal{U},1,p}\leq \frac{(C\widehat{M}%
_{p})^{2}}{r^{2}}{\mathcal{K}}_{1,p}(B,\Lambda ).
\end{equation*}%
And by iteration, we get (\ref{bis2}) for every $q$. So, let us prove (\ref%
{bis8}).

We have $\Vert \overline{B}_{k}\Vert _{\mathcal{H}(\mathcal{U}),p}=\Vert
DB_{k}\Vert _{\mathcal{H}(\mathcal{U}),p}\leq \left\Vert B_{k}\right\Vert _{%
\mathcal{U},1,p}$. We analyze now $\overline{\Lambda }_{m}.$ First, $\Vert
D\Lambda _{m}\Vert _{\mathcal{H}(\mathcal{U}),p}\leq \Vert \Lambda _{m}\Vert
_{\mathcal{U},1,p}$. Let $I_{m}:=\sum_{n=1}^{m}DLX_{n,j}^{k}B_{n-1}\in
\mathcal{H}(\mathcal{U})$. Since $D_{n^{\prime },j^{\prime }}LX_{n,j}^{k}=0$
if $(n^{\prime },j^{\prime })\neq (n,j)$ we obtain%
\begin{equation*}
\left\vert I_{m}\right\vert _{\mathcal{H}(\mathcal{U})}^{2}=\sum_{n=1}^{m}%
\vert D_{n,j}LX_{n,j}^{k}\vert ^{2}\left\vert B_{n-1}\right\vert _{\mathcal{U%
}}^{2}.
\end{equation*}%
Recalling that $D_{n,j}X_{n,j}^{k}$ and $B_{n-1}$ are independent and that $%
\Vert D_{n,j}LX_{n,j}^{k}\Vert _{p}^{2}\leq Cr^{-2}M_{2k_{\ast }p}^{k_{\ast
}}(X)$, we can write
\begin{align*}
\Vert I_{m}\Vert _{\mathcal{H}(\mathcal{U}),p}& =\Vert |I_{m}|_{\mathcal{H}(%
\mathcal{U})}^{2}\Vert _{p/2}^{1/2}\leq \Big(\sum_{n=1}^{m}\big\Vert \vert
D_{n,j}LX_{n,j}^{k}\vert ^{2}\,\vert B_{n-1}\vert _{\mathcal{U}}^{2}%
\big\Vert _{p/2}\Big)^{1/2} \\
& =\Big(\sum_{n=1}^{m}\Vert D_{n,j}LX_{n,j}^{k}\Vert _{p}^{2}\Vert
B_{n-1}\Vert _{\mathcal{U},p}^{2}\Big)^{1/2}\leq \frac{{\mathcal{K}}%
_{k_{\ast },p}\widehat{M}_{p}}{r}\times \Big(\sum_{n=1}^{m}\Vert
B_{n-1}\Vert _{\mathcal{U},p}^{2}\Big)^{1/2}.
\end{align*}%
By inserting all these estimates, we get (\ref{bis8}). So \textbf{A} is
proved. The proof of \textbf{B} is just identical so we skip it. $\square $

\begin{proposition}
\label{normaLQ} For every $q,N\in {\mathbb{N}}$ and $p\geq 2$ there exists a
universal constant $C$ depending on $k_{\ast },q$ and $p$ only such that
\begin{equation}
\left\Vert LQ_{N}(c,X)\right\Vert _{{\mathcal{U}},q,p}\leq \frac{C\widehat{M}%
_{p}^{N(q+1)}}{r^{q+1}}\mathcal{N}_{{\mathcal{U}},q+1}(c,\overline{M}_p),
\label{bis13}
\end{equation}
where $\overline{M}_p=b_pM_p\sqrt{k_*d_*}$ and $\widehat{M}_{p}$ is given in
(\ref{bis2'}).
\end{proposition}

\textbf{Proof. }We prove this by recurrence on $N$. The case $N=1$ is
straightforward, so we suppose $N>1$. We recall that, if $\left\vert \beta
\right\vert =m$ then $ c^{n,j}(\beta )=1_{\beta _{m}^{\prime }<n}c(\beta
,(n,j))$ and we write%
\begin{equation*}
Q_{N,k_\ast}(c,X)=c(\emptyset )+\sum_{n=1}^{\infty }\sum_{j=1}^{m_{\ast
}}Z_{n,j}Q_{N-1,k_\ast}(c^{n,j},X),
\end{equation*}%
where $Z=Z(X)$. Since $\left\langle
DZ_{n,j},DQ_{N-1,k_\ast}(c^{n,j},X)\right\rangle _{\mathcal{H(U)}}=0$ we get
(see (\ref{FD11}))
\begin{equation*}
LQ_{N,k_\ast}(c,X)=\sum_{n=1}^{\infty }\sum_{j=1}^{m_{\ast
}}Q_{N-1,k_\ast}(c^{n,j},X)LZ_{n,j}+\sum_{n=1}^{\infty }\sum_{j=1}^{m_{\ast
}}Z_{n,j}LQ_{N-1}(c^{n,j},X).
\end{equation*}%
So we are in the framework of the previous lemma with $B_{n-1}=Q_{N-1,k_%
\ast}(c^{n,j},X)$ and
\begin{equation*}
\Lambda =\sum_{n=1}^{\infty }\sum_{j=1}^{d_{\ast
}}Z_{n,j}LQ_{N-1,k_\ast}(c^{n,j},X).
\end{equation*}%
Notice that%
\begin{equation*}
M\sum_{n=1}^{\infty }\sum_{j=1}^{m_{\ast }}\mathcal{N}_{{\mathcal{U}}%
,q}(c^{n,j},M)\leq \mathcal{N}_{{\mathcal{U}},q}(c,M).
\end{equation*}%
Then, using (\ref{bis5}) (recall that $C_{p}(X)\leq \widehat{M}_{p})$ and
the recurrence hypothesis%
\begin{eqnarray*}
\left\Vert \Lambda \right\Vert _{\mathcal{U},q,p}^{2} &\leq &\widehat{M}%
_{p}^{2(q+1)}\sum_{n=1}^{\infty }\sum_{j=1}^{m_{\ast }}\left\Vert
LQ_{N-1,k_\ast}(c^{n,j},X)\right\Vert _{\mathcal{U},q,p}^{2} \\
&\leq &\widehat{M}_{p}^{2(q+1)}\sum_{n=1}^{\infty }\sum_{j=1}^{m_{\ast }}%
\widehat{M}_{p}^{2(q+1)(N-1)}\times \frac{C}{r^{q+1}}\mathcal{N}_{{\mathcal{U%
}},q+1}^{2}(c^{n,j},\overline{M}_p) \\
&=&\widehat{M}_{p}^{2(q+1)N}\times \frac{C}{r^{q+1}}\mathcal{N}_{{\mathcal{U}%
},q+1}^{2}(c,\overline{M}_p).
\end{eqnarray*}%
Moreover, by the estimates of the Sobolev norms given in (\ref{BIS1}), and
the same computations as above
\begin{equation*}
\sum_{n=1}^{\infty }\sum_{j=1}^{m_{\ast }}\left\Vert
Q_{N-1.k_\ast}(c^{n,j},X)\right\Vert _{{\mathcal{U}},q,p}^{2}\leq
C\sum_{n=1}^{\infty }\sum_{j=1}^{m_{\ast }}\mathcal{N}_{{\mathcal{U}}%
,q+1}^{2}(c^{n,j},\overline{M}_{p})\leq C\mathcal{N}_{{\mathcal{U}}%
,q+1}^{2}(c,\overline{M}_{p}).
\end{equation*}%
$\square $

\begin{remark}
By using Proposition \ref{M0} and \ref{normaLQ}, we give here an upper
estimate of the $L^2$-norm of the constant ${\mathcal{K}}_{q,0}(Q_{N,k_{\ast
}}(c,X))$ defined in (\ref{FD15'}). This will be very useful in the sequel.
By using the H\"older inequality we easily get
\begin{align*}
&\|{\mathcal{K}}_{q,0}(Q_{N,k_{\ast }}(c,X))\|_2 & \leq \big(\|Q_{N,k_{\ast
}}(c,X)\|_{q+1,4q}+\|LQ_{N,k_{\ast }}(c,X)\|_{q,4q}\big)^q \times\big(%
1+\|Q_{N,k_{\ast }}(c,X)\|_{q+1,16q}\big)^{4q}.
\end{align*}
By applying the estimates (\ref{BIS1}) and (\ref{bis13}) we obtain
\begin{equation}  \label{est-Cq0}
\|{\mathcal{K}}_{q,0}(Q_{N,k_{\ast }}(c,X))\|_2 \leq C |c|^q(1+|c|)^{4q},
\end{equation}
$C>0$ denoting a constant depending on $q,N,k_\ast$ and the moment bound $%
M_p(X)$ for a suitable $p>1$ and independent of the coefficients $c$.
\end{remark}

\subsection{Estimates of the covariance matrix}

\label{sect:cov}

In this section we give estimates for the Malliavin covariance matrix of $%
Q_{N,k_{\ast }}(c,X)$ which we shortly denote by $\sigma _{N}$. We restrict
ourself to the scalar case, so that $Q_{N,k_{\ast }}(c,X)\in {\mathbb{R}}=%
\mathcal{U}$ and $\sigma _{N}$ is just a scalar. We start from the formula
of the Malliavin derivative of $Q_{N,k_{\ast }}(c,X)$ already discussed in
the proof of Proposition \ref{M0}, that is,
\begin{equation}
D_{n_{0},j_{0}}Q_{N,k_{\ast
}}(c,X)=D_{n_{0},j_{0}}S_{N}(c,Z(X))=\sum_{m=0}^{N-1}\sum_{k=0}^{k_{\ast
}-1}\sum_{\beta \in \Lambda
_{n_{0},j_{0}}(m,k)}(k+1)(Dc)_{n_{0},j_{0},k}(\beta )\chi
_{n_{0}}V_{n_{0},j_{0}}^{k}Z^{\beta }(X)  \label{Dnj}
\end{equation}%
where $\Lambda _{n_{0},j_{0}}(m,k)$ denotes the multi-indexes of length $m$
which do not contain the pair $(n_{0},kd_{\ast }+j_{0})$ and where $%
(Dc)_{n_{0},j_{0},k}(\beta )=c((n_{0},kd_{\ast }+j_{0}))$ if $|\beta |=0$
and for $|\beta |=m\geq 1$,
\begin{equation}
\begin{array}{l}
\displaystyle(Dc)_{n_{0},j_{0},k}(\beta )=\sum_{i=1}^{m-1}c(\beta
_{1},\ldots ,\beta _{i},(n_{0},kd_{\ast }+j_{0}),\beta _{i+1},\ldots ,\beta
_{m})1_{\{\beta _{i}^{\prime }<n_{0}<\beta _{i+1}^{\prime }\}}+\smallskip
\\
\quad \quad \displaystyle+c((n_{0},kd_{\ast }+j_{0}),\beta _{1},\ldots
,\beta _{m})1_{\{n_{0}<\beta _{1}^{\prime }\}}+c(\beta _{1},\ldots ,\beta
_{m},(n_{0},kd_{\ast }+j_{0}))1_{\{n_{0}>\beta _{m}^{\prime }\}}.%
\end{array}
\label{Sc}
\end{equation}

The aim of this section is to prove the non-degeneracy estimate (\ref{e2})
in next Lemma \ref{E}. But we first need to study the conditional
expectation of $\sigma_N$ given the randomness from $\chi_n$ and $U_n$.

\begin{lemma}
Assume ${\mathfrak{D}}(\varepsilon, r,R)$. We denote by ${\mathbb{E}}%
_{U,\chi }$ the conditional expectation with respect to $\sigma (U_{n},\chi
_{n},$ $n\in {\mathbb{N}}).$ Then%
\begin{equation}
{\mathbb{E}}_{U,\chi }(\sigma_N)\geq \lambda _{R}^{N}\sum_{\left\vert \alpha
\right\vert =N}c^{2}(\alpha )\chi^{\alpha^{\prime }},  \label{LV5}
\end{equation}
where $\lambda_R>0$ is given in Lemma \ref{lambdaR} and for $%
\alpha=((\alpha_1^{\prime },\alpha_1^{\prime \prime
}),\ldots,(\alpha_m^{\prime },\alpha_m^{\prime \prime }))$, we set $%
\alpha^{\prime }=(\alpha^{\prime }_1,\ldots,\alpha^{\prime }_m)$ and $%
\chi^{\alpha^{\prime }}=\prod_{i=1}^m\chi_{\alpha^{\prime }_i}$.
\end{lemma}

\textbf{Proof.} We set here $Z=Z(X)$. We recall that $X_{n,j}=\chi
_{n}V_{n,j}+(1-\chi _{n})U_{n,j}$ and we define (with $k(l)$ and $j(l)$
defined in (\ref{kj}))
\begin{equation*}
\widetilde{V}_{n,l}=V_{n,j(l)}^{k(l)}-{\mathbb{E}}(V_{n,j(l)}^{k(l)}),\quad
\overline{U}_{n,l}=(1-\chi _{n})U_{n,j(l)}^{k(l)}+\chi _{n}{\mathbb{E}}%
(V_{n,j(l)}^{k(l)})-{\mathbb{E}}(X_{n,j(l)}^{k(l)}).
\end{equation*}%
Then%
\begin{equation*}
Z_{n,l}=X_{n,j(l)}^{k(l)}-{\mathbb{E}}(X_{n,j(l)}^{k(l)})=\chi
_{n}V_{n,j(l)}^{k(l)}+(1-\chi _{n})U_{n,j(l)}^{k(l)}-{\mathbb{E}}%
(X_{n,j(l)}^{k(l)})=\chi _{n}\widetilde{V}_{n,l}+\overline{U}_{n,l}.
\end{equation*}%
So, we have
\begin{equation*}
Z^{\alpha }=\overline{Z}^{\alpha }+\chi ^{\alpha ^{\prime }}\widetilde{V}%
^{\alpha },
\end{equation*}%
where
\begin{equation*}
\overline{Z}^{\alpha }=\sum_{\substack{ (\beta ,\gamma )=\alpha , \\ \gamma
\neq \emptyset }}\chi ^{\beta ^{\prime }}\widetilde{V}^{\beta }\times
\overline{U}^{\gamma }.
\end{equation*}%
One has
\begin{equation}
{\mathbb{E}}_{U,\chi }(\overline{Z}^{\alpha }\widetilde{V}^{\theta })=\sum
_{\substack{ (\beta ,\gamma )=\alpha , \\ \gamma \neq \emptyset }}\chi
^{\beta }{\mathbb{E}}_{U,\chi }(\widetilde{V}^{\beta }\widetilde{V}^{\theta
})\times \overline{U}^{\gamma }=0,\quad
\mbox{for every $\alpha ,\theta $
s.t. $\left\vert \alpha\right\vert \leq \left\vert \theta \right\vert $}.
\label{C1}
\end{equation}%
This is because $\left\vert \beta \right\vert <\left\vert \alpha \right\vert
\leq \left\vert \theta \right\vert $, so there is at least one $\theta
_{i}\notin \beta $ and ${\mathbb{E}}_{U,\chi }(\widetilde{V}^{\theta
_{i}})=0.$ For the same reason, one has
\begin{equation}
{\mathbb{E}}_{U,\chi }(\widetilde{V}^{\alpha }\widetilde{V}^{\theta
})=0\quad
\mbox{for every  $\alpha ,\theta $ s.t. $\left\vert \alpha \right\vert <\left\vert
\theta \right\vert $}.  \label{C2}
\end{equation}%
We recall that $V_{n_{0},j_{0}}^{k}=\widetilde{V}_{n_{0},kd_{\ast
}+j_{0}}+E(V_{n_{0},j_{0}}^{k})$ and we use (\ref{Dnj}) in order to we
write
\begin{align*}
D_{n_{0},j_{0}}S_{N}(c,Z)&
=\sum_{m=0}^{N-1}(A_{m,1}^{n_{0},j_{0}}+A_{m,2}^{n_{0},j_{0}}+A_{m,3}^{n_{0},j_{0}}),\quad %
\mbox{where} \\
A_{m,1}^{n_{0},j_{0}}& =\sum_{k=0}^{k_{\ast }-1}\sum_{\beta \in \Lambda
_{n_{0},j_{0}}(m,k)}(k+1)(Dc)_{n_{0},j_{0},k}(\beta )\chi _{n_{0}}\widetilde{%
V}_{n_{0},kd_{\ast }+j_{0}}\chi ^{\beta ^{\prime }}\widetilde{V}^{\beta }, \\
A_{m,2}^{n_{0},j_{0}}& =\sum_{k=0}^{k_{\ast }-1}\sum_{\beta \in \Lambda
_{n_{0},j_{0}}(m,k)}(k+1)(Dc)_{n_{0},j_{0},k}(\beta )\chi _{n_{0}}\widetilde{%
V}_{n_{0},kd_{\ast }+j_{0}}\overline{Z}^{\beta }, \\
A_{m,3}^{n_{0},j_{0}}& =\sum_{k=0}^{k_{\ast }-1}\sum_{\beta \in \Lambda
_{n_{0},j_{0}}(m,k)}(k+1)(Dc)_{n_{0},j_{0},k}(\beta )\chi _{n_{0}}{\mathbb{E}%
}(V_{n_{0},j_{0}}^{k})Z^{\beta },
\end{align*}%
$\Lambda _{n_{0},j_{0}}(m,k)$ denoting the multi-indexes of length $m$ which
do not contain the pair $(n_{0},kd_{\ast }+j_{0})$. By (\ref{C1}) and (\ref%
{C2}), one has ${\mathbb{E}}_{U,\chi
}(A_{N-1,1}^{n_{0},j_{0}}A_{m,i}^{n_{0},j_{0}})=0$ for every $m\leq N-1$ and
$i=2,3$ and ${\mathbb{E}}_{U,\chi
}(A_{N-1,1}^{n_{0},j_{0}}A_{m,1}^{n_{0},j_{0}})=0$ for every $m<N-1$. Thus, $%
A_{N-1,1}^{n_{0},j_{0}}$ is orthogonal (in $L^{2}({\mathbb{P}}_{U,\chi })$)
to $D_{n_{0},j_{0}}S_{N}(c,Z)-A_{N-1,1}^{n_{0},j_{0}}$, so that
\begin{equation*}
{\mathbb{E}}_{U,\chi }(|D_{n_{0},j_{0}}S_{N}(c,Z)|^{2})\geq {\mathbb{E}}%
_{U,\chi }(|A_{N-1,1}^{n_{0},j_{0}}|^{2}).
\end{equation*}%
Therefore,
\begin{equation*}
{\mathbb{E}}_{U,\chi }(\sigma _{N})=\sum_{n_{0}=1}^{\infty
}\sum_{j_{0}=1}^{d_{\ast }}{\mathbb{E}}(|D_{n_{0},j_{0}}S_{N}(c,Z)|^{2})\geq
\sum_{n_{0}=1}^{\infty }\sum_{j_{0}=1}^{d_{\ast }}{\mathbb{E}}_{U,\chi
}(|A_{N-1,1}^{n_{0},j_{0}}|^{2}).
\end{equation*}%
Now, we write
\begin{equation*}
\begin{array}{l}
\displaystyle A_{N-1,1}^{n_{0},j_{0}}=\prod_{i=1}^{N}\chi
_{n_{i}}\sum_{n_{1}<\cdots <n_{N}}\sum_{l_{1},\ldots ,l_{N}\in \lbrack
m_{\ast }]}d_{n_{0},j_{0}}((n_{1},l_{1}),\ldots
,(n_{N},l_{N}))\prod_{i=1}^{N}\widetilde{V}_{n_{i},l_{i}}\quad \mbox{with}%
\smallskip  \\
\displaystyle d_{n_{0},j_{0}}(\alpha )=\sum_{i=1}^{N}\sum_{k=0}^{k_{\ast
}-1}(k+1)c(\alpha )1_{\alpha _{i}=(n_{0},kd_{\ast }+j_{0})},\quad |\alpha
|=N.%
\end{array}%
\end{equation*}%
For every $\alpha $ there exists at most one $(k,i)$ such that $\alpha
_{i}=(n_{0},kd_{\ast }+j_{0})$ so that
\begin{equation*}
\displaystyle d_{n_{0},j_{0}}^{2}(\alpha )=\sum_{i=1}^{N}\sum_{k=0}^{k_{\ast
}-1}(k+1)c^{2}(\alpha )1_{\alpha _{i}=(n_{0},kd_{\ast }+j_{0})}.
\end{equation*}%
By using (\ref{LV4}),
\begin{align*}
\sum_{n_{0}=1}^{\infty }\sum_{j_{0}=1}^{d_{\ast }}{\mathbb{E}}_{U,\chi
}(|A_{N-1,1}^{n_{0},j_{0}}|^{2})& \geq \sum_{n_{0}=1}^{\infty
}\sum_{j_{0}=1}^{d_{\ast }}\lambda _{R}^{N}\sum_{n_{1}<\cdots
<n_{N}}\sum_{l_{1},\ldots ,l_{N}}d_{n_{0},j_{0}}^{2}((n_{1},l_{1}),\ldots
,(n_{N},l_{N}))\prod_{i=1}^{N}\chi _{n_{i}} \\
& \geq \lambda _{R}^{N}\sum_{|\alpha |=N}c^{2}(\alpha )\chi ^{\alpha
^{\prime }}
\end{align*}%
and the statement holds. $\square $

\medskip

We can now prove the main result of this section.

\begin{lemma}
\label{E} Assume ${\mathfrak{D}}(\varepsilon, r,R)$. Let $c\in \mathcal{C}({%
\mathbb{R}})$ with $|c|_{N}>0$. For every $\eta >0$,%
\begin{equation}  \label{e2}
{\mathbb{P}}(\sigma _{N}\leq \eta ) \leq \frac{2e^{3}}{9}N\exp \Big(-\Big(%
\frac{\varepsilon\mathfrak{m}_r}{2}\Big)^{2N}\frac{|c|_{N}^{2}}{\delta
_{\ast}^{2}(c)}\Big)+\frac{2Kk_{\ast }N}{\lambda_R\varepsilon \mathfrak{m}%
_r|c|_{N}^{2/(k_{\ast }N)}}\eta ^{1/(k_{\ast }N)}.
\end{equation}%
where $K$ a universal constant (the one in the Carbery Wright inequality)
and $\lambda_R$ is given in Lemma \ref{lambdaR}.
\end{lemma}

\begin{remark}
Sometimes $|c|_{N}$ is small and we would like to use $|c|_{m}$ instead,
with $m<N.$ We denote $\left\vert c\right\vert
_{m+1,N}^{2}=\sum_{k=m+1}^{N}c^{2}(\alpha )$. Then for every $h\geq 1$ there
exists $C>0$ such that
\begin{equation}
{\mathbb{P}}(\sigma _{N}\leq \eta )\leq C\frac{\left\vert c\right\vert
_{m+1,N}^{2h}}{\eta ^{h}}+\frac{2e^{3}}{9}m\exp \Big(-\Big(\frac{\varepsilon
\mathfrak{m}_{r}}{2}\Big)^{2m}\frac{|c|_{m}^{2}}{\delta _{\ast }^{2}(c)}\Big)%
+\frac{2Kk_{\ast }m}{\lambda _{R}\varepsilon \mathfrak{m}_{r}|c|_{m}^{2/(k_{%
\ast }m)}}(4\eta )^{1/(k_{\ast }m)}.  \label{e2'}
\end{equation}%
Indeed: we denote $Q_{m+1,N,k_{\ast }}(c,X)=Q_{N,k_{\ast
}}(c,X)-Q_{m,k_{\ast }}(c,X)$ and we use the inequality
\begin{equation*}
\sigma _{N}=\left\vert DQ_{N,k_{\ast }}(c,X)\right\vert _{\mathcal{H}%
}^{2}\geq \frac{1}{2}\left\vert DQ_{m,k_{\ast }}(c,X)\right\vert _{\mathcal{H%
}}^{2}-\left\vert DQ_{m+1,N,k_{\ast }}(c,X)\right\vert _{\mathcal{H}}^{2}
\end{equation*}%
in order to obtain%
\begin{align*}
& {\mathbb{P}}(\sigma _{N}\leq \eta )\leq {\mathbb{P}}(\left\vert
DQ_{m+1,N,k_{\ast }}(c,X)\right\vert _{\mathcal{H}}^{2}\geq \eta )+{\mathbb{P%
}}(\left\vert DQ_{m,k_{\ast }}(c,X)\right\vert _{\mathcal{H}}^{2}\leq 4\eta )
\\
& \quad \leq {\mathbb{P}}(\left\vert DQ_{m+1,N,k_{\ast }}(c,X)\right\vert _{%
\mathcal{H}}^{2}\geq \eta )+\frac{2e^{3}}{9}m\exp \Big(-\Big(\frac{%
\varepsilon \mathfrak{m}_{r}}{2}\Big)^{2m}\frac{|c|_{m}^{2}}{\delta _{\ast
}^{2}(c)}\Big)+\frac{2Kk_{\ast }m}{\lambda _{R}\varepsilon \mathfrak{m}%
_{r}|c|_{m}^{2/(k_{\ast }m)}}(4\eta )^{1/(k_{\ast }m)}.
\end{align*}%
Using Chebyshev's inequality and Lemma \ref{M0}, for every $h$,
\begin{equation*}
{\mathbb{P}}(\left\vert DQ_{m+1,N,k_{\ast }}(c,X)\right\vert _{\mathcal{H}%
}^{2}\geq \eta )\leq \eta ^{-h}\left\Vert DQ_{m+1,N,k_{\ast
}}(c,X)\right\Vert _{\mathcal{H},2h}^{2h}\leq C\eta ^{-h}\left\vert
c\right\vert _{m+1,N}^{2h},
\end{equation*}%
so the proof of (\ref{e2'}) is completed.
\end{remark}

\textbf{Proof of Lemma} \ref{E}. We will use the Carbery--Wright inequality
that we recall here (see Theorem 8 in \cite{[CW]}). Let $\mu $ be a
probability law on ${\mathbb{R}}^{J}$ which is absolutely continuous with
respect to the Lebesgue measure and has a log-concave density. There exists
a universal constant $K$ such that for every polynomial $Q(x)$ of order $%
k_{\ast }N$ and for every $\eta >0$ one has%
\begin{equation}
\mu (x:\left\vert Q(x)\right\vert \leq \eta )\leq Kk_{\ast }N\Big(\frac{\eta%
}{V_{\mu }(Q)}\Big)^{1/(k_{\ast }N)},\quad \mbox{where } V_{\mu }(Q)=\Big(%
\int Q^{2}(x)d\mu (x)\Big)^{1/2}.  \label{NCW1}
\end{equation}

We will use this result in the following framework. We recall that the
coefficients $c(\alpha )$ are null except a finite number of them. So we may
find $M$ such that, if $\left\vert \alpha \right\vert =m$ and $\alpha
_{m}^{\prime }>M$ then $c(\alpha )=0.$ It follows that we may write (see \ref%
{Dnj}))
\begin{equation*}
\sigma _{N}=q_{\chi ,\overline{U}}(V)
\end{equation*}%
where $q_{q,\overline{U}}(V)$ is a polynomial of order $k_{\ast }N$ with
unknowns $V_{n,j},n\leq M,j\leq d_{\ast }$ and coefficients depending on $%
\chi _{n}$ and $\overline{U}_{n,j,k}.$ Moreover we recall that ${\mathbb{P}}%
_{U,\chi }$ is the conditional probability with respect to $\sigma
(U_{i},\chi _{i},i\in {\mathbb{N}}).$ We denote by $\mu $ the law of $%
(V_{n,j},n\leq M,j\leq d_{\ast })$ under ${\mathbb{P}}_{U,\chi }$: this is a
product of laws of the form $c\psi _{r}(\left\vert x-\overline{x}\right\vert
^{2})dx$ so it is log-concave. So we are able to use (\ref{NCW1}). Using (%
\ref{LV5})%
\begin{equation*}
V_{\mu }(q_{\chi ,\overline{U}})\geq \int \vert q_{\chi ,\overline{U}%
}(x)\vert d\mu (x)={\mathbb{E}}_{U,\chi }(\sigma _{N})\geq
\lambda_R^N\sum_{|\beta|=N}c^{2}(\beta )\chi ^{\beta ^{\prime }}.
\end{equation*}%
We take now $\theta >0$ (to be chosen in a moment) and we use (\ref{NCW1})
in order to obtain%
\begin{eqnarray}
{\mathbb{P}}(\sigma _{N}\leq \eta ) &=&{\mathbb{P}}(q_{\chi ,\overline{U}%
}(V)\leq \eta )  \notag \\
&\leq &{\mathbb{P}}(V_{\mu }(q_{\chi ,\overline{U}})\leq \theta )+{\mathbb{E}%
}({\mathbb{P}}_{V,\chi }(q_{\chi ,\overline{U}}(V)\leq \eta )1_{\{V_{\mu
}(q_{\chi ,\overline{U}})\geq \theta \}})  \notag \\
&\leq &{\mathbb{P}}\Big(\sum_{\beta \in \Lambda _{N}}c^{2}(\beta )\chi
^{\beta ^{\prime }}\leq \frac{\theta }{\lambda_R^N}\Big)+Kk_{\ast }N(\eta
/\theta )^{1/(k_{\ast }N)}.  \label{landa}
\end{eqnarray}

The first term in the above inequality is estimated in Appendix \ref{hoef}.
In order to fit in the notation used there we denote $\Lambda _{N}(\beta
^{\prime })=\{\alpha : |\alpha|=N\mbox{ and }\alpha ^{\prime }=\beta
^{\prime }\}$ and $\overline{c}^{2}(\beta ^{\prime })=\sum_{\alpha \in
\Lambda _{N}(\beta ^{\prime })}c^{2}(\alpha ).$ Then
\begin{equation*}
\sum_{\beta \in \Lambda _{N}}c^{2}(\beta )\chi ^{\beta ^{\prime
}}=\sum_{\left\vert \beta ^{\prime }\right\vert =N}\overline{c}^{2}(\beta
^{\prime })\chi ^{\beta ^{\prime }}=\Psi _{N}(\overline{c}^2).
\end{equation*}%
Now we apply Lemma \ref{H} with $x=\theta /\lambda_R^N.$ Recall that $%
p=\varepsilon \mathfrak{m}_r$ and we have the restriction
\begin{equation}
\theta =\lambda_R^Nx< \lambda_R^N\Big(\frac{p}{2}\Big)^{N}\sum_{|\beta
^{\prime }| =N}\overline{c}^{2}(\beta ^{\prime })=\lambda_R^N\Big(\frac{%
\varepsilon \mathfrak{m}_r}{2}\Big)^{N}|\overline{c}|_{N}^{2}.  \label{theta}
\end{equation}%
We have $\left\vert \overline{c}\right\vert _{N}^{2}=|c|_{N}^{2}$ and
\begin{equation*}
\delta _{N}^{2}(\overline{c})=\max_{n}\sum_{n\in \beta ^{\prime },\left\vert
\beta ^{\prime }\right\vert =N}\overline{c}^{2}(\alpha ^{\prime
})=\max_{n}\sum_{n\in \alpha ^{\prime },|\alpha|=N}c^{2}(\alpha )=\delta
_{\ast}^{2}(c).
\end{equation*}%
Then (\ref{h3}) gives%
\begin{equation*}
{\mathbb{P}}\Big(\Psi _{N}(\overline{c}^2)\leq \frac \theta{\lambda_R^N}%
\Big) \leq \frac{2e^{3}}{9}N\exp \Big(-\frac{(\theta/\lambda_R^N)^{2}}{%
\delta _{\ast}^{2}(c)\left\vert c\right\vert _{N}^{2}}\Big).
\end{equation*}%
Inserting this in (\ref{landa}) we obtain
\begin{equation*}
{\mathbb{P}}(\sigma _{N}\leq \eta )\leq \frac{2e^{3}}{9}N\exp \Big(-\frac{%
(\theta/\lambda_R)^{2}}{\delta _{\ast}^{2}(c)\left\vert c\right\vert _{N}^{2}%
}\Big)+Kk_{\ast }N(\eta /\theta )^{1/(k_{\ast }N)}.
\end{equation*}
Now, $\theta$ is any constant satisfying the restriction (\ref{theta}). So,
by letting $\theta\uparrow \lambda_R^N((\varepsilon \mathfrak{m}_r)/2)^{N}|%
\overline{c}|_{N}^{2}=\lambda_R^N((\varepsilon \mathfrak{m}%
_r)/2)^{N}|c|_{N}^{2}$, we finally obtain (\ref{e2}). $\square $

\subsection{Proof of Theorem \protect\ref{D}}

\label{sect:proof} The goal of this section is to give the proof of Theorem %
\ref{D} so we use the notation from Section \ref{sect:3}.

We take $q\in{\mathbb{N}}$, $q\geq 1$, and we consider the sequence $%
\lambda_q=\frac{q}{q+k}$. Since $\lambda^2_q\uparrow 1$ as $q\to\infty$, we
can find $q$ such that such that $\lambda_q^2<\theta\leq\lambda^2_{q+1}$.
And since $\lambda^2_{q+1}\leq \lambda_q$, we get $\lambda_q^2<\theta\leq%
\lambda_q$. We work with this value of $q$ and we write simply $\lambda$ in
place of $\lambda_q$. Moreover, in the following, $C>0$ stands for a
constant which may vary from line to line and which depends on the
parameters in the statements but not on the coefficients $c,d\in\mathcal{C}({%
\mathbb{R}})$.

We define $a=\theta/\lambda$, so $\frac{1}{1+k}\leq \lambda<a\leq 1$. We
consider $\eta ,\delta \in (0,1)$, to be chosen in the sequel, and we use
the regularization Lemma \ref{R} (see (\ref{1b})) with the above choice of $%
q $ and $a$. This gives
\begin{eqnarray*}
&&\left\vert {\mathbb{E}}(f(Q_{N,k_{\ast }}(c,X)))-{\mathbb{E}}(f_{\delta
}(Q_{N,k_{\ast }}(c,X)))\right\vert \\
&\leq &C\left\Vert f\right\Vert _{\infty }\Big({\mathbb{P}}^{a}(\det \sigma
_{Q_{N,k_{\ast }}(c,X)}\leq \eta )+\frac{\delta ^{q}}{\eta ^{2q}}\left\Vert {%
\mathcal{K}}_{q,0}(Q_{N,k_{\ast }}(c,X))\right\Vert _{2}\Big) \\
&\leq &C\left\Vert f\right\Vert _{\infty }\Big({\mathbb{P}}^{a}(\det \sigma
_{Q_{N,k_{\ast }}(c,X)}\leq \eta )+\frac{\delta ^{q}}{\eta ^{2q}}%
\,|c|^q(1+|c|)^{4q}\Big),
\end{eqnarray*}%
the latter inequality following from (\ref{est-Cq0}). Moreover by (\ref{e2'}%
) (therein, $\sigma_N=\det\sigma _{Q_{N,k_{\ast }}}$), for every $h\geq 1$
(recall that $\overline{m}=m\vee m^{\prime })$%
\begin{equation*}
{\mathbb{P}}(\det\sigma _{Q_{N,k_{\ast }}(c,X)}\leq \eta )\leq C\Big(\frac{%
\left\vert c\right\vert _{m+1,N}^{2h}}{\eta ^{h}} +e_{m,N}(c) +\frac{1}{%
|c|_{m}^{2/(k_{\ast }m)}}\eta ^{1/(k_{\ast }\overline{m})}\Big).
\end{equation*}
So,
\begin{align*}
&\left\vert {\mathbb{E}}(f(Q_{N,k_{\ast }}(c,X)))-{\mathbb{E}}(f_{\delta
}(Q_{N,k_{\ast }}(c,X)))\right\vert \\
&\leq C\left\Vert f\right\Vert _{\infty }\Big(\frac{\left\vert c\right\vert
_{m+1,N}^{2ha}}{\eta ^{ha}}+e_{m,N}^{a}(c)+\frac{\eta ^{a/(k_{\ast }%
\overline{m})}}{|c|_{m}^{2a/(k_{\ast }m)}}+|c|^q(1+|c|)^{4q}\frac{\delta ^{q}%
}{\eta ^{2q}}\Big).
\end{align*}%
A similar estimate holds for $Q_{N,k_{\ast }}(d,Y).$ We use now $d_k$
defined in (\ref{D3''}). Since $\left\Vert f_{\delta }\right\Vert _{k,\infty
}\leq \delta ^{-k}\left\Vert f\right\Vert _{\infty }$ one has%
\begin{equation*}
\left\vert {\mathbb{E}}(f_{\delta }(Q_{N,k_{\ast }}(c,X)))-{\mathbb{E}}%
(f_{\delta }(Q_{N,k_{\ast }}(d,Y))\right\vert \leq k\delta
^{-k}d_{k}\left\Vert f\right\Vert _{\infty }.
\end{equation*}%
Putting this together,
we get%
\begin{align*}
&\left\vert {\mathbb{E}}(f(Q_{N,k_{\ast }}(c,X)))-{\mathbb{E}}%
(f(Q_{N,k_{\ast }}(d,Y))\right\vert \\
&\leq C\max \Big(1,\Big(|c|_m^{-\frac{2}{k_\ast m}}+|d|_{m^{\prime }}^{-%
\frac{2}{k_\ast m^{\prime }}}\Big)^a\Big)\left\Vert f\right\Vert _{\infty
}\times \\
&\quad\times \Big(\frac{\left\vert c\right\vert _{m+1,N}^{2ha}}{\eta ^{ha}}+%
\frac{\left\vert d\right\vert _{m^{\prime }+1,N}^{2ha}}{\eta ^{ha}}%
+e_{m,N}^{a}(c)+e_{m^{\prime },N}^{a}(d)+\eta ^{a/(k_{\ast }\overline{m}%
)}+(1+|c|+|d|)^{5q}\frac{\delta ^{q}}{\eta ^{2q}}+\delta ^{-k}d_{k}\Big).
\end{align*}%
We optimize first on $\delta :$ we take $\delta=d_k^{1/(q+k)}%
\eta^{2q/(q+k)}(1+|c|+|d|)^{-5q/(q+k)}$ and we obtain (recall that $%
\lambda=\frac q{q+k}\in(0,1)$),
\begin{align*}
(1+|c|+|d|)^{5q}\frac{\delta ^{q}}{\eta ^{2q}} &=\delta ^{-k}d_{k} = \eta
^{-2k\lambda}d_{k}^{\lambda }(1+|c|+|d|)^{5k\lambda} \leq \eta
^{-2k\lambda}d_{k}^{\lambda }(1+|c|+|d|)^{5k}.
\end{align*}%
It follows that
\begin{align*}
&\left\vert {\mathbb{E}}(f(Q_{N,k_{\ast }}(c,X)))-{\mathbb{E}}%
(f(Q_{N,k_{\ast }}(d,Y))\right\vert \\
&\leq C\max \Big(1,\Big(|c|_m^{-\frac{2}{k_\ast m}}+|d|_{m^{\prime }}^{-%
\frac{2}{k_\ast m^{\prime }}}\Big)^a\Big)\left\Vert f\right\Vert _{\infty
}(1+|c|+|d|)^{5k}\times \\
&\quad\times \Big(\frac{\left\vert c\right\vert _{m+1,N}^{2ha}}{\eta ^{ha}}+%
\frac{\left\vert d\right\vert _{m^{\prime }+1,N}^{2ha}}{\eta ^{ha}}%
+e_{m,N}^{a}(c)+e_{m^{\prime },N}^{a}(d)+\eta ^{a/(k_{\ast }\overline{m})}+
\eta ^{-2k\lambda }d_{k}^{\lambda } \Big).
\end{align*}%
We optimize now on $\eta:$ we take $\eta =d_k^{\lambda k_\ast\overline{m}%
/(a+2\lambda kk_\ast\overline{m})}$, so that
\begin{equation*}
\eta ^{-2k\lambda}d_{k}^{\lambda }= \eta ^{a/(k_{\ast }\overline{m})}=
d_k^{\lambda a/(a+2 \lambda kk_\ast\overline{m})} \leq d_k^{\lambda
a/(1+2kk_\ast\overline{m})},
\end{equation*}
the latter inequality follows from $d_k\leq 1$ and, since $a,\lambda\in(0,1)$%
, $a+2\lambda kk_\ast\overline{m}\leq 1+2kk_\ast\overline{m}$. By inserting,
\begin{align*}
&\left\vert {\mathbb{E}}(f(Q_{N,k_{\ast }}(c,X)))-{\mathbb{E}}%
(f(Q_{N,k_{\ast }}(d,Y))\right\vert \\
&\leq C\max \Big(1,\Big(|c|_m^{-\frac{2}{k_\ast m}}+|d|_{m^{\prime }}^{-%
\frac{2}{k_\ast m^{\prime }}}\Big)^a\Big)\left\Vert f\right\Vert _{\infty
}(1+|c|+|d|)^{5k }\times \\
&\quad\times \Big(\frac{\left\vert c\right\vert _{m+1,N}^{2ha}}{\eta ^{ha}}+%
\frac{\left\vert d\right\vert _{m^{\prime }+1,N}^{2ha}}{\eta ^{ha}}%
+e_{m,N}^{a}(c)+e_{m^{\prime },N}^{a}(d)+d_k^{\lambda a/(1+2kk_\ast\overline{%
m})}\Big).
\end{align*}%
Since $\left\vert c\right\vert^2 _{m+1,N}$ $\leq d_{k}^{\frac{k_\ast%
\overline{m}}{2kk_{\ast }\overline{m}+1}},$
\begin{equation*}
\frac{|c|_{m+1,N}^{2ha}}{\eta^{ha}} \leq d_k^{ah(\frac{k_\ast \overline{m}}{%
1+2kk_\ast \overline{m}}- \frac{\lambda k_\ast \overline{m}}{a+2\lambda
kk_\ast \overline{m}})}.
\end{equation*}
We note that the above exponent is positive because $a>\lambda$. 
So, we choose $h\geq 1$ and such that
\begin{equation*}
ah\Big(\frac{k_\ast \overline{m}}{1+2kk_\ast \overline{m}}- \frac{\lambda
k_\ast \overline{m}}{a+2\lambda kk_\ast \overline{m}}\Big)\geq \frac{\lambda
a}{(1+2kk_\ast\overline{m})},
\end{equation*}
so that
\begin{equation*}
\frac{|c|_{m+1,N}^{2ha}}{\eta^{ha}} \leq d_k^{\frac{\lambda a}{(1+2kk_\ast%
\overline{m})}}.
\end{equation*}
A similar estimate holds with $|c|_{m+1,N}^{2ha}$ replaced by $%
|d|_{m^{\prime }+1,N}^{2ha}$. We then obtain
\begin{align*}
&\left\vert {\mathbb{E}}(f(Q_{N,k_{\ast }}(c,X)))-{\mathbb{E}}%
(f(Q_{N,k_{\ast }}(d,Y))\right\vert \\
&\leq C\max \Big(1,\Big(|c|_m^{-\frac{2}{k_\ast m}}+|d|_{m^{\prime }}^{-%
\frac{2}{k_\ast m^{\prime }}}\Big)^a\Big)\left\Vert f\right\Vert _{\infty
}(1+|c|+|d|)^{5k} \Big(e_{m,N}^{a}(c)+e_{m^{\prime },N}^{a}(d)+d_k^{\lambda
a/(1+2kk_\ast\overline{m})}\Big).
\end{align*}%
The statement now follows by recalling that $\lambda a=\theta$ and, from (%
\ref{D3''}), $d_{k}\leq C(d_{k}(Q_{N,k_{\ast }}(c,X),Q_{N,k_{\ast }}(d,Y))$ $ +
\vert c\vert^{\frac{2(2kk_{\ast }\overline{m}+1)}{k_\ast\overline{m}}}
_{m+1,N}+|d|^{\frac{2(2 kk_{\ast }\overline{m}+1)}{k_\ast\overline{m}}}
_{m^{\prime }+1,N}).$ $\square$

\appendix

\section{An iterated Hoeffding's inequality}

\label{hoef}

In this section we work with multi-indexes $\alpha =(\alpha _{1},\ldots
,\alpha _{m})\in {\mathbb{N}}^{m}$ with $1\leq\alpha _{1}<\ldots <\alpha
_{m} $ and we look to
\begin{equation*}
\Psi _{m}(c^2)=\sum_{\left\vert \alpha \right\vert =m}c^{2}(\alpha )\chi
^{\alpha },
\end{equation*}%
where $\chi_n$, $n\in{\mathbb{N}}$, denote independent Bernoulli random
variables and $\chi^\alpha=\prod_{i=1}^m\chi_{\alpha_i}$. We denote%
\begin{equation*}
\left\vert c\right\vert _{m}^{2}=\sum_{\left\vert \alpha \right\vert
=m}c^{2}(\alpha ),\quad and\quad \delta _{m}^{2}(c)=\max_{n}\sum_{\left\vert
\alpha \right\vert =m,n\in \alpha }c^{2}(\alpha ).
\end{equation*}

\begin{lemma}
\label{H}Let $p={\mathbb{P}}(\chi _{j}=1)\in(0,1)$. If
\begin{equation}
x< \Big(\frac{p}{2}\Big)^{N}\left\vert c\right\vert _{N}^{2}  \label{hh2}
\end{equation}%
then%
\begin{equation}
{\mathbb{P}}(\Psi _{N}(c^{2})\leq x)\leq \frac{2e^{3}}{9}N\exp \Big(-\frac{%
x^{2}}{\delta _{N}^{2}(c)\left\vert c\right\vert _{N}^{2}}\Big).  \label{h3}
\end{equation}
\end{lemma}

\textbf{Proof.} We proceed by recurrence on $N.$ If $N=1$ we have
\begin{align*}
&{\mathbb{P}}(\Psi _{N}(c^{2}) \leq x) ={\mathbb{P}}\Big(\sum_{n}c^{2}(n)%
\chi _{n}\leq x\Big) \\
&\leq {\mathbb{P}}\Big(p\sum_{n}c^{2}(n)\leq 2x\Big)+{\mathbb{P}}\Big(%
\sum_{n}c^{2}(n)(p-\chi _{j})\geq x\Big) ={\mathbb{P}}\Big(%
\sum_{n}c^{2}(n)(p-\chi _{j})\geq x\Big),
\end{align*}%
the latter inequality following from (\ref{hh2}). And by Hoeffding's
inequality%
\begin{equation*}
{\mathbb{P}}\Big(\sum_{j}c^{2}(j)(p-\chi _{j})\geq x\Big)\leq \exp \Big(-%
\frac{2x^{2}}{\sum_{j}c^{4}(j)}\Big).
\end{equation*}%
Since
\begin{equation*}
\sum_{j}c^{4}(j)\leq \max_{j}c^{2}(j)\times \sum_{j}c^{2}(j)=\delta
_{1}^{2}(c)\left\vert c\right\vert _{1}^{2}
\end{equation*}%
(\ref{h3}) follows for $N=1$. We suppose now that (\ref{h3}) holds for $N-1$
and we prove it for $N.$ For $\beta $ with $\left\vert \beta \right\vert
=N-1 $ we define $c_{n}(\beta )=c(\beta ,n)1_{\{\beta _{N-1}<n\}}$ and we
write%
\begin{equation*}
\Psi _{N}(c^2)=\sum_{\left\vert \alpha \right\vert =N}c^{2}(\alpha )\chi
^{\alpha }=\sum_{n=N}^{\infty }\chi _{n}\sum_{\left\vert \beta \right\vert
=N-1,\beta _{N-1}<n}c^{2}(\beta ,n)\chi ^{\beta }=\sum_{n=N}^{\infty }\chi
_{n}\Psi _{N-1}(c^2_{n}).
\end{equation*}%
Then
\begin{equation*}
{\mathbb{P}}(\Psi _{N}(c^2)\leq x)\leq {\mathbb{P}}\Big(\sum_{n=N}^{\infty
}\Psi _{N-1}(c^2_{n})\leq \frac{2x}{p}\Big)+{\mathbb{P}}\Big(%
\sum_{n=N}^{\infty }(p-\chi _{n})\Psi _{N-1}(c^2_{n})\geq x\Big)=:a+b.
\end{equation*}%
We estimate first $b.$ We write%
\begin{equation*}
\sum_{n=N}^{\infty }\Psi _{N-1}(c^2_{n})=\sum_{\left\vert \beta \right\vert
=N-1}d_{n}^{2}(\beta )\chi ^{\beta }\quad \mbox{with}\quad d^{2}(\beta
)=\sum_{n>\beta _{N-1}}^{\infty }c^{2}(\beta ,n).
\end{equation*}%
Notice that
\begin{equation*}
\left\vert d\right\vert _{N-1}^{2}=\sum_{\left\vert \beta \right\vert
=N-1}\sum_{n>\beta _{N-1}}^{\infty }c^{2}(\beta ,n)=\sum_{\left\vert \alpha
\right\vert =N}c^{2}(\alpha )=\left\vert c\right\vert _{N}^{2}
\end{equation*}%
and%
\begin{equation*}
\delta _{N-1}^{2}(d)=\max_{k}\sum_{\left\vert \alpha \right\vert =N-1,k\in
\alpha }d^{2}(\alpha )=\max_{k}\sum_{\left\vert \alpha \right\vert =N-1,k\in
\alpha }\sum_{n>\alpha _{N-1}}^{\infty }c^{2}(\alpha ,n)\leq
\max_{k}\sum_{\left\vert \beta \right\vert =N,k\in \beta }c^{2}(\beta
)=\delta _{N}^{2}(c).
\end{equation*}%
We also have%
\begin{equation*}
\frac{2x}{p}< \frac{2}{p}\Big(\frac{p}{2}\Big)^{N}\left\vert c\right\vert
_{N}^{2}=\Big(\frac{p}{2}\Big)^{N-1}\left\vert d\right\vert _{N}^{2}
\end{equation*}%
so we can use the recurrence hypothesis and we get%
\begin{equation}
b={\mathbb{P}}\Big(\Psi _{N-1}(d^{2})\leq \frac{2x}{p}\Big)\leq \frac{2e^{3}%
}{9}(N-1)\exp \Big(-\frac{(2x/p)^{2}}{\delta _{N-1}^{2}(d)\left\vert
d\right\vert _{N-1}^{2}}\Big)\leq \frac{2e^{3}}{9}(N-1)\exp \Big(-\frac{x^{2}%
}{\delta _{N}^{2}(c)\left\vert c\right\vert _{N}^{2}}\Big).  \label{hh8}
\end{equation}%
We estimate now $a.$ We use Corollary 1.4 pg 1654 in Bentkus \cite{[Be]} which
asserts the following: if $M_{k},k\in {\mathbb{N}}$ is a martingale such
that $\left\vert M_{k}-M_{k-1}\right\vert \leq h_{k}$ almost surely, then,
for every $n\in {\mathbb{N}},$
\begin{equation}
{\mathbb{P}}(M_{n}\geq x)\leq \frac{2e^{3}}{9}\exp (-\frac{x^{2}}{%
\sum_{j=1}^{n}h_{j}^{2}}).  \label{hh1}
\end{equation}%
Since $0\leq \chi _{n}\leq 1$ we have
\begin{equation*}
\Psi _{N-1}(c^2_{n})\leq \sum_{\left\vert \beta \right\vert =n,\beta
_{N-1}<n}c^{2}(\beta ,n)=:h_{n}.
\end{equation*}%
Notice that $h_{n}\leq \delta _{N}^{2}(c)$ so that
\begin{equation*}
\sum_{j=1}^{n}h_{j}^{2}\leq \delta _{N}^{2}(c)\sum_{j=1}^{n}h_{j}=\delta
_{N}^{2}(c)\left\vert c\right\vert _{N}^{2}.
\end{equation*}%
So, using (\ref{hh1})%
\begin{equation*}
a={\mathbb{P}}\Big(\sum_{j=1}^{\infty }(p-\chi _{j})\Psi _{N-1}(c_{n})\geq x%
\Big)\leq \frac{2e^{3}}{9}\exp \Big(-\frac{x^{2}}{\delta
_{N}^{2}(c)\left\vert c\right\vert _{N}^{2}}\Big).
\end{equation*}%
This, together with (\ref{hh8}), gives (\ref{h3}). $\square $

\section{Norms}

\label{app:norms}

The aim of this section is to prove Lemma \ref{L1}. For $F=(F_1,\ldots,F_d)$
We work with the norms
\begin{eqnarray*}
\left\vert F\right\vert _{1,k} &=&\sum_{j=1}^d\sum_{1=1}^{k}\left\vert
D^{i}F_j\right\vert _{\mathcal{H}^{\otimes i}},\quad \quad \left\vert
F\right\vert _{k}=\left\vert F\right\vert +\left\vert F\right\vert _{1,k} \\
\left\Vert F\right\Vert _{1,k,p} &=&\Vert \left\vert F\right\vert
_{1,k}\Vert _{p},\quad \quad \left\Vert F\right\Vert _{k,p}=\left\Vert
F\right\Vert _{p}+\left\Vert F\right\Vert _{1,k,p}.
\end{eqnarray*}

To begin we give several easy computational rules:%
\begin{eqnarray}
\left\vert FG\right\vert _{k} &\leq &C\sum_{k_{1}+k_{2}=k}\left\vert
F\right\vert _{k_{1}}\left\vert G\right\vert _{k_{2}}  \label{Nn1} \\
\left\vert \left\langle DF,DG\right\rangle \right\vert _{k} &\leq
&C\sum_{k_{1}+k_{2}=k}\left\vert F\right\vert _{1,k_{1}+1}\left\vert
G\right\vert _{1,k_{2}+1}\quad and\quad  \label{Nn2} \\
\left\vert \frac{1}{G}\right\vert _{k} &\leq &\frac{C}{\left\vert
G\right\vert }\sum_{l=0}^{k}\frac{\left\vert G\right\vert _{k}^{l}}{%
\left\vert G\right\vert ^{l}}.  \label{Nn3}
\end{eqnarray}%
Now, for $F=(F_{1},\ldots ,F_{d})$ we consider the Malliavin covariance
matrix $\sigma _{F}^{i,j}=\left\langle DF^{i},DF^{j}\right\rangle $ and, if
$\det \sigma _{F}\neq 0,$ we denote $\gamma _{F}=\sigma _{F}^{-1}.$ We write%
\begin{equation*}
\gamma _{F}^{i,j}=\frac{\widehat{\sigma }_{F}^{i,j}}{\det \sigma _{F}}
\end{equation*}%
where $\widehat{\sigma }_{F}^{i,j}$ is the algebraic complement . Then,
using (\ref{Nn1})
\begin{equation*}
\big\vert \gamma _{F}^{i,j}\big\vert _{k}\leq C\sum_{k_{1}+k_{2}=k}\big\vert
\widehat{\sigma }_{F}^{i,j}\big\vert _{k_{1}}\left\vert \frac{1}{\det \sigma
_{F}}\right\vert _{k_{2}}.
\end{equation*}%
By (\ref{Nn1}) and (\ref{Nn2}), $\left\vert \widehat{\sigma }%
_{F}^{i,j}\right\vert _{k_{1}}\leq C\left\vert F\right\vert
_{1,k_{1}+1}^{2(d-1)}$ and $\left\vert \det \sigma _{F}\right\vert
_{k_{2}}\leq C\left\vert F\right\vert _{1,k_{2}+1}^{2d}.$ Then, using (\ref%
{Nn3})
\begin{equation*}
\left\vert \frac{1}{\det \sigma _{F}}\right\vert _{k_{2}}\leq \frac{C}{%
\left\vert \det \sigma _{F}\right\vert }\sum_{l=0}^{k_{2}}\frac{\left\vert
\det \sigma _{F}\right\vert _{k_{2}}^{l}}{\left\vert \det \sigma
_{F}\right\vert ^{l}}\leq \frac{C}{\left\vert \det \sigma _{F}\right\vert }%
\sum_{l=0}^{k_{2}}\frac{\left\vert F\right\vert _{1,k_{2}+1}^{2ld}}{%
\left\vert \det \sigma _{F}\right\vert ^{l}}
\end{equation*}%
so that
\begin{equation}
\big\vert \gamma _{F}^{i,j}\big\vert _{k}\leq C\frac{\left\vert F\right\vert
_{1,k+1}^{2(d-1)}}{\left\vert \det \sigma _{F}\right\vert }\sum_{l=0}^{k}%
\Big(\frac{\left\vert F\right\vert _{1,k+1}^{2d}}{\left\vert \det \sigma
_{F}\right\vert }\Big)^{l}\leq C\frac{\left\vert F\right\vert
_{1,k+1}^{2(d-1)}}{\left\vert \det \sigma _{F}\right\vert }\Big(1+\frac{%
\left\vert F\right\vert _{1,k+1}^{2d}}{\left\vert \det \sigma
_{F}\right\vert }\Big)^{k}  \label{Nn4}
\end{equation}

We denote%
\begin{equation}
\alpha _{k}=\frac{\left\vert F\right\vert _{1,k+1}^{2(d-1)}(\left\vert
F\right\vert _{1,k+1}+\left\vert LF\right\vert _{k})}{\left\vert \det \sigma
_{F}\right\vert },\quad \beta _{k}=\frac{\left\vert F\right\vert
_{1,k+1}^{2d}}{\left\vert \det \sigma _{F}\right\vert }  \label{Nn4'}
\end{equation}%
and%
\begin{equation}
{\mathcal{K}}_{n,k}(F)=(\left\vert F\right\vert _{1,k+n+1}+\left\vert
LF\right\vert _{k+n})^{n}(1+\left\vert F\right\vert _{1,k+n+1})^{2d(2n+k)}.
\label{Nn4''}
\end{equation}%
We also recall that for $\eta >0,$ we consider a function $\Psi _{\eta }\in
C^{\infty }({\mathbb{R}})$ such that $1_{(0,\eta )}\leq \Psi _{\eta }\leq
1_{(0,2\eta )} $ and $\Vert \Psi _{\eta }^{(k)}\Vert _{\infty }\leq
C_{k}\eta ^{-k},\forall k\in {\mathbb{N}}.$ Then we take $\Phi _{\eta
}=1-\Psi _{\eta }.$

\begin{lemma}
\textbf{A}. For every $k,n\in {\mathbb{N}}$ there exists a universal
constant $C$ (depending on $k$ and $n)$ such that, for $\omega $ such that $%
\det \sigma _{F}(\omega )>0,$
\begin{equation}
\left\vert H_{\rho }^{(n)}(F,G)\right\vert _{k}\leq C\alpha
_{k+n}^{n}\sum_{p_{1}+p_{2}=k+n}\left\vert G\right\vert _{p_{2}}(1+\beta
_{k+n})^{p_{1}}.  \label{Nn5}
\end{equation}%
\textbf{B}. For every $\eta >0$
\begin{equation}
\left\vert H_{\rho }^{(n)}(F,\Phi _{\eta }(\det \sigma _{F})G)\right\vert
_{k}\leq \frac{C}{\eta ^{2n+k}}\times {\mathcal{K}}_{n,k}(F)\times
\left\vert G\right\vert _{k+n}.  \label{Nn6}
\end{equation}
\end{lemma}

\bigskip

\textbf{Proof A}. We first prove (\ref{Nn5}) for $n=1$. We have
\begin{equation*}
H_{i}^{(1)}(F,G)=-\sum_{j=1}^{m}G\gamma _{F}^{i,j}LF^{j}+G\langle D\gamma
_{F}^{i,j},DF^{j}\rangle +\gamma _{F}^{i,j}\langle DG,DF^{j}\rangle.
\end{equation*}%
Using (\ref{Nn1})
\begin{eqnarray*}
\big\vert H_{i}^{(1)}(F,G)\big\vert _{k} &\leq
&C\sum_{k_{1}+k_{2}+k_{3}=k}\left( \left\vert \gamma _{F}\right\vert
_{k_{1}}\left\vert LF\right\vert _{k_{2}}\left\vert G\right\vert
_{k_{3}}+\left\vert \gamma _{F}\right\vert _{k_{1}+1}\left\vert F\right\vert
_{1,k_{2}+1}\left\vert G\right\vert _{k_{3}}+\left\vert \gamma
_{F}\right\vert _{k_{1}}\left\vert F\right\vert _{1,k_{2}+1}\left\vert
G\right\vert _{k_{3}+1}\right) \\
&\leq &C(\left\vert F\right\vert _{k+1}+\left\vert LF\right\vert
_{k})\sum_{p_{1}+p_{2}\leq k}\left( \left\vert \gamma _{F}\right\vert
_{p_{1}+1}\left\vert G\right\vert _{p_{2}}+\left\vert \gamma _{F}\right\vert
_{p_{1}}\left\vert G\right\vert _{p_{2}+1}\right) .
\end{eqnarray*}%
For $n>1$, we use recurrence and we obtain
\begin{equation*}
\left\vert H_{\gamma }^{(n)}(F,G)\right\vert _{k}\leq C(\left\vert
F\right\vert _{k+n+1}+\left\vert LF\right\vert _{k+n})^{n}\sum_{p_{1}+\ldots
+p_{n+1}\leq k+n-1}\prod_{i=1}^{n}\left\vert \gamma _{F}\right\vert
_{p_{i}}\times \left\vert G\right\vert _{p_{n+1}}.
\end{equation*}%
Then, using (\ref{Nn1}) first and (\ref{Nn4}) secondly, (\ref{Nn5}) follows.

\medskip

\textbf{B}. Let $G_{\eta }=\Phi _{\eta }(\det \sigma _{F})G).$ For every $%
p\in {\mathbb{N}}$ one has $\left\vert G_{\eta }\right\vert _{p}\leq C\eta
^{-p}\left\vert G\right\vert _{p}\left\vert F\right\vert _{1,p+1}^{d}.$
Moreover one has $H_{\rho }^{(n)}(F,G_{\eta })=1_{\{\det \sigma _{\Phi
}>\eta /2\}}H_{\rho }^{(n)}(F,G_{\eta }).$ So (\ref{Nn5}) implies (\ref{Nn6}%
). $\square $


\addcontentsline{toc}{section}{References}

\begin{thebibliography}{99}
\bibitem{[BC-CLT]} {\small Bally V., Caramellino L.: Asymptotic
development for the CLT in total variation distance. \emph{Bernoulli}, 22,
2442-2485 22.(2016).} 

\bibitem{[BC-EJP]} {\small Bally V., Caramellino L.: On the distances
between probability density functions. \emph{Electronic Journal of
Probability}, \textbf{19}, no. 110, 1-33 (2014).} 

\bibitem{[BCNon]} {\small Bally V., Caramellino L.: An Invariance
principle for Stochastic Series II. Non Gaussian limits. preprint arXiv
1607.04544 (2016)}.

\bibitem{[BCP]} {\small Bally V., Caramellino L., Poly G.: Convergence
in distribution norms in the CLT for non identical distributed random
variables. Preprint arXiv:1606.01629, (2016). }

\bibitem{[BR]} {\small Bally V., Ray C.:} {\small \emph{Approximation of
Markov semigroups in total variation distance. } Electronic J. of Probab.}
{\small 21, no 12.(2016).}

\bibitem{[BGL]} {\small Bakry D., Gentil I., Ledoux M.: \emph{%
Analysis and Geometry of Markov Diffusion Semigroups}. Springer (2014) }

\bibitem{[Be]} {\small Bentkus V.: On Hoeffding's inequalities. \emph{Ann. Probab.} \textbf{32}, 1650--1673 (2004)}
    
\bibitem{[Bo]} {\small Bogachev V.I., Kosov V.I., Zelenov G.I.: 
\emph{Fractional smoothness of distributions of polynomials and fractional
analog of the Hardy-Landau-Littelwod inequality.} \texttt{arXiv:1602.05207v2}}

\bibitem{[CW]}{\small  Carbery A.,  Wright J.: Distributional and $L^{q}$ norm
inequalities for polynomials over convex bodies in ${\mathbb{R}}^{n}$. \emph{%
Math. Research Lett.} \textbf{8}, 233--248, 2001.}

\bibitem{[D]} {\small Davydov Y.A., Martynova G.V.: \emph{Limit
behaviour of multiple stochastic integrals..} Stat. and Control of
Stochastic Processes, Nauka, Preila, Moscow, 1987, pp 55-57 (Russian) }

\bibitem{[dJ1]} {\small de Jong P.: A central limit theorem for generalized
quadratic forms. \emph{Probab. Th. Rel. Fields} \textbf{75}, 261-277 (1987).
}

\bibitem{[dJ2]} {\small de Jong P.: A central limit theorem for generalized
multilinear forms. \emph{Journal of Multivariate Analysis} \textbf{34},
275-289 (1990). }

\bibitem{[GR]} {\small Gamkrelidze N.G., Rotar' V.I.: On the rate of
convergence in the limit theorem for quadratic forms. \emph{Theory Probab.
Appl. } \textbf{22}, 394-397 (1977). }

\bibitem{[GT]} {\small G\"{o}tze F., Tikhomirov A.N.: Asymptotc
distributions of quadratc forms. \emph{Ann.of Probab. }\textbf{27},
1072-1098 (1999).}

\bibitem{[Ha]} {\small  Halmos P.R.: The theory of unbiaised estimation.
\emph{Ann. Math. Statist.} \textbf{17}, 34-43 (1946). }

\bibitem{[Hof1]} {\small Hoeffding W.: A class of statistics with
assymtotically normal distribution. \emph{Ann. Math. Statist.} \textbf{19},
293-325 (1948). }

\bibitem{[Hof2]} {\small Hoeffding W.: The strong law of large numbers for U
Statistics . \emph{Institute of Statistics, Mimeo-Series No 302, University
of North Carolina} (1961). }


\bibitem{[kb]} {\small Koroljuk V. S., Borovskich Yu. V.: \emph{Theory
of U-statistics}. Mathematics and its Applications, 273. Kluwer Academic
Publishers Group, Dordrecht, 1994.}


\bibitem{[Lee]} {\small Lee A.J.: \emph{U-Statistics. Theory and Practice.}
Statistics: textbooks and monographs, Vol 110, (1990).}

\bibitem{[LL]} {\small L\"{o}cherbach E., Loukianova D.: On Nummelin
splitting for continous time Harris reccurent Markov processes and
application to kernel estimation for multi-dimensional diffusion processes.
\emph{SPA. }\textbf{118}, 1301-1321 (2008).}


\bibitem{[MDO]} {\small Mossel E., O'Donnell R., Oleszkiewicz K.:
Noise stability of functions with low influences: Invariance and optimality.
\emph{Ann. Math.} \textbf{171}, pp. 295-341 (2010). }

\bibitem{[N]} {\small Nummelin E.: A Splitting Technique for Harris
Reccurent Markov Chains. \emph{Z. Wahrsch. verw. Gebiete}  \textbf{43} 309-318 (1978)}

\bibitem{[NN]} {\small Noreddine S., Nourdin I.: On the Gaussian
approximation of vector-valued multiple integrals. \emph{J. Multiv. Anal.}
\textbf{102}, no. 6, 1008-1017 (2011). }

\bibitem{[NP]} {\small Nourdin I., Peccati G.: \emph{Normal
Approximations Using Malliavin Calculus: from Stein's Method to Universality}%
. Cambridge Tracts in Mathematics, 192 (2012). }

\bibitem{[NP1]} {\small Nourdin I., Peccati G.: Stein's method
on Wiener chaos. \emph{Probab. Theory Related Fields} \textbf{145}, 75-118 (2009).}

\bibitem{[NPPS]} {\small Nourdin I., Peccati G., Poly G., Simone R.: 
Classical and free Forth Moment Theorems: universality and thresholds.  
\emph{J. Theoretical Probability} \textbf{29}, 653-680 (2016)}

\bibitem{[NPRein]} {\small Nourdin I., Peccati G., Reinert G.:
Invariance principles for homogeneous sums: universality of Wiener chaos.
\emph{Ann. Probab.} \textbf{38}, no. 5, 1947-1985 (2010). }

\bibitem{[NPRev]} {\small Nourdin I., Peccati G., R\'{e}veillac A.:
Multivariate normal approximation using Stein's method and Malliavin
calculus. \emph{Ann. Inst. H. Poincar\'{e} Probab. Statist.} \textbf{46},
no. 1, 45-58 (2010). }

\bibitem{[NPy]} {\small Nourdin I., Poly G.: Convergence in total
variation on Wiener chaos. \emph{Stochastic Process. Appl.} \textbf{123},
651--674 (2013). }

\bibitem{[NPy1]} {\small Nourdin I., Poly G.: An invariance principle
under the total variation distance. \emph{Stochastic Process. Appl.} vol
\textbf{125}, issue 6, p 2190-2205 (2015).}

\bibitem{bib:[N]} {\small Nualart D.: \emph{The Malliavin calculus and
related topics. Second Edition}. Springer-Verlag (2006). }

\bibitem{[NO]} {\small Nualart D., Ortiz-Latorre S.: Central limit
theorem for multiple stochastic integrals and Malliavin calculus. \emph{%
Stoch. Processes Appl.} \textbf{118}, 614-628 (2008). }

\bibitem{[PN]} {\small Nualart D., Peccati G.: Central limit theorems
for sequences of multiple stochastic integrals. \emph{Annals of Probability}
\textbf{33}, 177-193 (2005). }

\bibitem{[PT]} {\small Peccati G., Tudor C.A.: Gaussian limits for
vector-valued multiple stochastic integrals. \emph{S\'{e}minaire de
Probabilit\'{e}s XXXVIII}, 247-262 (2004). }

\bibitem{[Py]} {\small Poly, G.: Dirichlet forms and applications to the
ergodic theory of Markov chains. Phd thesis,
\texttt{htttps://tel.archives-ouvertes.fr/tel-00690724}, (2012) }

\bibitem{[PROH]} {\small Prohorov Y.: On a local limit theorem for
densities. Doklady Akad. Nauk SSSR (N.S.)83, 797-800 (1952). In Russian.}

\bibitem{[RS]} {\small Rotar' V.I., Shervvashidze T.L.: \emph{Some
estimates of distributions of quadratic forms. } Theory Probab. Appl. \textbf{%
6}, 738-751 (1985)}
\end{thebibliography}
\end{document}